\documentclass[12pt,a4paper,reqno]{amsart}
\usepackage{amssymb}
\usepackage{amscd}
\usepackage{enumerate}
\usepackage{graphicx}
\usepackage{siunitx}
\usepackage{tikz-cd}
\usepackage{color}
\usetikzlibrary{arrows}
\numberwithin{equation}{section}

\usepackage{mathabx}

\usepackage{mathtools}
\usepackage[tableposition=top]{caption}
\usepackage{booktabs,dcolumn}

\newtheorem{theorem}{Theorem}[section]
\newtheorem{lemma}[theorem]{Lemma}
\newtheorem{corollary}[theorem]{Corollary}
\newtheorem{proposition}[theorem]{Proposition}
\newtheorem{example}[theorem]{Example}

\theoremstyle{remark}
\newtheorem{remark}{Remark}[section]
\newtheorem{definition}[remark]{Definition}

\def\C{{\mbox{\rm\kern.24em
\vrule width.03em height1.43ex depth-.052ex \kern-.26em C}}}
\def\QSet{\mbox{\rm\kern.24em
\vrule width.03em height1.48ex depth-.051ex \kern-.26em Q}}
\def\E{{\mbox{\rm I\kern-.22em E}}}
\def\P{{\bf P}}
\def\D{{\bf D}}
\def\T{{\bf T}}
\def\S{{\mathbf{S}}}

\def\size{{\rm size}}
\def\diam{{\rm diam}}

\def\BMO{{\operatorname{BMO}}}
\def\M{{\operatorname{M}}}

\def\J{{\bf J}}

\def\X{{\mathbf X}}
\def\Z{{\mathbf Z}}
\def\R{{\mathbf R}}
\def\W{{\mathcal Z}}

\def\\P{{\mathcal P}}

\def\F{{\mathcal F}}

\def\I{{\mathcal I}}

\def\eps{\varepsilon}

\def\dist{{\rm dist}}

\def\rank{{\rm rank}}

\def\bas{\begin{align*}}
\def\eas{\end{align*}}
\def\bi{\begin{itemize}}
\def\ei{\end{itemize}}
\def \endprf{\hfill  {\vrule height6pt width6pt depth0pt}\medskip}
\def\emph#1{{\it #1}}

\begin{document}
\title{Maximal multilinear operators}

\author{Ciprian Demeter}
\address{Department of Mathematics, UCLA, Los Angeles CA 90095-1555}
\email{demeter@@math.ucla.edu}

\author{Terence Tao}
\address{Department of Mathematics, UCLA, Los Angeles CA 90095-1555}
\email{tao@@math.ucla.edu}

\author{Christoph Thiele}
\address{Department of Mathematics, UCLA, Los Angeles CA 90095-1555}
\email{thiele@@math.ucla.edu}

\thanks{The first author was supported by NSF Grant DMS-0556389}
\thanks{The third author was supported by NSF Grant DMS-0400879}
\begin{abstract}
We establish multilinear $L^p$ bounds for a class of maximal multilinear averages of functions on one variable, reproving and generalizing the bilinear maximal function bounds of Lacey \cite{La}.  As an application we obtain almost everywhere convergence results for these averages, and in some cases we also obtain almost everywhere convergence for their ergodic counterparts on a dynamical system.
\end{abstract}

\maketitle

\section{Introduction}

Let $n>1$, $m\ge 1$ and consider an $(n-1)\times {m}$ real-valued matrix $A={(a_{i,j})}_{i=1\:j=1}^{n-1\:m}$. This naturally gives rise to the multilinear averages:
\begin{equation}
\label{eq:avhar}
T_{A,\R,r}(f_1,\ldots,f_{n-1})(x) := \frac{1}{(2r)^{m}}\int_{|t_1|,\ldots,|t_m|\le r}\prod_{i=1}^{n-1}f_i(x+\sum_{j=1}^{m}a_{i,j}t_j)d\vec{t},  
\end{equation} 
where $r > 0$ and $f_1,\ldots,f_{n-1}$ are arbitrary measurable functions on $\R$.
Part of the motivation for considering such averages comes from ergodic theory. Let $\X=(X, \Sigma, m, S)$ be a dynamical system, i.e. a complete probability space $(X,\Sigma,m)$ endowed with an invertible bimeasurable transformation $S:X\to X$ such that $m{S^{-1}}=m$. We define the iterates $S^n: X \to X$ for $n \in \Z$ in the usual manner. In case the matrix $A$ has integer entries, one can consider the following ergodic averages:
 \begin{equation}
\label{eq:averg}
T_{A,\X,L}(f_1,\ldots,f_{n-1})(x) := \frac{1}{(2L+1)^{m}}\sum_{|l_1|,\ldots,|l_m|\le L}\prod_{i=1}^{n-1}f_i(S^{\sum_{j=1}^{m}a_{i,j}l_j}x).  
\end{equation}
We use $L^p(\R)$ to denote the usual Lebesgue spaces on $\R$, and $L^p(\X)$ to denote the Lebesgue spaces on the dynamical system $\X$.

In this paper we shall be primarily concerned with the problem of almost everywhere convergence of these averages as $r \to 0$ or $L \to \infty$
in the case that the $f_i$ obey some $L^{p_i}$ type integrability condition.  As it is well known, such problems are related to the boundedness properties of the maximal operator 
\begin{equation}
\label{eq:avmax}
\begin{split}
T_{A,\R}^*(f_1,\ldots,f_{n-1})(x) &:= \sup_{r > 0} |T_{A,\R,r}(f_1,\ldots,f_{n-1})(x)| \\
&=
\sup_{r > 0} |\frac{1}{(2r)^{m}}\int_{|t_1|,\ldots,|t_m|\le r}\prod_{i=1}^{n-1} f_i (x+\sum_{j=1}^{m}a_{i,j}t_j)d\vec{t}|  
\end{split}
\end{equation} 
or the closely related maximal operator
\begin{equation}\label{tax-ergodic}
\begin{split}
T^*_{A,\X}(f_1,\ldots,f_{n-1})(x) &:= 
\sup_{L > 0} |T_{A,\X,L}(f_1,\ldots,f_{n-1})(x)| \\
&=
\sup_{L > 0} |\frac{1}{(2L+1)^{m}}\sum_{|l_1|,\ldots,|l_m|\le L}\prod_{i=1}^{n-1}f_i(S^{\sum_{j=1}^{m}a_{i,j}l_j}x)|.  
\end{split}
\end{equation}
It turns out that standard transference arguments allow one to convert any positive or negative boundedness result for
$T^*_{A,\R}$ to one for $T^*_{A,\X}$ and vice versa; see Proposition \ref{correspondence}.  Thus we shall view the boundedness problems for these two
maximal operators as being equivalent.

Since one can easily establish convergence for \eqref{eq:avhar} in any reasonable topology 
when the $f_1,\ldots,f_{n-1}$ are smooth, compactly supported functions, a standard density argument then shows that as soon as
the maximal operator $T_{A,\R}^*$ maps $L^{p_1}(\R) \times \ldots \times L^{p_{n-1}}(\R)$ to weak $L^q(\R)$ for some $0 < q < \infty$, then the averages \eqref{eq:avhar} will converge pointwise almost everywhere when $f_i \in L^{p_i}(\R)$ for $1 \leq i \leq n-1$, at least in the case when all the $p_1,\ldots,p_{n-1}$ are finite\footnote{When one or more of the exponents is $\infty$ one can proceed by localization arguments, exploiting the fact that an $L^\infty$ function is locally in $L^p$ for any $p < \infty$.  This costs us an epsilon in the exponents but in most of our results the range of exponents will be open and so this will not make any difference.}.  In fact
these averages will converge almost everywhere to the pointwise product $f_1 \ldots f_{n-1}$.  In the converse direction, Stein's maximal principle
\cite{stein:maximal} shows that in many cases, almost everywhere convergence of \eqref{eq:avhar} can only be established via such weak $L^q$ bounds on
the maximal operator $T_{A,\R}^*$.

For the ergodic averages \eqref{eq:averg}, the situation is more difficult because there is no obvious counterpart of the
class $C^\infty_c(\R)$ of  smooth compactly supported functions on which the convergence is easy to establish\footnote{An alternate approach would be to establish either a $V^q$ variational estimate on $T_{A,\X,L}$ in $L$ for some $q < \infty$, or an oscillation inequality, since any of these automatically implies convergence as $L \to \infty$, in the spirit of Doob's inequality or Lepingle's inequality.  We will not pursue such an approach here, but see for instance \cite{De1}, \cite{DLTT}.}.  However, one can use
the class $L^\infty(\X)$ as a substitute, in the sense that once almost everywhere convergence for $T_{A,\X,L}$ is established for
$f_1,\ldots,f_{n-1} \in L^\infty(\X)$, one can extend this convergence result to the case when $f_i \in L^{p_i}(\X)$ provided that
one knows that the maximal operator $T_{A,\R}^*$ maps $L^{p_1}(\R) \times \ldots \times L^{p_{n-1}}(\R)$ to weak $L^q(\R)$ for some $0 < q < \infty$, since transference arguments then give an analogous boundedness statement for $T_{A,\X}^*$. Thus the problem of almost everywhere convergence of $T_{A,\X,L}$
for functions $f_i \in L^{p_i}(\X)$ factors into two rather distinct problems, namely
 establishing convergence for $L^\infty(\X)$ functions (which is a problem in ergodic theory), and establishing
 a bound for $T^*_A$ (which is a problem in multilinear harmonic analysis).  In this paper we shall focus almost exclusively on the latter problem. The former problem is quite difficult, except when $n=2$; the $n=3$ case already requires
 a deep result of Bourgain \cite{Bo2}, and convergence for higher $n$ is only proven for very special averages (see e.g. \cite{As:1}) or with  additional spectral assumptions on the shift $S$ (see \cite{As:2}, \cite{Le:1}, \cite{Le:2}), and we will not make progress on these issues here.

Let us now discuss some important special cases of the above general setup.

\subsection{Linear averages} If $n=2$, $m=1$, and $A = (a)$ for some non-zero integer $a$, then we have
$$ T_{A,\R,r} f_1(x) = \frac{1}{2r} \int_{-r}^r f_1(x+at)\ dt$$
and
$$ T_{A,\X,L} f_1(x) = \frac{1}{2L+1} \sum_{l=-L}^L f_1(S^{al} x).$$
If $f_1$ is in $L^p(\R)$ (resp. $L^p(\X)$) for some $1 \leq p \leq \infty$, then the Lebesgue differentiation theorem (resp. the Birkhoff ergodic
theorem) shows that $T_{A,\R,r} f_1$ (resp. $T_{A,\X,L}$) are almost everywhere convergent.  Both of these results require
the \emph{Hardy-Littlewood maximal inequality}, which asserts that the Hardy-Littlewood maximal operator
$$ \M f_1(x) := \sup_{r > 0} \frac{1}{r} \int_{-r}^r |f_1|(x+t)\ dt$$
maps $L^1$ to weak $L^1$.  The Lebesgue differentiation theorem follows immediately from the maximal inequality, whereas the Birkhoff ergodic theorem requires that one first establish almost everywhere convergence for a dense class such as $L^\infty(\X)$.

\subsection{Bilinear averages}\label{bilinear} 
Let $n=3$, $m=1$, and $\displaystyle A = \left( \begin{array}{l} a_1 \\ a_2 \end{array} \right)$
for some distinct non-zero integers $a_1,a_2$, thus
$$ T_{A,\R,r}(f_1,f_2)(x) = \frac{1}{2r} \int_{-r}^r f_1(x+ a_1 t) f_2(x + a_2 t)\ dt$$
and 
$$ T_{A,\X,L}(f_1,f_2)(x) = \frac{1}{2L+1} \sum_{l=-L}^L f_1(S^{a_1 l} x) f_2(S^{a_2 l} x).$$
As a consequence of a deep theorem of Bourgain \cite{Bo2} (relying on Fourier analysis on the torus), it is known that the averages
$T_{A,\X,L}(f_1,f_2)$ converge almost everywhere whenever $f_1, f_2 \in L^\infty(X)$.  Almost everywhere convergence in other classes then pivots on understanding the bilinear maximal operator
$$ T^*_{A,\R}(f_1,f_2)(x) = \sup_{r > 0} |\frac{1}{2r} \int_{-r}^r f_1(x+ a_1 t) f_2(x + a_2 t)\ dt|.$$
This operator clearly maps $L^\infty(\R) \times L^\infty(\R) \to L^\infty(\R)$, and from the Hardy-Littlewood maximal inequality
it also maps $L^\infty(\R) \times L^1(\R)$ or $L^1(\R) \times L^\infty(\R)$ to weak $L^1$.  This, combined with bilinear interpolation, is enough to establish
almost everywhere convergence of the ergodic averages $T_{A,\X,L}$ for $f_1 \in L^{p_1}(X), f_2 \in L^{p_2}(X)$ when $1/p_1 + 1/p_2
< 1$ (one also obtains the edge $1/p_1 +1/p_2 = 1$ from this argument as long as $p_1, p_2 < \infty$).  It was shown by Lacey \cite{La}, using time-frequency analysis, that  $T^*_{A,\R}$ in fact
maps $L^{p_1}(\R) \times L^{p_2}(\R)$ to $L^q(\R)$ whenever $\frac{1}{q}=\frac{1}{p_1}+\frac{1}{p_2}$ and $q>\frac{2}{3}$.  
This allows one to extend the almost everywhere convergence result to the larger range $1/p_1 + 1/p_2 < 3/2$.
It is an interesting question as to whether this is the true limit for these results.  Certainly one has boundedness
for a single-scale operator $T_{A,\R,r}$ or $T_{A,\X,L}$ all the way up to the range $1/p_1 + 1/p_2 \leq 2$.  On the other hand,
the time-frequency approach is known to break down at $1/p_1 + 1/p_2 = 3/2$ (see \cite{La}).

\subsection{Furstenberg averages}\label{furstenberg} 
Let $n \geq 2$, $m=1$, and let $A$ be the matrix
$$ A := \left( \begin{array}{l}
1 \\
2 \\
\vdots \\
n-1
\end{array}\right).$$
Then \eqref{eq:avhar} becomes the multilinear average
$$ T_{A,\R,r}(f_1,\ldots,f_{n-1})(x) = \frac{1}{2r} \int_{-r}^r \prod_{i=1}^{n-1} f_i( x + it )\ dt$$
and \eqref{eq:averg} becomes the Furstenberg average
$$ T_{A,\X,L}(f_1,\ldots,f_{n-1})(x) = \frac{1}{2L+1} \sum_{l=-L}^L f_i(S^{il} x).$$
Note the cases $n=2, n=3$ are special cases of the linear and bilinear averages considered earlier.  These averages are related to
the Furstenberg recurrence theorem \cite{Fu} and to Szemer\'edi's theorem on arithmetic progressions \cite{S}, and are also connected to the 
recent result in \cite{GT} that the primes contain arbitrarily long progressions.  For instance, the Furstenberg recurrence theorem is essentially the assertion that
$$ \liminf_{L \to \infty} \int_X T_{A,\X,L}(f,\ldots,f) f\ dm > 0$$
whenever $f$ is non-negative and does not vanish almost everywhere.  The question of norm convergence of $T_{A,\X,L}$ is more difficult and has only
been recently treated in the independent works of Host and Kra \cite{KH} and Ziegler \cite{Z}.  They showed that
 if $f_1,\ldots,f_{n-1} \in L^\infty(X)$ then $T_{A,\X,L}(f_1,\ldots,f_{n-1})$ converges in $L^2(X)$ norm (and hence in $L^p(X)$ norm for any $1 \leq p < \infty$).  Their approach relies on the reduction to convergence for functions in a sub-$\sigma$-algebra $\W_{n-1}$ of $\Sigma$, known as a \emph{characteristic factor}, on which $T$ can be represented as an inverse limit of translations on nilmanifolds. The advantage of such a concrete representation is that  this  particular type of translations is quite well understood. In particular, $\W_0$ is the  $\sigma$-algebra spanned by the invariant sets of powers of $T$, while the action of $T$ on the Kronecker factor $\W_1$ is isomorphic with a rotation on some abelian group. The $\sigma$-algebras $\W_{k}$ with $k\ge 2$ give rise to noncommutative factors which require a more delicate analysis.  The work in this paper will however proceed in a different direction, focusing on the quantitative bounds of various operators associated to these averages rather than analyzing characteristic factors.  It is of course possible to extend these norm convergence results to functions $f_i$ in other spaces $L^{p_i}(X)$ by exploiting boundedness properties $T_{A,\X,L}$ or $T_{A,\R,r}$, but we will not pursue this issue here, though we
will mention that some surprising subtleties in this problem in the case $1/p_1 + \ldots + 1/p_n > 1$ have been 
 uncovered by Christ \cite{CH}.
 
The problem of almost everywhere convergence, as opposed to norm convergence, for the Furstenberg averages remains open
even for $n=4$.  One can obtain some bounds of the corresponding maximal operators in $L^p$ spaces by leveraging the corresponding
bounds in the bilinear setting.  For instance one can extend Lacey's bilinear estimates mentioned earlier to the multilinear setting by estimating all but two of the functions in $L^\infty$.  This ultimately leads to a bound on $T^*_{A,\R}$ from $L^{p_1}(\R) \times \ldots \times L^{p_n}(\R)$ to $L^q(\R)$ whenever $1 < p_1,\ldots,p_n \leq \infty$ and $1/q = 1/p_1 + \ldots + 1/p_n < 3/2$.  

\subsection{Averages along cubes}\label{cubes}

The work of Host and Kra \cite{KH} related the norm convergence of the above Furstenberg averages to the norm convergence of
averages of cubes, which is a special case of \eqref{eq:averg} with $n=2^m$. To define them, let $V_m$ be the index
set $V_m:=\{0,1\}^m\setminus \{0\}^m$. The averages on the $m$-dimensional cubes are 
 \begin{equation}
\label{eq:cubes}
\frac{1}{(2L+1)^{m}}\sum_{\vec{i}\in \{-L, \ldots, L\}^m}\prod_{\epsilon\in V_m}f_{\epsilon}(S^{\vec{i}\cdot\epsilon}x).
\end{equation} 
For example, when $m=1$ (so $n=2$) we just have a linear averaging operator.  
When $m=2$ (and so $n=4$), this averaging operator along squares is essentially the same as $T_{A,\X,L}$ with
$$ A := \left( \begin{array}{ll}
0 & 1 \\
1 & 0 \\
1 & 1
\end{array} \right)$$
while when $m=3$ (and $n=8$) the averaging operator along cubes is essentially $T_{A,\X,L}$ with
$$ A := \left( \begin{array}{lll}
0 & 0 & 1 \\
0 & 1 & 0 \\
0 & 1 & 1 \\
1 & 0 & 0 \\
1 & 0 & 1 \\
1 & 1 & 0 \\
1 & 1 & 1
\end{array} \right).$$
It is proved in \cite{KH} that the averages\footnote{Actually, a more general class of averages is shown in \cite{KH}  to have $\W_{m-1}$ as a characteristic factor; we refer the reader to \cite{KH} for the details.} in ~\eqref{eq:cubes} have $\W_{m-1}$ as a characteristic factor for $L^2$-norm convergence, and as a consequence that these averages converge in $L^2(\X)$ whenever
$f_\epsilon \in L^\infty(\X)$.  Using these characteristic factors, Assani \cite{As:1} showed that these averages also converged
pointwise almost everywhere when $f_\epsilon \in L^\infty(\X)$. It is somehow peculiar that these techniques do not seem to be able to give an alternative (non-Fourier analytical) proof to Bourgain's pointwise result mentioned earlier. 

To extend the latter $L^\infty(\X)$ convergence result to an $L^p(\X)$ convergence result requires control of a maximal function.
For sake of concreteness let us just focus on the case $m=2$, where the relevant maximal function is
$$ \sup_{r > 0} |\frac{1}{(2r)^2} \int_{-r}^r \int_{-r}^r f_{10}(x + t_1) f_{01}(x+t_2) f_{11}(x+t_1+t_2)\ dt_1 dt_2|.$$
One can deduce a certain number of bounds on this maximal function from the Hardy-Littlewood maximal inequality and multilinear
interpolation.  Indeed, the maximal inequality and H\"older's inequality implies that this maximal function lies in weak $L^{1/2}$ whenever two of $f_{10}, f_{01}, f_{11}$ lie in $L^1$ and the other one lies in $L^\infty$, while this maximal operator is trivially in $L^\infty$ when
all three of $f_{10}, f_{01}, f_{11}$ lie in $L^\infty$.  Interpolation then gives bounds (and hence almost everywhere
convergence of the associated  averages along squares) when $f_{01} \in L^{p_{01}}$, $f_{10} \in L^{p_{10}}$, $f_{11} \in L^{p_{11}}$ with
$1/p_{01} + 1/p_{10} + 1/p_{11} < 2$, with an extension to the boundary $1/p_{01} + 1/p_{10} + 1/p_{11} = 2$ when all of
the exponents are finite.  As a corollary of our main result (which is proven using time-frequency techniques)
 we shall be able to extend this range to $1/p_{01} + 1/p_{10} + 1/p_{11} < 5/2$, in analogy with the situation for bilinear
 averages discussed earlier (see Corollaries \ref{geas0}, \ref{geas} below).

\subsection{Main results}

We now study the maximal operator $T^*_{A,\R}$ defined in \eqref{eq:avmax} for a general $(n-1)\times m$ matrix $A=(a_{i,j})$; we will allow the $a_{i,j}$ here to be non-integer as one can still define $T^*_{A,\R}$ in this case.  To state the main result we need some notation.  We introduce the extended matrix $\E(A)$, which is the $n \times (m+1)$ matrix
$$\E(A):=\begin{pmatrix}
a_{1,1}&a_{1,2}&\dots & a_{1,m}&1\\
a_{2,1}&a_{2,2}&\dots & a_{2,m}&1\\
\hdotsfor[2.0]4\\
a_{n-1,1}&a_{n-1,2}&\dots & a_{n-1,m}&1\\
0&0&\dots &0&1
\end{pmatrix}.
$$
Note that the range of this matrix consists of all $n$-tuples of the form
$$ (x+\sum_{j=1}^{m}a_{1,j}t_j, \ldots, x+\sum_{j=1}^{m}a_{n-1,j}t_j, x)$$
for $x, t_1, \ldots, t_m \in \R$.

A set of row indices $i$ is said to be a \emph{set of
linear independence} for a matrix $B$, if the set of corresponding 
rows of $B$ is linearly independent.
Given a matrix $A$, let $S_{A,\epsilon}$ for $0<\epsilon<1/4$ be the set
of all tuples 
$(x_1,\dots,x_{n-1})$ where $x_i\in \{0,1/2+\epsilon,1-\epsilon\}$ for all $i$, there 
is at most one index $i$ with $x_i=1/2+\epsilon$, the indices $i$ with $x_i=1-\epsilon$ 
form a set of linear independence for $A$, and the indices $i$ 
with $x_i\in \{1/2+\epsilon,1-\epsilon\}$ form a set of linear independence for $\E(A)$.
Let $H_{A,\epsilon}$ be the convex hull of $S_{A,\epsilon}$ and let $H_A$ be 
the union of all $H_{A,\epsilon}$ with $0<\epsilon<1/4$.

The following is our main theorem:
\begin{theorem}
\label{thm:ct1}
Assume $n\ge 3$ and let $A$ be a  matrix as above. 
Let $(p_1,\ldots,p_{n-1})$ be a tuple of real numbers with 
\begin{equation}
\label{z1}
1<p_i\le \infty
\end{equation}
for $1\le i\le n-1$ and set
\begin{equation}
\label{z2}
\frac 1{{p_n}'}=\sum_{i=1}^{n-1}\frac{1}{p_i}.
\end{equation}
If 
\begin{equation*}
(1/p_1,\ldots 1/p_{n-1})\in H_A,
\end{equation*}
then the operator $T_{A,\R}^{*}$ 
$$T_{A,\R}^{*}:L^{p_1}\times \dots\times L^{p_{n-1}}\to L^{p_n'}$$
is bounded.
\end{theorem}

\begin{remark}  The condition \eqref{z2} is mandated by scaling considerations (i.e. dimensional analysis).
As we shall see shortly, the theorem is trivial if one restricts the tuples $(1/p_i)$ to the convex   
hull of those points in $S_{A,\epsilon}$ which do not have a component equal to $1/2+\epsilon$.
This happens in particular when $n=2$. Thus, in a nutshell, we are gaining $1/2-\epsilon$ over the trivial 
estimates. 
\end{remark}

	\begin{remark}
For some matrices $A$ we can obtain a
better range of exponents than stated in the theorem. Namely, when the matrix $A$ is a diagonal block matrix, 
we may gain $1/2$ for every block. More precisely, the argument works for $A$ 
upper block triangular and $\E(A)$ modulo the last column and restricted to the rows other than
the last row is block diagonal. The argument involves only separation of variables and H\"older's 
inequality, so we shall not elaborate on this.
\end{remark}

The following corollary is weaker than the theorem, but has the advantage of an easy description
of the range of exponents and covers many of the cases of interest. Define the \emph{nondegeneracy rank} of the matrix $A$, denoted by $\rank^{*}(A)$, 
to be the largest integer $r$ such that any $r$ rows of $A$ are linearly independent. It is an immediate observation that $\rank^*(A)+1\ge \rank^*(\E(A))\ge \rank^*(A)$.

\begin{corollary}
\label{cor:ct2}
Assume $n\ge 2$ and let $A$ be a  matrix as above.
Define the complexity parameter $k=n-\rank^*(\E(A))$.
Let $(p_1,\ldots,p_{n-1})$ be a tuple with 
$$
1<p_i\le \infty
$$
for $1\le i\le n-1$ and set
$$
\frac 1{{p_n}'}=\sum_{i=1}^{n-1}\frac{1}{p_i}.
$$
If 
\begin{equation}
\label{z4}
\frac{1}{p_{1}}+\dots+\frac{1}{p_{n-1}}<n-k-\frac{1}{2},
\end{equation}
then the operator $T_{A,\R}^{*}$ 
$$T_{A,\R}^{*}:L^{p_1}\times \dots\times L^{p_{n-1}}\to L^{p_n'}$$
is bounded.
\end{corollary}
\begin{proof}
The closure of the region of tuples $(1/p_i)$ in the corollary is the intersection
of the cube $[0,1]^{n-1}$ with a half space. All extremal points of this set are on an
edge of the cube and thus have all but at most one coordinate in $\{0,1\}$. The only 
possible value for the exceptional coordinate is $1/2$ as the right- 
hand-side of (\ref{z4}) is equal to $1/2$ modulo the integers. 
Thus the region in the corollary is the convex hull of all tuples $(x_1,\dots,x_2)$
with at most $n-k-1$ components equal to $1-\epsilon$, at most one component equal to
$1/2+\epsilon$ and the remaining components equal to $0$. The corollary then follows from the rank conditions on $A$ and $\E(A)$ and the fact that $\rank^*(A)\ge n-k-1$.
\end{proof}

\begin{remark} As discussed earlier, the boundedness results in Theorem \ref{thm:ct1} and Corollary \ref{cor:ct2}
immediately imply almost everywhere convergence for $T_{A,\R,\eps}(f_1,\ldots,f_{n-1})$ as $\eps \to 0$ 
when $f_i \in L^{p_i}(\R)$ if all the $p_i$ are finite, since this convergence is trivial for $f_i$ in the dense 
class $C^\infty_c(\R)$.  The $p_i=\infty$ cases can also be handled by a localization argument and exploiting some open-ness properties of $H_A$.  The situation for the ergodic averages is however substantially more difficult.
\end{remark}  

\begin{remark} If $\rank^*(\E(A))=\rank(\E(A))$,
then the regions described in Theorem \ref{thm:ct1} and Corollary \ref{cor:ct2} are equal.
\end{remark}

\begin{remark}
It is worth noting that $p_n'$ can be less than $1$, indeed it is less than $1$ in all nontrivial cases.
In some cases one can get below 1 by using just H\"older's inequality and interpolation, see for instance the discussion
in Section \ref{cubes}. 
\end{remark}

\begin{remark} Theorem ~\ref{thm:ct1} is a direct analog of the singular integral version in \cite[Theorem 1.1]{MTT1}, which roughly speaking replaces $T_{A,\R}^*$ with the related expression
$$
p.v. \int_{\R^m}\prod_{i=1}^{n-1} f_i (x+\sum_{j=1}^{m}a_{i,j}t_j) K(\vec{t}) d\vec{t}$$
for some Calder\'on-Zygmund kernel $K$.  As a consequence the methods of proof are quite similar. The parameter $k$ in 
Corollary \ref{cor:ct2} plays the same role as the parameter $k$ appearing in \cite[Theorem 1.1]{MTT1}, measuring the 
complexity of the averages under investigation. The case $k=0$ for the singular integral version can be solved with classical methods, namely Littlewood-Paley theory or wavelets, just as the case $k=0$ for the maximal version can be solved using the classical Hardy-Littlewood maximal inequality. 

Readers familiar with \cite{MTT1} will observe that the range of exponents in Theorem \ref{thm:ct1} is somewhat more permissive
than that in \cite{MTT1}.  More precisely, the restriction $k<\frac{n}2$ as well as several restrictions on the exponents $p_i$
from \cite{MTT1} are not needed in Theorem ~\ref{thm:ct1}. 
This is a consequence of the fact that there are trivial reductions in the maximal operator case if there are
exponents $p_i=\infty$, while in the singular integral setting there are no such trivial reductions.  This explains why for instance
we can obtain nontrivial estimates for the trilinear maximal operator ($n=4,k=2$)
\begin{equation}\label{trilinear}
T_{\vec{a}}^{*}(f_1,f_2,f_3) := \sup_{\epsilon>0}\frac{1}{\epsilon}\int_{|t|\le \epsilon}|f_1(x+a_1t)f_2(x+a_2t)f_3(x+a_3t)|\ dt
\end{equation}
with $a_1$, $a_2$, $a_3$, $0$ pairwise different (see Example \ref{4ex} below), despite the fact that no $L^p$ bounds of any
sort are known for the trilinear Hilbert transform
$$ p.v. \int_\R f_1(x+a_1t)f_2(x+a_2t)f_3(x+a_3t)\ \frac{dt}{t}.$$
\end{remark}

\begin{remark}
It should be emphasized that the nontrivial estimates from the $k>1$ cases are all obtained by such trivial reductions 
to the case $k=1$ and multilinear interpolation. In other words, there is no special theory developed yet to address the case $k>2$. It is quite probable that more sophisticated techniques will extend the range of the exponents in this case.
An interesting connection concerns the fact that  averages corresponding to some $k\ge 0$ appear to have $\W_{k}$ as a characteristic factor for $L^2$-norm convergence. In particular, it is an exercise based on the techniques from  \cite{As:1} and  from \cite{Bo2} to show that  $\W_{k}$ is  the characteristic factor  even for a.e. convergence, when $k=0,1$. This would support the evidence that, as in the case of norm convergence,  $k$ is the only parameter which dictates the complexity of the averages and of the techniques needed for the proof. 

Similar difficulties are encountered when dealing with polynomial maximal operators such as $P^{*}(f_1,f_2)(x) :=\sup_{\epsilon>0}\frac{1}{\epsilon}\int_{|t|\le \epsilon}|f_1(x+t)f_2(x+t^2)|\ dt$. In all these instances, the decomposition of the maximal operator, as explained in third section below, gives rise to a summation over a larger family of multidimensional cubes, each of which is indexed by more than just one parameter. Curiously, the boundedness of the maximal operator associated with polynomial averages, unlike the multilinear averages studied here (see Proposition \ref{correspondence}), does not in general transfer from harmonic analysis to ergodic theory. It is really that the results in  these two contexts have different meaning and most probably distinct ideas behind their proofs. An illuminating contrast comes from the fact that   $\sup_{\epsilon>0}\frac{1}{\epsilon}\int_{|t|\le \epsilon}|f(x+t^2)|d t$ can be easily  bounded by the Hardy-Littlewood maximal function, while Bourgain showed that the convergence of the ergodic averages along squares needs completely new ideas \cite{Bo1}.  
\end{remark}

Let us illustrate Theorem \ref{thm:ct1} and Corollary \ref{cor:ct2} with some examples.

\begin{example} Consider the bilinear averages from Section \ref{bilinear}.  Here the extended matrix is
$$\E(A)=\begin{pmatrix}
a_1 & 1 \\
a_2 & 1 \\
0   & 1
\end{pmatrix}.
$$
One can check that $\rank(A) = 1$ and $\rank^*(\E(A)) = \rank(\E(A)) = 2$, and 
$$S_{A,\epsilon} = \{ (0,0), (0,1/2+\epsilon), (1/2+\epsilon,0), (0, 1-\epsilon), (1-\epsilon, 0), (1/2+\epsilon,1-\epsilon), (1-\epsilon,1/2+\epsilon) \}$$
and hence
$$H_A = \{ (a,b): 0 \leq a,b < 1; a+b < 3/2 \}.$$
In this case Theorem \ref{thm:ct1} and Corollary \ref{cor:ct2} give the same results, namely recovering the bilinear maximal function
estimates of Lacey \cite{La} described earlier.  Indeed we give a reasonably self-contained\footnote{We will require some results from other papers, notably the multilinear interpolation theory from \cite{MTT1}, the weak Bessel inequality for forests (see e.g. \cite{MTT2}), a maximal Fourier inequality of Bourgain \cite{Bo1}, and an interval selection lemma of Lacey \cite{La}.} proof of the main results from \cite{La} here, following Lacey's approach.
\end{example}

\begin{example} Consider the $n=4$ Furstenberg average from Section \ref{furstenberg}.  Here the extended matrix is
$$\E(A):=\begin{pmatrix}
1 & 1 \\
2 & 1 \\
3 & 1 \\
0   & 1
\end{pmatrix}.
$$
One can check that $\rank(A) = 1$ and $\rank^*(\E(A)) = \rank(\E(A)) = 2$, and $S_{A,\epsilon}$ consists of those triples
$(a,b,c)$ with $a,b,c \in \{0,1/2+\epsilon,1-\epsilon\}$, at most one of $a,b,c$ equal to $1/2+\epsilon$, and at most one of $a,b,c$ equal to $1-\epsilon$. This gives
$$H_A = \{ (a,b,c): 0 \leq a,b,c < 1; a+b < 3/2 \}.$$
In this case, Theorem \ref{thm:ct1} and Corollary \ref{cor:ct2} recover the multilinear estimates mentioned at the end of
Section \ref{furstenberg} that can be trivially obtained from Lacey's bilinear result.  Similar considerations apply to higher
values of $n$.
\end{example}

\begin{example}\label{4ex} Consider the $m=2$ average along squares from Section \ref{cubes}.  Here the extended matrix is
$$\E(A):=\begin{pmatrix}
0 & 1 & 1 \\
1 & 0 & 1 \\
1 & 1 & 1 \\
0 & 0 & 1
\end{pmatrix}.
$$
One can check that $\rank(A) = 2$ and $\rank^*(\E(A)) = \rank(\E(A)) = 3$, and $S_{A,\epsilon}$ consists of those triples
$(a,b,c)$ with $a,b,c \in \{0,1/2+\epsilon,1-\epsilon\}$, at most one of $a,b,c$ equal to $1/2+\epsilon$, and at most two of $a,b,c$ equal or $1-\epsilon$. This gives
$$H_A = \{ (a,b,c): 0 \leq a,b,c < 1; a+b < 5/2 \}.$$
\end{example}

Combining the above example with Proposition \ref{correspondence} from Appendix \ref{correspond-sec}, and the result of Assani \cite{As:1}, we obtain the following corollary.

\begin{corollary}
\label{geas0}
Let $1<p_1,p_2,p_3\le \infty$ be such that $\frac1{p_1}+\frac1{p_2}+\frac1{p_3}<\frac52$. 
For every dynamical system ${\bf X}=(X,\Sigma,m,S)$, the averages on squares
$$\frac{1}{N^2}\sum_{i=-N}^{N}\sum_{j=-N}^{N}f_1(S^ix)f_2(S^jx)f_3(S^{i+j}x)$$ 
converge a.e. $x$, for each $f_i\in L^{p_i}({\bf X})$. 
\end{corollary}

\begin{remark}
A version of Corollary \ref{geas0} holds for  all averages with $k=1$. The convergence for $L^{\infty}$ functions follows by  using the aforementioned fact that these averages have characteristic factor $\W_1$ for pointwise convergence.  We omit the details.
\end{remark}

\begin{remark}
In \cite{DTT0} we use combinatorial methods involving sum set estimates to get nontrivial positive results in Corollary \ref{geas0}. This completely different approach  gives the result only in a small range, $p'_4>\frac{1}{2+\epsilon}$ for some unspecified $\epsilon$, and does not seem to extend to the case when  $p'_4$ is smaller then or even close to $\frac25$. 
\end{remark}

\begin{remark}
An interesting contrast to the results of Theorem ~\ref{thm:ct1} is provided by the constructions from \cite{De}, showing that some  maximal operators fail to be bounded when the indices $p_i,1\le i\le n-1$ are sufficiently close to 1. As a consequence, both Furstenberg's averages with $n\ge 4$ and the averages on cubes with $m\ge 3$ are proved to diverge a.e. in some range of $L^p$ spaces. The trilinear maximal operator from ~\eqref{trilinear}
has been proved in \cite{CH} to be unbounded for $p_1=p_2=p_3=p$, $1\le p<\frac32$, for appropriate choices of $\vec{a}$ depending on $p$.
The main ingredient behind these negative results is the fact that the polynomials $x+\sum_{j=1}^{m}a_{i,j}t_j,1\le i\le n-1$ are linearly dependent in $\R[x,t_1,\ldots,t_l]$ and hence $\rank^*(\E(A))\le n-2$ and $k\ge 2$.
In other words, our tools provide negative results only when $k\ge 2$, and all positive results are trivially deduced
from positive results when $k=0,1$. Further progress would require to break this barrier in the complexity $k$
either for positive or for negative results.
\end{remark}

The following is the straight-forward application of Corollary \ref{cor:ct2} to averages on cubes.
In this case, while $\rank(A)$ is the dimension of the cube, we have $\rank^*(A)=2$, an obstruction for higher 
nondegeneracy rank being the linear dependence of the polynomials $t_1$, $t_2$, and $t_1+t_2$. 
On the other hand, $\rank^*(\E(A)) = 3$, and hence:
\begin{corollary}
\label{geas}
Let $1<p_{\epsilon}\le \infty, \epsilon\in V_m$, be such that 
\begin{equation}
\label{masfg}
\sum_{\epsilon\in V_m}\frac1{p_{\epsilon}}<\frac52.
\end{equation}
For every dynamical system ${\bf X}=(X,\Sigma,m,S)$, the averages  on $m$-dimensional cubes ~\eqref{eq:cubes} converge a.e. for each $f_{\epsilon}\in L^{p_{\epsilon}}({\bf X})$.
\end{corollary}

This of course generalizes Corollary \ref{geas0}.
It would be interesting to know whether one can improve over $5/2$ on the right-hand-side of (\ref{masfg}).
Certainly the methods of this paper do not yield such an improvement, and \cite{De} provides an upper bound of
$28/5$ for the right-hand-side of (\ref{masfg}) for three dimensional cubes.

Theorem \ref{thm:ct1} is proven using standard time-frequency strategies, and in particular follows the approach of
Lacey \cite{La}, though it is more self-contained and employs some technical simplifications over that in \cite{La}.
In Section \ref{interp-sec} we use the theory of multilinear interpolation to
reduce Theorem \ref{thm:ct1} to a model case, Theorem \ref{reduct-00}, in which the matrix $A$ is in a simplified normal form,
the functions $f_1,\ldots,f_{n-1}$ have become $L^2$-normalized functions adapted to certain sets $E_1,\ldots,E_{n-1}$, and
the output is being measured in another set $E'_n$ which excludes a certain exceptional set determined by the Hardy-Littlewood maximal
function.  In Sections \ref{Fourier-sec}, \ref{sec:3} we use the Fourier transform and wave packet decomposition to reduce matters to bounding
a certain model sum (Theorem \ref{reduct-5}) involving the inner product of the functions $f_1,\ldots,f_n$ with various wave packets (and maximal wave packets) associated to a certain ``rank one'' collection of multitiles.  To estimate this model sum, we organize the collection of multitiles into trees; after obtaining an upper bound for the contribution of a single tree (see Proposition \ref{thm:2} and Section \ref{sec:6}) one quickly
reduces (essentially by summing a geometric series; see Section \ref{overview-sec}) to that of proving estimates for a tree selection algorithm (Lemma \ref{l:18}), which in turn reduces to
a certain maximal Bessel inequality concerning wave packets in a forest (Theorem \ref{thm:5a}, slightly improving and simplifying a similar result from \cite{La}).  This Bessel inequality will involve a certain logarithmic-type loss involving the size parameter $2^m$, but by some ``good-$\lambda$'' type reductions in Section \ref{bmo-sec} we can replace this factor with another logarithmic factor involving instead the multiplicity $\|N_\F\|_{L^\infty}$ of the forest (Theorem \ref{thm:3}).  After some
sparsification of the tile set, some elimination of exceptional tiles, and duality, one reduces to establishing a certain maximal Bessel inequality
on two families of tiles (see \eqref{ptop} and \eqref{pnotop}).  These inequalities are proven by using the time localization properties of 
wave packets, a non-maximal Bessel inequality (proven in Section \ref{nonmax-sec}), and the Radamacher-Menshov inequality.  In the case of
one of these inequalities \eqref{pnotop}, one also needs a maximal inequality of Bourgain \cite{Bo1}.  Finally, in an Appendix (Section \ref{correspond-sec}) we present a standard correspondence principle equating boundedness of maximal functions on $\R$ with maximal functions on measure-preserving systems.

\section{Interpolation reductions}\label{interp-sec}

The rest of the paper is devoted to the proof of Theorem \ref{thm:ct1}.  We shall use the methods of multilinear time-frequency
analysis and work entirely on $\R$, thus we will not make any further reference to the dynamical system $\X$.

In this section we use some multilinear interpolation techniques to reduce the operator $T^*_{A,\R}$ and the exponents
$p_1,\ldots,p_n$ to a standard form, and
then to also reduce the input functions $f_1,\ldots,f_{n-1}$ (and an additional output function $f_n$ arising from duality) to another standard form.

We first introduce some basic notation.
If $E$ is a measurable subset of $\R$, we use $1_E$ to denote the indicator function of $E$ and $|E|$ to denote the Lebesgue measure.
Also $\M f(x) := \sup_{r>0}\frac{1}{2r}\int_{x-r}^{x+r} |f|(y)dy$ denotes the classical Hardy-Littlewood maximal function.  
The notation $a\lesssim b$ or $a = O(b)$ means that $a\le cb$ for some universal constant $C$ (which will be allowed to depend on parameters such as $n$ and $p_1,\ldots,p_n$), and $a\sim b$ means that  $a\lesssim b$ and  $b\lesssim a$.  In some cases we will subscript the $\lesssim$ notation by a parameter to emphasize the fact that the constant $C$ involved can depend on that parameter, thus for instance $a \lesssim_\mu b$ means that $C$ can depend on $\mu$.  If $x \in \R^n$ we use $\|x\|$ to denote the Euclidean norm of $x$.

Now we can reduce the operator $T^*_{A,\R}$ and the exponents $p_1,\ldots,p_n$ to a standard form.

\begin{theorem}[First reduction]\label{reduct} Let $n \geq 3$, let $\Sigma$ be a hyperplane in $\R^{n-1}$ containing the origin but not containing any of the $n-1$ coordinate vectors $e_1,\ldots,e_{n-1}$ or the vector $(1,\ldots,1)$.
Then the $(n-1)$-linear operator $T^*$ defined by
\begin{equation}\label{tstar}
\begin{split}
&T^{*}(f_1,\ldots,f_{n-1})(x)= \\
&\sup_{r>0} \frac{1}{r^{n-2}}\int_{\vec t \in \Sigma: \|\vec t\| \le r} |f_1(x+t_1)\cdot\ldots\cdot f_{n-1}(x+t_{n-1})|d\vec{t}
\end{split}
\end{equation}
is bounded from $L^{p_1}(\R) \times \ldots \times L^{p_{n-1}}(\R)$ to $L^{p'_n}(\R)$
whenever $1 < p_1,\ldots,p_{n-1} < 2$, $\frac{5}{2}-n < \frac{1}{p_n} < 3-n$, and
$$ \frac{1}{p'_n} = \frac{1}{p_1} + \ldots + \frac{1}{p_{n-1}}.$$
The bound of course depends on $p_1,\ldots,p_n$ and the $\lambda_i$.
\end{theorem}

\begin{remark}
Note that $\rank(\E(A))=\rank^*(\E(A))=n-1$.
Hence we are in the case $k=1$ of Corollary \ref{cor:ct2} and the corollary is 
equivalent to Theorem \ref{thm:ct1} in this case.  The condition that $\Sigma$ does not contain $e_1,\ldots,e_{n-1}$ or $(1,\ldots,1)$ corresponds to the nondegeneracy condition in \cite{MTT1}.
\end{remark}

\begin{proof}[of Theorem \ref{thm:ct1} assuming Theorem \ref{reduct}] 
By multilinear interpolation as in \cite{MTT1}
it suffices to prove the estimate for tuples $(1/p_i)$ in $S_{A,\epsilon}$ for some $0 < \epsilon < 1/2$, so in particular $1/p_i = \{ 1/2+\epsilon, 1-\epsilon, 0\}$ for all $i$.  We may of course assume the $f_i$ are non-negative.
For each index $i$ with $p_i=\infty$ we can trivially estimate $f_i$ by its supremum norm
and remove it from the maximal operator:
$$\sup_{r>0}\frac{1}{(2r)^{m}}\int_{|t_1|,\ldots,|t_m|\le r}\prod_{i=1}^{n-1}
f_i(x+\sum_{j=1}^{m}a_{i,j}t_j)\ d\vec{t}$$
$$\lesssim \|f_{j}\|_{L^\infty} 
\sup_{r>0}\frac{1}{(2r)^{m}}\int_{|t_1|,\ldots,|t_m|\le r}\prod_{i\neq j}
f_i(x+\sum_{j=1}^{m}a_{i,j}t_j)\ d\vec{t}.$$
Doing this to each such exponent, we may assume without loss of generality that $1/p_i\in \{1/2+\epsilon, 1-\epsilon\}$ for all $i$.

If $1/p_i=1-\epsilon$ for all $i$, then by definition of $S_{A,\epsilon}$ 
the rows of the matrix $A$ are linearly independent and we may
do a change of variables so that $a_{i,j}$ is the Kronecker delta for $1\le i,j\le n-1$. Of course
the cube of integration in the parameter space $\{(t_1,\ldots,t_n)\}$ will be a parallelepiped in the new variables, 
but we may use the positivity of the $f_i$ and estimate the characteristic function of the parallelepiped by that of a cube, conceding a bounded loss in the estimates.
We may also assume that $A$ is a square matrix of dimension $m=n-1$, since in the case $m > n-1$ we may fix the 
variables $t_j$ with $j>n-1$
and apply the result in the square matrix case to fixed translates of the function $f_i$ obtaining an
$L^p(\R)$ bound independently of the translation. Then we perform a dummy average in
the variable $t_j$ with $j>n-1$ to obtain the desired estimate. In the square matrix case 
we estimate 
$$\sup_{r>0}\frac{1}{(2r)^{m}}\int_{|t_1|,\ldots,|t_{n-1}|\le r}\prod_{i=1}^{n-1}
|f_i(x+t_i)|d\vec{t}$$
$$\le \prod_{i=1}^{n-1}\sup_{\epsilon>0}\frac{1}{\epsilon}\int_{|t|\le \epsilon}|f_i(x+t)|\, dt$$
and then apply the Hardy Littlewood maximal theorem for $L^{1+\epsilon}$ and H\"older's 
inequality to obtain the desired estimate.

It remains to consider the case when $1/p_j=1/2+\epsilon$ for one index $j$ and $1/p_i=1-\epsilon$ for $i\neq j$; note
that this places these exponents in the situation of Theorem \ref{reduct}.
We may assume that $n \geq 3$ since the $n=2$ case follows from the Hardy-Littlewood maximal inequality.
By symmetry we may assume that $j=n-1$. The first $n-2$ rows of $A$ are linearly independent
and we may assume that $(a_{i,j})_{1\le i,j,\le n-2}$ is the Kronecker delta.
We may assume that the last row of $A$ is a linear combination of the other rows, or otherwise we 
can apply the reasoning of the previous paragraph. 
By a reasoning as in the previous paragraph we may also assume that $m\le n-2$.
Thus after a change of variables if necessary (and covering the resulting parallelepiped by a ball)
the operator $T_{A,\R}^*$ takes the form
\eqref{tstar} for some hyperplane $\Sigma$.
If $\Sigma$ contains $e_i$, then we perform the $t_i$ average first, estimate
the average using the Hardy Littlewood maximal function of $f_i$, and use H\"older's inequality 
to reduce matters to the case with one function less. We may thus assume that $\Sigma$ does not contain
any of the $e_i$.
 Finally, the hypothesis that the first $n-1$ rows of $\E(A)$ are linearly independent implies that $\Sigma$ does not
 contain $(1,\ldots,1$),
and the claim now follows from Theorem \ref{reduct}.
\end{proof}

\begin{remark} If the hypothesis $\frac{5}{2}-n < \frac{1}{p_n} < 3-n$ is replaced by $1 < p_n < \infty$ then Theorem \ref{reduct}
is easy to prove.  Indeed, in this case we can use H\"older's inequality to obtain the pointwise estimate
$$ T^*(f_1,\ldots,f_{n-1})(x) \lesssim (\prod_{i=1}^{n-1} \M |f_i|^{p_i/p'_n})^{p'_n/p_i}(x),$$
at which point the claim follows from the Hardy-Littlewood maximal inequality.
\end{remark}

To prove Theorem \ref{reduct}, it suffices to prove the following ``restricted weak-type''
analogue.  For any measurable $E \subset \R$, let $X(E)$ denote the space of functions supported on $E$ which are bounded in magnitude by $1$.

\begin{theorem}[Second reduction]\label{reduct-0} Let $n \geq 3$, and let $\Sigma$ and $T^*$ be as in Theorem \ref{reduct}.  Let $E_1,\ldots,E_n$ be subsets of $\R$ of positive finite measure.  Let $p_1,\ldots,p_n$ be such that
$1 < p_1,\ldots,p_{n-1} < 2$, $\frac{5}{2}-n < 1/p_n < 3-n$, and
$$ \frac{1}{p'_n} = \frac{1}{p_1} + \ldots + \frac{1}{p_{n-1}}.$$
Then there exists a subset $E'_n$ of $E_n$ with $|E'_n| \geq \frac{1}{2} |E_n|$ such that one has
$$ |\int T^*(g_1,\ldots,g_{n-1}) g_n| \lesssim |E_1|^{1/p_1} \ldots |E_n|^{1/p_n}$$
for all $g_1 \in X(E_1), \ldots, g_{n-1} \in X(E_{n-1}), g_n \in X(E'_n)$.  Here the implied constant is allowed to depend on $n, p_1,\ldots,p_n$ and $\Sigma$.
\end{theorem}

In the notation of \cite{MTT1}, Theorem \ref{reduct-0} asserts that the $n$-sublinear form $\int T^*(f_1,\ldots,f_{n-1}) f_n$
is of restricted type $(1/p_1,\ldots,1/p_n)$ with $n$ as the bad index.  The deduction of Theorem \ref{reduct} from Theorem \ref{reduct-0} follows from a variant of the Marcinkiewicz interpolation theorem and is a minor modification of the argument in \cite[Lemma 3.11]{MTT1}; the details will be omitted here.  The point of Theorem \ref{reduct-0} is that the functions $g_1,\ldots,g_n$ have been normalized,
indeed $g_j$ can be thought of as essentially the indicator function of $E_j$ (or $E'_n$ when $j=n$).

By a limiting argument we may take $E_1,\ldots,E_n$ to be finite unions of intervals, and $g_1,\ldots,g_n$ to be smooth;
this allows us to justify a number of formal computations in the sequel without difficulty, and we shall do so without any further comment.

To prove Theorem \ref{reduct-0}, we may apply a rescaling argument to normalize $|E_n| = 1$.  From the Hardy-Littlewood maximal inequality we may then set 
\begin{equation}\label{enp-def}
E'_n := E_n \backslash \Omega
\end{equation}
where $\Omega$ is the exceptional set
\begin{equation}\label{omega-def}
\Omega := \bigcup_{i=1}^n \{ \M 1_{E_i} \geq C |E_i| \}
\end{equation}
for a sufficiently large absolute constant $C$, so that $|E'_n| \sim 1$.  It is convenient to renormalize for each $i<n$ $\alpha_i := 1/p_i - 1/2$ and $f_i := g_i / |E_i|^{1/2}$, thus $f_i$ lives in the $L^2$-normalized space $X_2(E_i)$ of functions supported on $E_i$ and bounded in magnitude by
$1/|E_i|^{1/2}$.  We also set $\alpha_n := \frac{n-2}{2} - \alpha_1 -\ldots-\alpha_{n-1}$, thus $0 < \alpha_n < 1/2$.
Theorem \ref{reduct-0} now reduces to

\begin{theorem}[Third reduction]\label{reduct-00} Let $n \geq 3$, and let $\Sigma$ and $T^*$ be as in Theorem \ref{reduct}.  Let $E_1,\ldots,E_n$ be finite unions of intervals with $|E_n|=1$, and let $E'_n$ be defined by \eqref{enp-def}, \eqref{omega-def}, so that $|E'_n| \sim 1$.
Then one has
$$ |\int T^*(f_1,\ldots,f_{n-1}) f_n| \lesssim |E_1|^{\alpha_1} \ldots |E_{n-1}|^{\alpha_{n-1}}$$
for all smooth $f_1 \in X_2(E_1), \ldots, f_{n-1} \in X_2(E_{n-1}), f_n \in X_2(E'_n)$ and any $0 < \alpha_1,\ldots,\alpha_n < 1/2$ with
$\alpha_1 + \ldots + \alpha_n = \frac{n-2}{2}$.  The implied constant can depend on $n, \alpha_1,\ldots,\alpha_n, \Sigma$.
\end{theorem}

This reduction is slightly more convenient to work in as the $L^2$ normalization of $f_1,\ldots,f_n$ will be useful
for a certain ``(maximal) Bessel inequality'' which is crucial to a later stage of the argument.

\section{Fourier representation}\label{Fourier-sec}

Our task is now to prove Theorem \ref{reduct-00}.   
As in \cite{La}, we begin by replacing the rather rough truncation in \eqref{tstar} by a smoother one which has a more tractable Fourier representation.  As is customary, for any $f \in L^1(\R)$, we define the Fourier transform
$$ \hat f(\xi) := \int_\R e^{-2\pi \sqrt{-1} x \xi} f(x)\ dx$$
and the inverse Fourier transform
$$ \check f(x) := \int_\R e^{2\pi \sqrt{-1} x \xi} f(\xi)\ d\xi;$$
we use $\sqrt{-1}$ here instead of $i$ in order to free up the letter $i$ for use as an integer-valued index. 

Let us fix the hyperplane $\Sigma$.
We view the hyperplane $\Sigma$ as an $n-2$-dimensional Euclidean space with Lebesgue measure $d\vec{t}$, 
and thus endowed with its own Fourier transform; thus if $\theta$ is a Schwartz
function on $\Sigma$ we have the inverse Fourier transform
$$ \check{\theta}(\vec t) := \int_\Sigma e^{2\pi \sqrt{-1} \vec t \cdot \vec \xi} \theta(\vec \xi)\ d\xi.$$
We now introduce the multilinear operator
$$ T_{\theta}(f_1,\ldots,f_{n-1})(x) := \int_{\Sigma} (\prod_{j=1}^{n-1}f_j(x+t_i)) \check{\theta}(\vec{t})d\vec{t};$$
this operator can also be written in Fourier space as
$$
T_{\theta}(f_1,\ldots,f_{n-1}) 
 = C_\Sigma \int_{\R^{n-1}} (\prod_{j=1}^{n-1}\hat{f_j}(\xi_j)) 
\theta(\pi(\xi)) e^{2\pi \sqrt{-1} x(\xi_1+\ldots\xi_{n-1})}d\vec{\xi},
$$
where $\pi: \R^{n-1} \to \Sigma$ is the orthogonal projection onto $\Sigma$ and $C_\Sigma > 0$ is a normalization constant depending only on $\Sigma$.
For any integer $k$, write $\theta_k(\xi) := \theta(2^k \xi)$.
We define the associated maximal function $T_\theta^*$ as
$$T_{\theta}^{*}(f_1,\ldots,f_{n-1})(x) := \sup_{k\in\Z} |T_{\theta_k}(f_1,\ldots,f_{n-1})(x)|$$

We shall deduce Theorem \ref{reduct-00} from

\begin{theorem}[Fourth reduction]\label{reduct-2}  Let $n \geq 3$, and let $\Sigma$ be as in Theorem \ref{reduct}.
$0 < \alpha_1,\ldots,\alpha_n < 1/2$ with $\alpha_1 + \ldots + \alpha_n = \frac{n-2}{2}$.
Let $\theta$ be a smooth function supported on a ball $\{ \xi \in \Sigma: \| \xi \| \leq 4 \}$ which
is constant on a ball $\{ \xi \in \Sigma: \| \xi \| \leq 1/4 \}$, and obeys the estimate
\begin{equation}
\label{e:fieq}
|\check{\theta}(t)| \lesssim \frac{1}{(1+\|t\|)^{N^3}} \hbox{ for all } t \in \Sigma
\end{equation}
for some large integer $N$ depending on $\alpha_1,\ldots,\alpha_n$.
Let $E_1,\ldots,E_n$ be finite unions of intervals with $|E_n|=1$, 
and let $E'_n$ be defined by \eqref{enp-def}, \eqref{omega-def}.
Then one has
$$ |\int T^*_\theta(f_1,\ldots,f_{n-1}) f_n| \lesssim |E_1|^{\alpha_1} \ldots |E_{n-1}|^{\alpha_{n-1}}$$
for all smooth $f_1 \in X_2(E_1), \ldots, f_{n-1} \in X_2(E_{n-1}), f_n \in X_2(E'_n)$.  The implied constant can depend on $\Sigma, \alpha_1,\ldots,\alpha_{n-1},N$ and on the implicit constant in \eqref{e:fieq}.
\end{theorem}

\begin{proof}[of Theorem \ref{reduct-00} assuming Theorem \ref{reduct-2}] 
We may take $f_1,\ldots,f_{n-1}$ non-negative.

Let $\eta$ be a fixed real-valued symmetric 
Schwartz function\footnote{Such a function can be constructed by starting with a real-valued symmetric function on the
ball $\{ \| \xi \| \leq 1/2 \}$, then convolving it with itself and normalizing it.}
 on $\Sigma$ supported on the ball $\{ \| \xi \| \leq 1 \}$ 
whose Fourier transform is non-negative and $\check{\eta}(0) = 1$.  Observe that
$$ T_{\eta_k}(f_1,\ldots,f_{n-1})(x) := \frac{1}{2^{k(n-2)}} \int_{\Sigma} (\prod_{j=1}^{n-1}f_j(x+t_i)) \check{\eta}
(\vec{t}/2^k)d\vec{t}.$$ 
From this, the positivity of the $f_i$ and $\check \eta$ it is easy to establish the pointwise estimate
$$ T^{*}(f_1,\ldots,f_{n-1})(x) \lesssim T^*_\eta(f_1,\ldots,f_{n-1})(x)$$
(where the implied constant depends on $\eta$)
so it suffices to show that
$$ |\int T^*_\eta(f_1,\ldots,f_{n-1}) f_n| \lesssim |E_1|^{\alpha_1} \ldots |E_{n-1}|^{\alpha_{n-1}}.$$
We cannot yet apply Theorem \ref{reduct-2}, because $\eta$ is not constant near the origin.  Indeed the requirement that $\check \eta$ be non-negative forces $\eta$ to have a negative Laplacian at the origin.  Fortunately, we can rectify this by a a further dyadic decomposition.
More precisely, we split
$$
\eta(\xi)=\eta_2(\xi)+\sum_{l=-\infty}^{0} \phi_l(\xi)
$$
with $\eta_2$ smooth, symmetric, supported in $\|\xi\|\le 11/10$ and equal to $1$ on $\|\xi\|\le 1$, while 
$\phi_l(\xi) := (\eta-\eta_2)(\xi) (\eta_2(\xi/2^l) - \eta_2(\xi/2^{l-1}))$.  
One can easily verify that the function $\eta_2$ is already 
of the form required for Theorem \ref{reduct-2} and so $T^*_{\eta_2}$ gives an acceptable contribution to $T^*_\eta$.
As for the tail terms $\phi_l$, we observe the Fourier estimates
$$|\frac1{2^l}\check{\phi_l}(\frac{\xi}{2^l})|\lesssim 2^{-|l|}\frac{1}{(1+\|\xi\|)^{N^3}}
$$ 
uniformly in $l$. Also, $\phi_l$ is constant on $\|\xi\|\le 2^l/4$ and zero when $\|\xi\|\ge 4 \times 2^l$.  
A simple rescaling argument using Theorem \ref{reduct-2} (noting that $T^*_\phi$ is unchanged if one replaces $\phi$ by $\phi(2^l \cdot)$) then shows that
$$ |\int T^*_{\phi_l}(f_1,\ldots,f_{n-1}) f_n| \lesssim 2^{-|l|} |E_1|^{\alpha_1} \ldots |E_{n-1}|^{\alpha_{n-1}}.$$
The claim now follows from the triangle inequality.
\end{proof}

\section{Discretization}\label{sec:3}

It remains to prove Theorem \ref{reduct-2}.  We now perform the usual dyadic decompositions to reduce matters to estimating a certain sum over dyadic objects, namely a collection of ``multitiles'', after first doing some additional refinements to ensure that these multitiles obey some good geometrical properties (specifically, a rank one condition).

We introduce two large constants $1 \ll C_0 \ll C_1$ (depending on $\Sigma$, and $C_1$ assumed to be large compared to $C_0$) 
that will be used to sparsify the time-frequency geometry.  We will take some care to specify how the implied constants in the $\lesssim$
notation depend on $C_0$ and $C_1$; however we will allow these constants to depend freely on $n, \alpha_1,\ldots,\alpha_n,N, \Sigma$.

It will be convenient to dilate $\theta$ by $C_0$,
so that $\theta$ is now supported on $\{ \xi \in \Sigma: \| \xi \| \leq 4C_0 \}$ which
is constant on a ball $\{ \xi \in \Sigma: \| \xi \| \leq C_0/4 \}$; this affects our final bounds by some factor depending on $C_0$, but
as we shall eventually choose $C_0$ to be a quantity depending only on existing parameters such as $n,\alpha_1,\ldots,\alpha_n,N,\Sigma$, this shall
be of no consequence.  We perform the dyadic decomposition
$$\theta_k(\xi)=\sum_{i\ge k}\varphi_i(\xi)
$$ 
where $\varphi(\xi):=\theta(\xi)-\theta(2\xi)$ is a smooth function  
supported an annulus $\|\xi\| \sim C_0$, and $\varphi_i(\xi) := \varphi(2^i \xi)$. 
Thus
$$ T_{\theta_k} = \sum_{i \ge k} T_{\varphi_i}$$
and hence for any $f_1,\ldots,f_n$ we have
$$ \int T^*_\theta(f_1,\ldots,f_{n-1}) f_n = \sum_i \int T_{\varphi_i}(f_1,\ldots,f_{n-1})(x) f_n(x) 1_{i \geq k(x)}\ dx$$
for some integer-valued measurable function $k: \R \to \Z$.  Thus it suffices to establish the multilinearized estimate
\begin{equation}\label{decomp-1} 
|\sum_i \int T_{\varphi_i}(f_1,\ldots,f_{n-1})(x) f_n(x) 1_{i \geq k(x)}\ dx| \lesssim_{C_0,C_1}
|E_1|^{\alpha_1} \ldots |E_{n-1}|^{\alpha_{n-1}}
\end{equation}
for each such function $k: \R \to \Z$, which we now fix.  Note that we can write the left-hand side as
\begin{equation}\label{decomp-2}
\int_\R \int_{\Sigma} (\prod_{j=1}^n f_{j,i}(x+t_i)) \check{\varphi_i}(\vec{t})d\vec{t}\ dx
\end{equation}
where $f_{j,i} := f_j$ for $1 \leq j \leq n-1$ and $f_{n,i}(x) := f_n(x) 1_{i \geq k(x)}$, and we adopt the convention that $t_n = 0$.  One should think of $i$ as a scale parameter, corresponding to the terms with frequency uncertainty $\sim 2^{-i}$ and time uncertainty $\sim 2^{i}$.  Note that the annulus that $\varphi_i$ is supported in has thickness $\sim C_0 2^{-i}$ and can thus tolerate the frequency uncertainty associated to the scale $i$.

The next (standard) step is wave packet decomposition.  We shall adopt the usual trick of covering the time domain $\R$ by three overlapping
dyadic grids to eliminate some artificial boundary effects caused by dyadicity.

For each $1 \leq j \leq n$, let us pick a Schwartz function $\psi_j$ such that $\hat{\psi_j}$ is supported in $[0.1,0.9]$, and
that $\psi_j$ is rapidly decreasing; in particular we have the bounds
\begin{equation}\label{psij}
 |\psi_j(x)| \lesssim (1 + |x|)^{-10N} \hbox{ for all } x \in \R
\end{equation}
and we have the following property for every $\xi\in\R$:
$$\sum_{l\in \Z}\left|\hat{\psi_j}\left(\xi-\frac{l}3\right)\right|^2=1.
$$
This is possible because the translates of $[0.1,0.9]$ by integer multiples of $\frac{1}{3}$ cover the real line $\R$ with
some room to spare for smooth cutoffs.
For each scale $i \in \Z$ we can then decompose
$$f_{i,j}=\sum_{m,l\in \Z}\langle f_{i,j},\psi_{j,i,m,\frac{l}3}\rangle\psi_{j,i,m,\frac{l}3},$$
where 
$$\psi_{j,i,m,l}(x) := 2^{-\frac{i}{2}}\psi_j(2^{-i}x-m)e^{2\pi \sqrt{-1}2^{-i}xl}$$
and $\langle f, g \rangle := \int f\overline{g}$ is the usual inner product.  Inserting this decomposition into
\eqref{decomp-1}, \eqref{decomp-2} and using the triangle inequality, we reduce to showing that
\begin{equation}
\label{e:23765}
\sum_{i\in \Z}\sum_{\vec{m},\vec{l}\in\Z^n}C_{\vec{m},\vec{l},i} 
2^{i(1-\frac{n}{2})}\prod_{j=1}^{n}|\langle f_{j,i},\psi_{j,i,m_j,l_j} \rangle|
\lesssim_{C_0,C_1} |E_1|^{\alpha_1} \ldots |E_{n-1}|^{\alpha_{n-1}}
\end{equation}
where $\vec m = (m_1,\ldots,m_n)$, $\vec l = (l_1,\ldots,l_n)$, and $C_{\vec{m},\vec{l},i}$ are the operator coefficients
\begin{equation}\label{cdef}
 C_{\vec{m},\vec{l},i}
:=\frac{1}{2^{i(1-\frac{n}{2})}}|\int_\R \int_\Sigma \prod_{j=1}^{n}\psi_{j,i,m_j,l_j}(x+t_j) 
 \check{\theta_i}(\vec{t}) d\vec{t} dx|. 
\end{equation}
One should think of $\vec m$ as containing the time location, and $\vec l$ as containing the frequency location
information; roughly speaking, the summand in \eqref{e:23765} is
the contribution when $f_j$ is localized in space to $2^{i} m_j + O(2^{i})$ and localized in frequency to $2^{-i} l_j + O(2^{-i})$.

We now use the geometry of the hyperplane $\Sigma$ to obtain localization estimates on the coefficients $C_{\vec{m}, \vec{l}, i}$.
We let $\Gamma \subset \R^n$ denote the hyperplane $\Gamma:=\{(\xi_1,\ldots,\xi_n):\xi_1+\ldots+\xi_n=0\}$.

\begin{lemma}  We have the estimate
\begin{equation}\label{cbound}
C_{\vec{m},\vec{l},i} \lesssim_{C_0,C_1} (1+\diam\{m_1,\ldots, m_n\})^{-N^2}.
\end{equation}
Furthermore, if $C_{\vec{m},\vec{l},i}$ is non-zero, then 
\begin{equation}
\label{eq:41}
l_1 + \ldots + l_n = O(1)
\end{equation}
and
\begin{equation}
\label{eq:1749}
\| \pi(l_1,\ldots,l_{n-1}) \| \sim C_0.
\end{equation}
\end{lemma}

\begin{remark} In the notation of \cite{MTT1}, these conditions are essentially asserting that the tuples $(\vec{m}, \vec{l}, i)$ with a sizeable coefficient $C_{\vec{m},\vec{l},i}$ form a collection of multitiles of rank one (which is also the situation with the bilinear Hilbert transform).  See also Definition \ref{rankdef} below.
\end{remark}

\begin{proof}  We first observe that by rescaling by $2^i$ that $C_{\vec{m}, \vec{l}, i}$ is actually independent of $i$.  Thus we may assume $i=0$ throughout the proof.

To prove \eqref{cbound}, we then use the physical space representation 
\eqref{cdef} of $C_{\vec{m},\vec{l},0}$, followed by the triangle inequality, to obtain
$$
 C_{\vec{m},\vec{l},0}
\lesssim \int_\R \int_\Sigma \left| (\prod_{j=1}^n \psi_j(x+t_j-m_j)) \check{\theta_0}(\vec{t})\right|d\vec{t} dx.
$$
Now as $\psi$ is rapidly decreasing, we conclude from \eqref{e:fieq} that
$$
 C_{\vec{m},\vec{l},0}
\lesssim_{C_0,C_1} \int_\R \int_\Sigma (1+\|t\|)^{-2N^2} \prod_{j=1}^n (1 + |x+t_j-m_j|)^{-2N^2} d\vec{t} dx
$$
(say), and the claim \eqref{cbound} follows from the pointwise estimate
$$ \prod_{j=1}^n (1 + |x+t_j-m_j|)^{-N^2} \lesssim (1+\diam\{m_1,\ldots, m_n\})^{-N^2} (1+\|t\|)^{N^2}.$$

Now suppose that $C_{\vec{m},\vec{l},i}$ is non-zero.  To exploit this we use the Fourier representation, converting \eqref{cdef} to
$$
 C_{\vec{m},\vec{l},0} \sim
| \int_{\Gamma} (\prod_{j=1}^{n}\hat{\psi}(\xi_j-l_j))  \theta(\pi(\xi_1,\ldots,\xi_{n-1}))\ d\vec{\xi}|.
$$
Thus there exists $\vec \xi \in \Gamma$ such that $\xi_j-l_j$ lies in the support of $\hat \psi$ for all $1 \leq j \leq n$
and $\pi(\xi_1,\ldots,\xi_{n-1})$ lies in the support of $\theta$.  From the former property we have $l_j = \xi_j + O(1)$, and
\eqref{eq:41} follows from the definition of $\Gamma$.  From the latter property we have $\| \pi(\xi_1,\ldots,\xi_{n-1})\|\sim C_0$, and the claim follows by using the approximation $l_j = \xi_j + O(1)$ and the homogeneity of $\pi$. 
\end{proof}

In view of the above lemma, it now suffices to show that
\begin{equation}\label{summarize}
\sum_{(\vec{m},\vec{l},i) \in \Omega} (1+\diam\{m_1,\ldots, m_n\})^{-100n^2}
2^{i(1-\frac{n}{2})}\prod_{j=1}^{n}|\langle f_{j,i},\psi_{j,i,m_j,l_j} \rangle|
\lesssim_{C_0,C_1} |E_1|^{\alpha_1} \ldots |E_{n-1}|^{\alpha_{n-1}}
\end{equation}
where $\Omega$ is a collection of triples $(\vec{m}, \vec{l}, i) \in \Z^n \times \Z^n \times \Z$ obeying \eqref{eq:41} and
\eqref{eq:1749}.

We now perform a number of refinements to improve the nesting properties of the set $\Omega$.  First we observe 
that for each $(\vec{m},\vec{l},i) \in \Omega$ and $1 \leq j \leq n$, the Fourier transform of $\psi_{j,i,m_j,l_j}$
is contained in the interval $[2^{-i} \frac{l_j}{3}, 2^{-i} (\frac{l_j}{3}+1)]$ (in fact they are contained in the slightly smaller interval $[2^{-i} (\frac{l_j}{3} + 0.1), 2^{-i} (\frac{l_j}{3}+0.9)]$).  These intervals are almost dyadic, but for
the denominator of $3$.  However this factor of $3$ can be eliminated in the following standard manner.  Let ${\mathcal D}_0, {\mathcal D}_1, {\mathcal D}_2$ be the dyadic grids
\begin{equation}\label{dyad}
\begin{split}
{\mathcal D}_0 &:= \{ [ 2^{-i} l, 2^{-i}(l+1)]: i, l \in \Z \} \\
{\mathcal D}_1 &:= \{ [ 2^{-i} (l + (-1)^i/3), 2^{-i}(l+1 + (-1)^i/3)]: i, l \in \Z \} \\
{\mathcal D}_2 &:= \{ [ 2^{-i} (l - (-1)^i/3), 2^{-i}(l+1 - (-1)^i/3)]: i, l \in \Z \}.
\end{split}
\end{equation}
Thus ${\mathcal D}_0$ is the standard dyadic grid, and the other two grids are essentially similar (one can view the latter two grids as translates of the first by the non-terminating $2$-adic $\pm 1/3$).  In particular, within a single grid we have the nesting property that if two intervals intersect, then the shorter one is contained by the longer one.
Observe that every interval $[2^{-i} \frac{l_j}{3}, 2^{-i} (\frac{l_j}{3}+1)]$ belongs
to one of these three grids.  By pigeonholing once for each $j$ (conceding a factor of $3^n$ in the estimates), we can assume that for fixed $j$,
the intervals $[2^{-i} \frac{l_j}{3}, 2^{-i} (\frac{l_j}{3}+1)]$ belong to a single dyadic grid.   For ease of exposition we shall assume
that these intervals always lie in the standard dyadic grid ${\mathcal D}_0$, thus the intervals $[2^{-i} \frac{l_j}{3}, 2^{-i} (\frac{l_j}{3}+1)]$ 
are genuine dyadic intervals.  The other cases are handled similarly but with some minor changes in notation.

Morally speaking, the localizing factor $(1+\diam\{m_1,\ldots, m_n\})^{-N^2}$ in \eqref{summarize} implies that the 
diagonal contribution $m_1 = \ldots = m_n$ is the dominant contribution.  Again to simplify the exposition, we shall focus
entirely on this diagonal case $m_1 = \ldots = m_n$.  We now briefly sketch how to pass from the diagonal case to the general case. Write $m_j = m_1 + r_j$.  For each fixed $n-1$-tuple of integers $r_2,\ldots,r_n$, one can convert the case $m_j = m_1 + r_j$ to the diagonal case $m_j = m_1$ by shifting the function $\psi_j$ by $r_j$.  This affects the bounds \eqref{psij} but only by $(1+|r_j|)^{10N}$ at worst.  This gives a total loss of $\prod_{j=1}^n (1 + |r_j|)^{10nN}$ for this contribution, but one is
also gaining a factor of $(1 + \diam(0,r_2,\ldots,r_n))^{-N^2}$, and the product is then summable in $r$ if $N$ is large enough.  Thus it suffices to treat the diagonal case.

Another application of the pigeonhole principle (giving up a constant factor of $C_1$ in the estimates)
allows one to refine the scale parameter $i$ to not take values in the  integers, but to instead take values in a 
residue class $\{ i = c \mod C_1 \}$ for some residue $c$.  This ``sparsification'' of the scales will be useful in obtaining a certain rank separation condition in the frequencies below.

Finally, we analyze the conditions \eqref{eq:41} and \eqref{eq:1749}.  Observe that if we instead had the exact
constraints $l_1 + \ldots + l_n = 0$ and $\pi(l_1,\ldots,l_{n-1})=0$, then $(l_1,\ldots,l_n)$ would be restricted to a one-dimensional subspace of $\R^n$.  Since $\Sigma$ did not contain $e_1,\ldots,e_{n-1}$ or $(1,\ldots,1)$, it is easy to see
that the non-zero vectors in this one-dimensional subspace have no zero coordinates; thus we have $l_j = c_{j,j'} l_{j'}$ for all
$1 \leq j,j' \leq n$ and some explicit non-zero finite constants $c_{j,j'}$ depending only on $\Sigma$; furthermore we have $c_{j,j} = 1$, $c_{1,j'} + \ldots + c_{n,j'} = 0$ and $\pi( c_{1,j'},\ldots,c_{n-1,j'}) = 0$.  Returning now to the inexact constraints \eqref{eq:41}, \eqref{eq:1749}, we conclude that
$$ l_j = c_{j,j'} l_{j'} + O( C_0 )$$
for all $1 \leq j,j' \leq n$.
By pigeonholing (and conceding a factor of $C_0^{n^2}$ at worst) we may thus assume that
\begin{equation}\label{lll}
l_j = \lfloor c_{j,j'} l_{j'} \rfloor + a_{j,j'}
\end{equation}
on $\Omega$ for all $1 \leq j,j' \leq n$ and some fixed integers $a_{j,j'} = O(C_0)$; note that $a_{j,j}$ is necessarily zero.  
Thus each frequency $l_j$ 
is now uniquely determined by any of the other frequencies $l_{j'}$.  Furthermore, from \eqref{eq:41}, \eqref{eq:1749} we have
$$ a_{1,j'} + \ldots + a_{n,j'} = O(1) \hbox{ and } \| \pi( a_{1,j'},\ldots,a_{n-1,j'} ) \| \sim C_0.$$
If $C_0$ is large enough, this implies the following basic fact:

\begin{lemma}\label{rank-cond}  For each $j'$, there exist at least two $j$ distinct from $j'$ such that $|a_{j,j'}| \sim C_0$.
\end{lemma}

The upshot of this lemma is that whenever we fix one of the frequencies of $f_1,\ldots,f_n$, at least two other frequencies depend in a ``lacunary'' manner on the scale parameter $i$.  This fact will be crucial in controlling the geometry of certain ``trees'' which will appear later.

The estimate \eqref{summarize} has now been reduced to
\begin{equation}\label{summarize-2}
\sum_{(\vec{m},\vec{l},i) \in \Omega} 
2^{i(1-\frac{n}{2})}\prod_{j=1}^{n}|\langle f_{j,i},\psi_{j,i,m_j,l_j} \rangle|
\lesssim |E_1|^{\alpha_1} \ldots |E_{n-1}|^{\alpha_{n-1}}
\end{equation}
We now convert \eqref{summarize-2} into the more traditional language of multitiles and wave packets. 

\begin{definition}[Tiles]
A \emph{tile} $P$ is a rectangle $P=I_P\times\omega_{P}$ with both $I_P$ and $\omega_P$ dyadic intervals, obeying the Heisenberg relation $|I_P|\cdot|\omega_P|=1$; we
refer to $I_P$ as the \emph{time interval} of $P$ and $\omega_P$ as the \emph{frequency interval}. A \emph{multitile} $s$ is an $n$-tuple $s = (s_1,\ldots,s_n)$ of tiles with the same time interval $I_s := I_{s_1} = \ldots = I_{s_n}$.
If $I$ is an interval and $C > 0$ is a number, we let $CI$ denote the interval with the same center as $I$ but $C$ times the length (note that this interval will most likely not be dyadic).
Let us say that a function $\psi_P$ is a \emph{wave packet adapted to a tile} $P$ if $\hat \psi_P$ is supported in $0.8 \omega_P$
 and we have the pointwise estimate
\begin{equation}\label{psisj}
 |\psi_P(x)| \lesssim |I_P|^{-1/2} \chi_{I_P}^{10N}(x) \hbox{ for all } x \in \R
\end{equation}
where for any interval $I$, $\chi_I$ is the weight function
$$ \chi_I(x) := (1 + \frac{(x - c(I))^2}{|I|^2})^{-1/2}$$
and $c(I)$ is the center of $I$; in particular observe that $\psi_P$ is normalized to have an $L^2$ norm of $O(1)$.  
\end{definition}

Note that because of all the reductions we have already achieved, every triple $(\vec{m},\vec{l},i)$ in $\Omega$ gives rise to
a multitile $s$ with $s_j := [2^i m_j, 2^i (m_j+1)] \times \omega_{s_j} := [2^{-i} \frac{l_j}{3}, 2^{-i} (\frac{l_j}{3}+1)]$.
In particular we have $|I_s| = 2^i$.  
Let $\S_{\max}$ denote the collection of all multitiles obtained this way.  
For each multitile $s \in \S_{\max}$ arising from a triple $(\vec{m}, \vec{l}, i)$, define the functions $\psi_{s,j}$ for $1 \leq j \leq n$ by setting
$$ \psi_{s,j}(x) := \psi_{j,i,m_j,l_j}(x).$$
Observe that for each $1 \leq j \leq n$, $\psi_{s,j}$ is a wave packet adapted to $s_j$.
We also observe the following important consequence of Lemma \ref{rank-cond}.  

\begin{definition}[Rank one]\label{rankdef}  A collection $\S$ of multitiles is said to have \emph{rank one} if
for every $j \in \{1,\ldots,n\}$ there exists distinct $j_1(j), j_2(j) \in \{1,\ldots,n\} \backslash \{j\}$ and signs $\epsilon_1(j), \epsilon_2(j) \in \{-1,+1\}$ (not necessarily distinct) with the following properties.
\begin{itemize}
\item (Scale separation) If $s, s' \in \S$ are such that $|\omega_{s_j}| > |\omega_{s'_j}|$, then $|\omega_{s_j}| \geq 2^{C_1} |\omega_{s',j}|$.
\item (One independent frequency parameter) If $s, s' \in \S$ are such that $\omega_{s_j} = \omega_{s'_j}$, then $\omega_{s_{j'}}=\omega_{s'_{j'}}$ for all $1 \leq j' \leq n$.
\item (Nearby $j$-frequencies implies nearby $j'$-frequencies) If $s, s' \in \S$ are such that $10\omega_{s,j} \cap 10\omega_{s',j} \neq \emptyset$ and $|I_s| \geq |I_{s'}|$, then $\dist(\omega_{s,j'}, \omega_{s',j'}) \lesssim C_0 |I_{s'}|^{-1}$ for 
all $1 \leq j' \leq n$.
\item (Lacunarity property) If $s, s' \in \S$ are such that $10\omega_{s_j} \cap 10\omega_{s',j} \neq \emptyset$ and
$|I_s| > |I_{s'}|$, then $\dist(\omega_{s,j_t}, \omega_{s',j_t}) \sim C_0 |I_{s'}|^{-1}$ for $t=1,2$.  In particular $10\omega_{s_{j_t}}$ and $10\omega_{s'_{j_t}}$ are disjoint.  Furthermore we require $\epsilon_t(j) (\xi' - \xi) \geq 0$ for all $\xi \in 10 \omega_{s,j_t}$ and $\xi' \in 10\omega_{s',j_t}$.
\end{itemize}
\end{definition}

\begin{remark} For the definition of higher order rank (which we will not need here), see \cite{MTT1}.  Actually our definition of rank one is slightly stronger than that in \cite{MTT1} in that we require the indices $j_1,j_2, \epsilon_1(j), \epsilon_2(j)$ to depend only on $j$, and not be dependent on $s, s'$, but this is only a minor technical change.
\end{remark}

\begin{lemma}[Rank one property]\label{rank-2}  $\S_{\max}$ has rank one.
\end{lemma}

This lemma shows, among other things, that the multitiles in $\S_{\max}$ have essentially one independent frequency parameter.  Note that if $\S_{\max}$ has rank one, then so does any subset $\S$ of $\S_{\max}$.

\begin{proof} The scale separation property follows since for each multitile $s \in \S_{\max}$, we have $|\omega_{s_j}| = |I_s|^{-1} = 2^{-i}$
for all $1 \leq j \leq n$ and some integer $i= C \mod C_1$.  The remaining properties follow from \eqref{lll} and
Lemma \ref{rank-cond}, setting $j_1,j_2$ to be the indices distinct from $j$ such that $|a_{j_1,j}|, |a_{j_2,j}| \sim C_0$,
and $\epsilon_t(j)$ to be the sign of $a_{j_t,j}$.
\end{proof}

We also define the modified wave packets $\phi_{s,j}$ by setting 
\begin{equation}\label{phij-def}
\phi_{s,j} := \psi_{s,j} \hbox{ for } 1 \leq j \leq n-1
\end{equation}
and 
\begin{equation}\label{phij2-def}
 \phi_{s,n}(x) := \psi_{s,n}(x) 1_{|I_s| > 2^{k(x)}}.
\end{equation}
The estimate \eqref{summarize-2} can now be rewritten as
$$
\sum_{s \in \S_{\max}} 
|I_s|^{(1-\frac{n}{2})} \prod_{j=1}^{n} |\langle f_j,\phi_{s,j} \rangle| \lesssim |E_1|^{\alpha_1} \ldots |E_{n-1}|^{\alpha_{n-1}}.$$
By the monotone convergence theorem we can replace $\S_{\max}$ by a finite subset $\S$ of $\S_{\max}$, so long as our estimates are uniform in $\S$.  Note that the properties in Lemma \ref{rank-2} will be preserved if we pass from $\S_{\max}$ to $\S$.
We can now deduce \eqref{summarize-2} (and hence Theorem \ref{thm:ct1}) from the following more abstract result.

\begin{theorem}[Fifth reduction]\label{reduct-5}  Let $n \geq 3$, let
$0 < \alpha_1,\ldots,\alpha_n < 1/2$ with $\alpha_1 + \ldots + \alpha_n = \frac{n-2}{2}$, and let $N$
be a sufficiently large integer depending on $\alpha_1,\ldots,\alpha_n$.
and let $\S$ be a finite collection of multitiles which is rank one.  For each $s \in \S$ and $1 \leq j \leq n$, let $\psi_{s,j}$ be a wave packet adapted to $s_j$.  Let $k: \R \to \Z$ be an arbitrary measurable function, 
and let $\phi_{s,j}$ be defined by \eqref{phij-def}, \eqref{phij2-def}.  
Let $E_1,\ldots,E_n$ be finite unions of intervals with $|E_n|=1$, and let $E'_n$ be defined by
\eqref{enp-def}, \eqref{omega-def}.
Then one has
$$
\sum_{s \in \S} 
|I_s|^{(1-\frac{n}{2})} \prod_{j=1}^{n} |\langle f_j,\phi_{s,j} \rangle| \lesssim |E_1|^{\alpha_1} \ldots |E_n|^{\alpha_n}.$$
for all smooth $f_1 \in X_2(E_1), \ldots, f_{n-1} \in X_2(E_{n-1}), f_n \in X_2(E'_n)$.  The implied constant can depend on $\alpha_1,\ldots,\alpha_{n-1},N$ and on the bounds in the rank one condition and \eqref{psisj} but is uniform in $\S$.
\end{theorem}

\begin{remark} If the $\phi_{s,j}$ were replaced by $\psi_{s,j}$ (i.e. if the cutoff $|I_s| > 2^{k(x)}$ were not present) then
this result would follow from the results in \cite{MTT1}.  Thus the novelty (which is also present in \cite{La}) is
the cutoff $|I_s| > 2^{k(x)}$, which ultimately arises from the maximal function nature of $T^*_{A,\R}$.
\end{remark}

\section{Trees}
\label{sec:4}

It remains to prove Theorem \ref{reduct-5}.  To do this we use the standard strategy of organizing the multitiles into trees,
estimating the contribution of each tree separately, controlling the total number of trees of a certain ``size'', and then summing up.  

Henceforth we fix the tile collection $\S$ and the functions $f_1,\ldots,f_n$ and sets $E_1,\ldots,E_n$, as well as the exponents $\alpha_1,\ldots,\alpha_n$ and $N$, the 
function $k(x)$ and the wave packet functions $\psi_{s,j}$ (which of course determine $\phi_{s,j}$).
We now recall a standard notion of tile order.

\begin{definition}[Tile order]  For any two tiles $P$ and $P'$, we write $P < P'$ if $I_P \subsetneq I_{P'}$ and $3\omega_P \supsetneq 3\omega_{P'}$,
and $P \leq P'$ if $P < P'$ or $P = P'$.
\end{definition}

Note that this is a partial order on tiles.  The factor of $3$ is convenient for technical reasons to provide a little more frequency separation; the presence of the large constants $C_0$ and $C_1$ in the rank condition will allow us to have this additional factor.

\begin{definition}[Trees]  A \emph{multitile tree}, or \emph{tree} for short,
is a triplet $(\T, T, i)$ where $1 \leq i \leq n$ is the \emph{index}
of the tree, $T \in \S$ is a multitile, and $\T \subset \S$ is a collection of multitiles such
that $s_i \le T_i$ for all $s \in \T$.  We shall often abuse notation and abbreviate a tree $(\T,T,i)$ as $\T$.
We refer to $I_\T := I_T$ as the \emph{time interval} of the tree.  If $1 \leq j \leq n$ and $\epsilon \in \{-1,+1\}$, we say
that a tree $(\T,T,i)$ is \emph{$(j,\epsilon)$-separated} if $j = j_t(i)$ and $\epsilon = \epsilon_t(i)$ for some $t\in\{1,2\}$.
We say that a tree is \emph{$j$-separated} if it is $(j,\epsilon)$-separated for some $\epsilon \in \{-1,+1\}$.
\end{definition}

\begin{example} For any tile $T$ and $1 \leq i \leq n$, the singleton tree $(\{T\}, T, i)$ is a multitile tree.
\end{example}

\begin{remark} We use the rather clumsy terminology \emph{multitile tree} to distinguish from the notion of a \emph{lacunary tree}, which consists of tiles rather than multitiles, that we will introduce in Section \ref{bmo-sec}.
Note that we do not require that the tree $\T$ contains its top $T$, although this is often the case; also note that if $(\T,T,i)$ is a tree then so
is $(\T \cup \{T\},T,i)$ and $(\T \backslash \{T\},T,i)$ (so one can always add or remove the top from a tree).  
This additional flexibility in our definition
of tree (not present in some other literature) is convenient because it makes the notion of tree more stable with respect to passage to subsets.  In particular, if $(\T,T,i)$ is a $(j,\epsilon)$-separated tree and $\T' \subset \T$, then $(\T',T,i)$ is also a $j$-separated 
tree.  Furthermore, if $\T'$ takes the form $\T' := \{ s \in \T: s_i \le T'_i \}$ for some multitile $T'$, then $(\T',T',i)$ is also a $j$-separated tree.
\end{remark}

The rank one condition implies certain geometric facts about trees, which we collect below for the reader's convenience.

\begin{lemma}\label{quasi} Let $(\T,T,i)$ and $(\T', T', i)$ be $(j,\epsilon)$-separated multitile trees.
\begin{itemize}
\item[(i)] The frequency intervals of a multitile in $\T$ are determined entirely by the size of the spatial interval.  In other words, if
$s, s' \in \T$ and $|I_s| = |I_{s'}|$, then $\omega_{s_k} = \omega_{s'_k}$ for all $1 \leq k \leq n$. 
\item[(ii)] Each multitile in $\T$ has a distinct time interval: if $s, s' \in \T$ and $s \neq s'$, then $I_s \neq I_{s'}$.
\item[(iii)] If $s \in \T$ and $s \neq T$, then $\dist( 10\omega_{s_j}, 10\omega_{T_j} ) \sim C_0 |I_s|^{-1}$; in particular,
$10\omega_{s_j}$ and $10\omega_{T_j}$ are disjoint.
\item[(iv)] Suppose that $s \in \T$ and $s' \in \T'$ are such that $\omega_{s_j} \subsetneq \omega_{s'_j}$ and $I_{s'} \cap I_T \neq \emptyset$.  
Then $s'_j < T_j$, and furthermore we have
$\epsilon (\xi - \xi') > 0$ whenever $\xi \in \omega_{T_i}$ and $\xi' \in \omega_{T'_i}$.
\end{itemize}
\end{lemma}

\begin{proof}  If $|I_s| = |I_{s'}|$ then $|\omega_{s_j}| = |\omega_{s'_j}|$; since these intervals intersect, we must have
$\omega_{s_j} = \omega_{s'_j}$ and then (i) follows from the rank one condition. Property (ii) follows immediately from (i).
Now we show (iii).  From (ii) we see that $I_s$ is strictly smaller than $I_\T$, and so $\omega_{s_i}$ strictly contains $\omega_{T_i}$.  The claim then follows from lacunarity property of the rank condition.  Finally, we show (iv).  We have
$|\omega_{s_j}| < |\omega_{s'_j}|$ and $|I_{s}| \leq |I_T|$ and hence $|I_{s'}| < |I_T|$.  By dyadic nesting this means
that $I_{s'} \subsetneq I_T$, and to show that $s'_j < T_j$ 
it will suffice to show that $3\omega_{s'_j}$ intersects $3\omega_{T_j}$.  But $\omega_{T_j}$
lies within $\lesssim C_0 |\omega_{s_j}|^{-1}$ of $\omega_{s_j}$, which is contained inside $\omega_{s'_j}$.  Since
$|\omega_{s'_j}| \geq 2^{C_1} |\omega_{s_j}|$ by scale separation, the claim $s'_j < T_j$ follows if $C_1$ is sufficiently large depending on $C_0$.  To show the remaining claim in (iv), we observe from the rank separation condition that $\dist(\omega_{s'_j},\omega_{T'_j}) \sim C_0 |\omega_{s'_j}|$, with $\omega_{T'_j}$ lying below $\omega_{s'_j}$ if $\epsilon = +1$
and above if $\epsilon=-1$.  The claim follows.
\end{proof}




We can now introduce the concept of size.  There will be one size for each of the functions $f_1,\ldots,f_n$.

\begin{definition}[Size]
For a  set of multitiles $\S'\subset\S$ and $1 \leq j \leq n$ define its \emph{$j$-size} as 
$$\size_j(\S') := \sup_{\T}\left(\frac{1}{|I_\T|}\sum_{s\in \T}|\langle f_j, \phi_{s,j}\rangle|^2\right)^{\frac{1}{2}}$$ 
where the supremum is taken over all the $j$-separated trees $(\T,T,i)$ with $\T \subset \S'$.
\end{definition} 

\begin{remark} In the above definition the trees $\T$ are not required to contain their top $T$.  However it is easy to see
that a tree without a top can be partitioned into trees with tops that have disjoint time intervals, and because of this one
could replace the supremum in the definition of size by a supremum over trees that contain their tops without affecting the size.  However
we will not need to do this in this paper.
\end{remark}

%


\section{High-level overview of proof}\label{overview-sec}

Following the usual time-frequency approach, we can now reduce the task of proving Theorem \ref{reduct-5} to that of
verifying a number of lemmas concerning trees.

The first lemma is easy to state and prove:

\begin{lemma}[Contribution of a single tree]
\label{l:328}
If $(\T,T,i_0)$ is a tree then
$$\sum_{s\in \T}|I_s|^{1-\frac{n}2}\prod_{i=1}^{n}|\langle f_i, \phi_{s,i}\rangle|\leq |I_\T|\prod_{i=1}^{n}\size_i(\T).$$
\end{lemma}

\begin{proof}  
By definition of size we have
$$ ( \sum_{s\in \T} |\langle f_i, \phi_{s,i}\rangle|^2 )^{1/2} \leq |I_\T|^{1/2} \size_i(\T)$$
when $i = j_1(i_0)$ or $i = j_2(i_0)$.  Also, since a singleton multitile is always a tree, we also have
$$ |\langle f_i, \phi_{s,i} \rangle| \leq |I_s|^{1/2} \size_i(\T)$$
for the other $n-2$ values of $i$.  The claim then follows from H\"older's inequality.
\end{proof}

In light of this lemma, the task is now to subdivide the collection $\S$ into distinct trees $\T$ for which one has the bound
\begin{equation}\label{tsum}
 \sum_{\T} |I_\T| \prod_{i=1}^{n}\size_i(\T) \lesssim |E_1|^{\alpha_1} \ldots |E_{n-1}|^{\alpha_{n-1}}.
\end{equation}

This will be accomplished via a number of propositions.  First we need a basic upper bound on the size of a tree, which we
prove in Section \ref{sec:6}.

\begin{proposition}[Size estimate]
\label{thm:2}
Let $1 \leq j \leq n$, let $\S'$ be a collection of multitiles in $\S$, and let 
$$\\P_{\S'}:=\{I\:\:\text{dyadic}:I_s\subseteq I\subseteq I_{s'}\:\:\text{for some}\:\: s,s'\in \S'\}$$ be the time convexification of $\S'$.  Then
$$\size_j(\S') \lesssim |E_j|^{-1/2} \sup_{I\in \\P_{\S'}}\frac{1}{|I|} \int_{E_j} \chi_{I}^{N}.$$
\end{proposition}

Note that this bound is consistent with the hypothesis $f_j \in X_2(E_j)$ and the intuition that the $j$-size is something like
a $\BMO$ average of $f_j$.

To decompose the collection of multitiles $\S$ into trees, we need the following result.

\begin{lemma}[Splitting lemma]
\label{l:18}
Let $\S'$ be a finite collection of multitiles, $1\le j\le n$ and  suppose that $\size_j(\S')\le 2^{m+1}$. 
Let $\mu > 0$ and suppose that $N$ is sufficiently large depending on $\mu$.
 Then  $\S'$ can be written as a  disjoint union 
 \begin{equation}\label{sp-decomp}
 \S'=(\bigcup_{\T \in \F} \T) \cup\S_2
\end{equation}
where $\F$ is a collection of trees
 such that
\begin{equation}
\label{eq:671}
 \sum_{\T\in\F}|I_T|\lesssim_\mu  2^{-2m}\left(\frac{1}{|E_j|^{1/2}2^m}\right)^{\frac{2}{\mu}},
\end{equation}
while  
\begin{equation}\label{2msize}
\size_j(\S_2)\le 2^{m}.
\end{equation}
\end{lemma}

This lemma is quite difficult and will be proven in Sections \ref{bessel-sec}-\ref{nonmax-sec}.  Assuming the lemma for the moment, 
we may iterate it in the standard way (see e.g. \cite{MTT2}) we conclude

\begin{corollary}[Tree selection algorithm]
\label{thm:6}
Let $\S'$ be a finite collection of multitiles and $1\le j\le n$. Let $\mu > 0$ and suppose that $N$ is sufficiently large
depending on $\mu$. 
Then, after discarding tiles $s$ of $j$-size zero (in the sense that $\langle f_j, \phi_{s,j} \rangle = 0$), 
there exists a partition 
$$\S'=\bigcup_{m: 2^m\le \size_j(\S')} \bigcup_{\T \in \F^{m,j}} \T$$
 where for each $m$, $\F^{m,j}$ is a collection of trees such that $\size_j(\T)\le 2^{m+1}$ and
\begin{equation}
\label{eq:3er}
\sum_{\T\in \F^{m,j}}|I_T|\lesssim_\mu 2^{-2m}{\left(\frac{1}{|E_j|^{1/2}2^m}\right)}^{\frac{2}{\mu}}.\end{equation}
\end{corollary}

Now we prove \eqref{tsum}.  It will suffice for each $l \geq 0$ to prove the stronger estimate
\begin{equation}\label{tsum-alt}
 \sum_{\T} |I_\T| \prod_{i=1}^{n}\size_i(\T) \lesssim 2^{-l} |E_1|^{\alpha_1} \ldots |E_n|^{\alpha_n}.
\end{equation}
under the additional assumption that
\begin{equation}
\label{eq:1297}
2^l\le 1+\frac{\text{dist}(I_s,\R\setminus\Omega)}{|I_s|}<2^{l+1}
\end{equation}
for all tiles $s\in \S$, since the original claim \eqref{tsum} then follows by dyadic decomposition of $\S$.

From \eqref{eq:1297} and Proposition \ref{thm:2} we have
\begin{equation}\label{size0}
\size_i(\S)\lesssim |E_i|^{\frac12}2^l \hbox{ for } 1 \leq i < n
\end{equation}
and
\begin{equation}\label{size1}
 \size_n(\S) \lesssim 2^{(1-N)l}.
 \end{equation}
Now use the selection algorithm in Theorem ~\ref{thm:6} for $\S$ to get for each $i$ the collections of trees $\F^{m,i}$; the tiles
of $i$-size zero can be safely discarded (viewing them as singleton trees)
as they make no contribution to \eqref{tsum-alt}.  We can then partition
$$ \S = \bigcup_{m_1,\ldots,m_n} \S^{m_1,1} \cap \ldots \cap \S^{m_n,n}$$
where $\S^{m,i} := \bigcup_{\T \in \F^{m,i}} \T$ and we implicitly assume that 
\begin{equation}\label{2mi}
2^{m_i}\le \size_i(\S).
\end{equation}
By pigeonholing we can restrict to the case when $m_j = \max(m_1,\ldots,m_n)$ for some fixed $1 \leq j \leq n$.  We
then have the partition
$$ \S = \bigcup_{m_1,\ldots,m_n: m_j = \max(m_1,\ldots,m_n)} \bigcup_{\T \in \F^{m_j,j}} (\T \cap \S^{m_1,1} \cap \ldots \cap \S^{m_n,n}).$$
Note that $\T \cap \S^{m_1,1} \cap \ldots \cap \S^{m_n,n}$ is a tree with the same top as $\T$, and with $j$-size at most $2^{m_j+1}$;
this tree need not contain its top, but this is of no consequence for us.
To verify \eqref{tsum-alt} it thus suffices to show that
\begin{equation}\label{mmm}
\sum_{m_1,\ldots,m_n: m_j = \max(m_1,\ldots,m_n)} \sum_{\T \in \F^{m_j,j}} |I_\T| 2^{m_1} \ldots 2^{m_n}
\lesssim 2^{-l} |E_1|^{\alpha_1} \ldots |E_n|^{\alpha_n}.
\end{equation}
Meanwhile, from \eqref{eq:3er} we have
$$ \sum_{\T\in \F^{m_j,j}}|I_T|\lesssim_\mu 2^{-2m_j}{\left(\frac{1}{|E_j|^{1/2}2^{m_j}}\right)}^{\frac{2}{\mu}}$$
where $\mu$ is a large parameter to be chosen later.  Also, from \eqref{2mi} we have
$$ 2^{m_1} \ldots 2^{m_n} \leq 2^{m_j} \prod_{i \neq j} \size_i(\S)^{2\alpha_i} 2^{(1-2\alpha_i)m_i}.$$
From these bounds and summing the geometric series in all the $m_i$ for $i \neq j$, we have
\begin{align*}
\sum_{m_1,\ldots,m_n: m_j = \max(m_1,\ldots,m_n)} &\sum_{\T \in \F^{m_j,j}} |I_\T| 2^{m_1} \ldots 2^{m_n}
\lesssim_\mu \\
&\prod_{i \neq j} \size_i(\S)^{2\alpha_i}
\sum_{m_j} 2^{m_j} (\prod_{i \neq j} 2^{(1-2\alpha_i) m_j}) 2^{-2m_j}{\left(\frac{1}{|E_j|^{1/2}2^{m_j}}\right)}^{\frac{2}{\mu}}.
\end{align*}
Since $\alpha_1 + \ldots + \alpha_n = (n-2)/2$, we can rewrite the right-hand side as
$$ \prod_{i \neq j} \size_i(\S)^{2\alpha_i}
\sum_{m_j} 2^{2\alpha_j m_j} {\left(\frac{1}{|E_j|^{1/2}2^{m_j}}\right)}^{\frac{2}{\mu}}.$$
Summing the geometric series, we can bound this (for $\mu$ sufficiently large) by
$$ \lesssim_\mu \prod_i \size_i(\S)^{2\alpha_i} {\left(\frac{1}{|E_j|^{1/2}\size_j(\S)}\right)}^{\frac{2}{\mu}}.$$
Applying \eqref{size0}, \eqref{size1} we obtain \eqref{mmm} as desired, if $N$ and $\mu$ are chosen sufficiently large.
 This proves Theorem \ref{reduct-5} and hence Theorem \ref{thm:ct1}.
 
It remains to prove Theorem \ref{thm:2} and Lemma \ref{l:18}.  This will occupy the remainder of the paper.
  
\section{Single tree size estimate}
\label{sec:6}

In this section we prove Theorem \ref{thm:2}.
This estimate is well known in the case $j\le n-1$, when the cutoff $|I_s| > 2^{k(x)}$ has no effect; see \cite[Lemma 6.8]{MTT2}.  Thus we shall focus instead on the more difficult case $j=n$.

It suffices to show that for each $f\in X(E)$, each $N$ and each  $n$-separated multitile tree $(\T,T,i)$
satisfying
$$ \sum_{s\in \T}|\langle f, \phi_{s,n}\rangle|^2
=\alpha^2 |I_\T|,$$
$$|\langle f, \phi_{s,n}\rangle|\le \alpha |I_s|^{1/2}$$
for each $s\in\T$, we have 

$$
\alpha
\lesssim \beta:=\sup_{I\in \\P_{\T}}\frac{1}{|I|} \int_{E} \chi_{I}^{N}.
$$

Fix such a $(\T,T,i)$, $N$ and $f$.  By frequency translation invariance we may assume that $0 \in \omega_{T,n}$. Assume for contradiction that the above inequality does not hold, more precisely, assume
\begin{equation}
\label{newlyaddedeq2}
\alpha\ge K \beta
\end{equation}
for some sufficiently large $K$, whose (implicit) value will become clear after a few lines of argument. 

We first note that if $f_1$ denotes the restriction of $f$ to the complement of $2I_\T$ then from the decay of $\phi_{s,n}$ we get 
$$|\langle f_1,\phi_{s,n}\rangle|\lesssim\left(\frac{|I_s|}{|I_\T|}\right)^{N} |I_\T|^{-\frac{1}{2}}\int_{E} \chi_{I_s}^{N}$$
for all $s \in \T$. So, if $f_2$ denotes the restriction of $f$ to  $2I_\T$ we conclude that if $K$ is large enough
$$ \sum_{s\in \T}|\langle f_2, \phi_{s,n}\rangle|^2
\ge \frac{\alpha^2}{2} |I_\T|,$$
$$|\langle f_2, \phi_{s,n}\rangle|\le 2\alpha |I_s|^{1/2}$$
for each $s\in\T$

We will next prove that
$$\frac{\alpha^2}{2}\le\frac{1}{|I_\T|} \int f_2 \sum_{s\in\T}\langle f_2,\phi_{s,n}\rangle\phi_{s,n} \lesssim \alpha\beta,$$
thus contradicting \eqref{newlyaddedeq2}. 

Denote $a_s:=\langle f_2,\phi_{s,n}\rangle$.  We can estimate
$$ |\sum_{s\in\T}a_s\phi_{s,n}(x)| \leq \sup_{k}|\sum_{s\in\T\atop{|I_s| > 2^k}}a_s\psi_{s,n}(x)|.$$
Since $\T$ is $n$-separated, we see from Lemma \ref{quasi}
that the tiles $s \in \T$ with $|I_s| > 2^k$ have a disjoint frequency support
from the tiles $s \in \T$ with $|I_s| \leq 2^k$.  Next, note that we can write $\sum_{s\in\T\atop{|I_s| > 2^k}}a_s\psi_{s,n}(x)$
as a Fourier multiplier applied to the function $F := \sum_{s\in\T}a_s\psi_{s,n}(x)$, where the symbol of the multiplier is a cutoff smoothly adapted
to an interval of length $\sim C_0 2^{-k}$.  From this and standard kernel estimates, we conclude that
$$ \sup_{k}|\sum_{s\in\T\atop{|I_s| > 2^k}}a_s\psi_{s,n}(x)| \lesssim \M(F),$$
and so it will suffice to show that
\begin{equation}\label{fmf}
\frac{1}{|I_\T|} \int f_2\M(F)\lesssim \alpha\beta.
\end{equation}

For a dyadic interval $J$ denote by $J_1,J_2,J_3$ the three dyadic intervals of the same length with $J$, sitting at the left of $J$, with $J_3$ being adjacent to $J$. Similarly  let $J_5,J_6,J_7$ be the three dyadic intervals of the same length with $J$, sitting at the right of $J$, with $J_5$ being adjacent to $J$. Also define $J_4:=J$. Let $\mathcal J$  be the set of all dyadic intervals $J$ with the following properties:
\begin{enumerate}
\item[(a)]$ J\cap 2I_T\not=\emptyset$
\item[(b)]$ \nexists \;I\in\\P_{\T}: |I|<|J|\:\:and\:\: I\subset 3J$
\item[(c)]$ J_i\in \\P_{\T}\:\: for\:\: some\:\:1\le i\le 7.$ 
\end{enumerate}

We claim that $2I_\T\subset\cup_{J\in\mathcal J}J$. Indeed, assume by contradiction that there exists some $x\in 2I_\T\setminus\cup_{J\in\mathcal J}J$. Let $J^{(0)}\subset J^{(1)}\subset J^{(2)}\subset\ldots$ be the sequence of dyadic intervals of consecutive lengths containing $x$, with $|J^{(0)}|=\min_{I\in\\P_{\T}}|I|.$ Since $J^{(0)}\notin\mathcal J$ and since (a) and (b) are certainly satisfied for  $J^{(0)}$, it follows that $J_i^{(0)}\notin \\P_{\T}$ for each $1\le i\le 7$. Moreover, note that for each $1\le i\le 7$ there is no $I\in\\P_{\T}$ with $I\subset J_i^{(0)}$. We proceed now by induction. Assume that for some $j\ge 0$ we proved that for each $1\le i\le 7$ we have  $J_i^{(j)}\notin \\P_{\T}$  and also that there is no $I\in\\P_{\T}$ with $I\subset J_i^{(j)}$. Note that this implies the same for $j+1$. Indeed, since $3J^{(j+1)}\subset 7J^{(j)}$ and by induction hypothesis, it follows that (b) is satisfied for $J^{(j+1)}$. Hence $J_i^{(j+1)}\notin\\P_{\T}$ for each $1\le i\le 7$. We verify now  the second statement of the induction. Note that  if there was an $I\in\\P_{\T}$ with $I\subset J_i^{(j+1)}$ than the hypothesis of the induction and the fact that $3J^{(j+1)}\subset 7J^{(j)}$ would imply that $i\in\{1,2,6,7\}$. Hence $I\subset J_i^{(j+1)}\subset I_T$, and by convexity of $\\P_{\T}$ it would follow that $J_i^{(j+1)}\in\\P_{\T}$, impossible. This closes the induction. To see how the claim follows from here, observe that $I_T=J_i^{(j)}$ for some $i,j$, which certainly contradicts the fact that $I_T\in\\P_{\T}$. 

Next thing we prove is that on each interval $2J$ with $J\in \mathcal J$, the function $F$ is ``essentially'' constant. More exactly we will show that for each $x\in 2J$, $|F(x)-F(c(J))|\lesssim \alpha.$ We have
\begin{equation*}
|F(x)-F(c(J))|\le \sum_{s\in\T\atop{|I_s|\ge |J|}}|a_s|\sup_{x\in 2J}|\psi_{s,n}(x)-\psi_{s,n}(c(J))|+2\sum_{s\in\T\atop{|I_s|<|J|,I_s\cap 3J=\emptyset}}|a_s|\sup_{x\in 2J}|\psi_{s,n}(x)|
\end{equation*}

The first term on the right hand side can be further bounded by
\begin{equation*}
|J|\sum_{s\in\T\atop{|I_s|\ge |J|}}\sup_{x\in 2J}|\psi_{s,n}^{\prime}(x)||a_s|\lesssim |J|\alpha\sum_{s\in\T\atop{|I_s|\ge |J|}}\frac{1}{|I_s|}\chi_{I_s}^{3}(c(J))\lesssim \alpha.
\end{equation*}
while the second term can be bounded by
\begin{equation*}
\sum_{s\in\T\atop{|I_s|<|J|,I_s\cap 3J=\emptyset}}\frac{1}{|I_s|^{\frac12}}\chi_{I_s}^{3}(c(J))|a_s|\lesssim \alpha\sum_{s\in\T\atop{|I_s|<|J|,I_s\cap 3J=\emptyset}}\frac{1}{|I_s|^{\frac12}}\chi_{I_s}^{3}(c(J))\lesssim \alpha
\end{equation*}

Define now the measure space $X=\cup_{J\in\mathcal J}J$ and its $\sigma$-algebra $\Upsilon$ generated by the maximal intervals $ J\in\mathcal J$. Recall that $2I_T\subset\cup_{J\in\mathcal J}J=X\subset 10I_T$. We will see that for each $x\in J$ 
\begin{equation}
\label{eq:50}
\M(F)(x)\lesssim \frac1{|J|}\int_{J} \M(F)(z)dz+\alpha.
\end{equation}
Indeed, if $r>\frac12|J|$,  
\begin{align*}
\frac{1}{2r}\int_{x-r}^{x+r}|F|(z)dz&\lesssim \inf_{y\in J}\M(F)(y)\\&\lesssim \frac1{|J|}\int_{J} \M(F)(z)dz.
\end{align*}
 On the other hand, if $r\le\frac12|J|$,
\begin{align*}
\frac{1}{2r}\int_{x-r}^{x+r}|F|(z)dz&\lesssim \sup_{y\in 2J}|F|(y)\\&\lesssim \inf_{y\in J}|F|(y)+\alpha\\&\lesssim \frac1{|J|}\int_{J} \M(F)(z)dz +\alpha.
\end{align*}
From ~\eqref{eq:50} we can write
\begin{align*}
\frac{1}{|I_\T|}\int f_2\M(F)&\lesssim \frac{1}{|I_\T|}\int_{X} f_2\E(\M(F)|\Upsilon)+\alpha\sup_{J\in \mathcal J}\frac {1}{|J|}\int_{J}f_2
\\&=\frac{1}{|I_\T|}\int_{X} \E(f_2|\Upsilon)\E(\M(F)|\Upsilon)+\alpha\sup_{J\in \mathcal J}\frac1{|J|}\int_{J}f_2\\&\lesssim \frac{1}{|I_\T|}\|\E(f_2|\Upsilon)\|_{L^{\infty}}\int_{X} \E(\M(F)|\Upsilon)+\alpha
\sup_{J\in  \mathcal J}\frac{1}{|J|}\int_{E} \chi_{J}^N\\&\lesssim \left(\alpha+\left[\frac{1}{|I_\T|}\int_{X} \E(\M(F)|\Upsilon)^2\right]^{\frac12}\right)\sup_{J\in  \mathcal J}\frac{1}{|J|}\int_{E} \chi_{J}^N\\&\lesssim \alpha\sup_{J\in  \mathcal J}\frac{1}{|J|}\int_{E} \chi_{J}^N,
\end{align*}
where $\E(\cdot|\Upsilon)$ denotes the conditional expectation relative to $\Upsilon.$ Finally, note that since for each $J\in\mathcal J$, $J_i\in\\P_{\T}$ for some $i$, we have that 
$$\sup_{J\in  \mathcal J}\frac{1}{|J|}\int_{E} \chi_{J}^N\lesssim \sup_{I\in \\P_{\T}}\frac{1}{|I|}\int_{E} \chi_{I}^N,$$ 
which yields \eqref{fmf}.  This concludes the proof of Theorem \ref{thm:2}.
\endprf

\section{Reduction to Bessel inequality}\label{bessel-sec}

We still have to prove Lemma \ref{l:18}.  This will be achieved by means of a certain maximal Bessel inequality and
a stopping time argument. We first recall a definition.

\begin{definition}\cite{MTT1}
Let $j\in\{1,2,\ldots,n\}$. Two $j$-separated multitile trees $(\T,T,i)$ and $(\T',T',i)$ with the same index are said to be 
\emph{strongly $j$-disjoint} if $\T \cap \T' = \emptyset$, and furthermore
whenever $s\in \T, s'\in \T'$ are such that $\omega_{s_j}\subsetneq \omega_{s'_j}$, then one has $I_T\cap I_{s'}=\emptyset$, and similarly with $\T$ and $\T'$ reversed.
A collection of $j$-separated multitile trees is called \emph{mutually strongly $j$-disjoint} if each two multitile trees in the collection are strongly $j$-disjoint.
\end{definition}

\begin{remark} If two $j$-separated multitile trees $(\T,T,i)$ and $(\T',T',i)$ are strongly $j$-disjoint, then one
has $s_j\cap s_j'=\emptyset$ for each $s\in \T, s'\in \T'$.  This is because if $s_j$ and $s'_j$ intersect, then since $\T \cap \T' = \emptyset$, we must have either $\omega_{s_j} \subsetneq \omega_{s'_j}$ or $\omega_{s'_j} \subsetneq \omega_{s_j}$, and the claim then follows from the definition of strong $j$-disjointness.  This may help explain the terminology ``strong $j$-disjointness''.
\end{remark}

The next estimate controls the extent to which disjoint trees can each absorb a lot of $L^2$ energy.  It is the main 
technical estimate used in the proof, and the core of Lacey's original argument in \cite{La}.   The proof is rather difficult 
and will occupy the remainder of this paper.

\begin{theorem}[Maximal Bessel inequality, multitile version]
\label{thm:5}
Let $\F$ be a finite collection of strongly $j$-disjoint, $j$-separated multitile trees. 
Let $\mu > 0$ and suppose that $N$ is sufficiently large depending on $\mu$ (recall the definition of $N$ from \eqref{e:fieq}).
Assume also that 
\begin{equation}
\label{eq:70}
2^m\le\left(\frac{1}{|I_\T|}\sum_{s\in \T}|\langle f_j,\phi_{s,j}\rangle|^2\right)^{\frac12}\le 2^{m+1}
\end{equation}
and
\begin{equation}
\label{eq:71}
\left(\frac{1}{|I_{T'}|}\sum_{s\in \T\atop{I_s\subset I_{T'}}}|\langle f_j,\phi_{s,j}\rangle|^2\right)^{\frac12}\le 2^{m+1}
\end{equation}
 for each $T'\in \T\in\F$.
Then if $\mu >0$  we have
\begin{equation}
\label{eq:12}
\sum_{\T \in \F} |I_{\T}| \lesssim_\mu 2^{-2m}{\left(\frac{1}{|E_j|^{1/2}2^m}\right)}^{\frac{2}{\mu}}.\end{equation}
\end{theorem}

\begin{remark}  The factor ${\left(\frac{1}{|E_j|^{1/2}2^m}\right)}^{\frac{2}{\mu}}$ in \eqref{eq:12} is technical and should
be ignored.  Intuitively, the condition \eqref{eq:70} asserts that the function $f_j$, when ``restricted'' to a tree $\T$ in $\F$, has
$L^2$ norm roughly comparable to $2^m |I_T|^{1/2}$.  The strong disjointness of the trees is an assertion that these restrictions are in 
some sense ``almost orthogonal''.  Since $f_j$ has an $L^2$ norm of $O(1)$, we see that \eqref{eq:12} is indeed a kind of Bessel inequality.
This estimate is standard (and fairly straightforward) when $j \neq n$, but when $j=n$ the presence of the cutoff $1_{|I_s| > 2^{k(x)}}$ in the
modified wave packet $\phi_{s,n}$ presents some significant difficulties (already encountered in \cite{La}).
\end{remark}

Let us now show how Lemma \ref{l:18} follows from Theorem \ref{thm:5}.  This will be a standard stopping time argument of the type which has been commonly used in time-frequency analysis, see for instance  \cite{MTT1}, but for sake of completeness we present
the argument here.

We perform the following algorithm to construct $\S_2$ and $\F$.

\begin{itemize}

\item Step 0.  Initialize $\F$ to be empty, and $\S_2$ to equal $\S'$.

\item Step 1.  If $\size_j(\S_2) \leq 2^m$, then we terminate the algorithm.  Otherwise, we have
$$ 2^m < \size_j(\S_2) \leq \size_j(\S') \leq 2^{m+1}.$$
By definition of size, we can find a $j$-separated multitile tree $\T = (\T,T,i)$ in $\S_2$ obeying \eqref{eq:71}.

\item Step 2.  The multitile tree $(\T,T,i)$ mentioned above is a $(j,\epsilon)$-separated tree for some $\epsilon = \pm 1$.  For fixed $i$ and $\epsilon$, we may assume that this tree maximizes the quantity $\epsilon \xi_{T_i}$, where $\xi_{T_i}$ is the center of the frequency tile $T_i$, subject to the constraints \eqref{eq:71} and $\T \subseteq \S_2$.

\item Step 3.  Clearly the multitile tree $\T$ is non-empty, since it has positive size.
Add the multitile tree $\T$ to the collection $\F$, and delete the multitiles in $\T$ from $\S_2$.  Note that this removes at least one multitile from $\S_2$.

\item Step 4.  Next, define the (possibly empty) companion 
tree $(\tilde \T, T, j)$ where $\tilde \T := \{ s \in \S_2: s_j \leq T_j \}$, add this tree $\tilde \T$ to $\F$ also, 
and delete the multitiles in $\tilde \T$ from $\S_2$. Then return to Step 1.

\end{itemize}

This algorithm terminates in finite time since $\S_2$ was initially finite, and every iteration of the algorithm removes at least
one multitile from $\S_2$.  It is also clear that this algorithm will obtain a decomposition \eqref{sp-decomp} obeying \eqref{2msize}.  The only remaining task is to verify the bound \eqref{eq:671}.  It suffices to do this for each fixed $1 \leq i \leq n$, thus restricting the summation to those trees $\T = (\T,T,i)$ in $\F$ with index $i$.  This in turn fixes the quantity
$\epsilon$ appearing in Step 2 above, so if one indexes the trees $\T$ in the order that they are added to $\F$ then $\epsilon \xi_{T_j}$ will be non-increasing.  We also only need to focus on those trees selected using Step 3 rather than Step 4, since
the trees in Step 4 have the same time interval as those in Step 3 and so we are only giving up a factor of $2$ by doing this.

To prove \eqref{eq:671}, it suffices by Theorem \ref{thm:5} to show that the trees in $\F$ with fixed $i$ and $\epsilon$ arising
from Step 3 are mutually strongly $j$-disjoint.  Suppose for contradiction that there were two trees $\T, \T'$ in $\F$ of this type
which were not strongly $j$-disjoint.  Since these trees have distinct multitiles by construction, the only way that strong $j$-disjointness can fail (up to swapping $\T$ and $\T'$) is if there exist $s \in \T$, $s' \in \T'$ with
$\omega_{s,j}\subsetneq \omega_{s',j}$ and $I_T\cap I_{s'} \neq \emptyset$.  From Lemma \ref{quasi} we conclude that
$s'_j \leq T_j$ and $\epsilon (\xi_{T_i} - \xi_{T'_i}) \geq 0$.  The latter condition, combined with the non-increasing
nature of the $\epsilon \xi_{T_j}$, ensures that $T$ was selected earlier in the algorithm than $T'$.  But then $s'$ would have been selected in the companion tree $\tilde T$ and could not have remained in $\S_2$ by the time $T'$ was selected, a contradiction.
This ensures the strong $j$-disjointness and concludes the deduction of Lemma \ref{l:18} from Theorem \ref{thm:5}.

It remains to prove Theorem \ref{thm:5}.  This will occupy the remainder of the paper.

\section{Good-$\lambda$ reduction}\label{bmo-sec}

The only remaining task in the proof of Theorem \ref{thm:ct1} is the maximal Bessel inequality in Theorem \ref{thm:5}.  
This will be accomplished in stages.  In this section we rephrase the inequality as an inequality concerning tiles rather than
multitiles, and use some ``$\BMO$ theory'' for tiles to replace the ${\left(\frac{|E_j|^{1/2}}{2^m}\right)}^{\frac{2}{\mu}}$
factor in \eqref{eq:12} by a factor which depends instead on the counting function $N_\F$.  This $\BMO$ theory is quite elementary
and may have some independent interest.

We will focus on the hardest case $j=n$, in which one must deal with the presence of the cutoff $1_{|I_s| > 2^{k(x)}}$ in the
modified wave packet $\phi_{s,n}$.  The cases $j \neq n$ are significantly simpler (see for instance \cite{MTT1}) and in any event
can be handled by the argument here (e.g. by the artificial expedient of setting $k(x)$ to be so low that the cutoff $1_{|I_s| > 2^{k(x)}}$ disappears).

The Bessel inequality is now really only a statement about the $n$-tiles of the multitiles in $\S$, and so we shall introduce new notation to focus only on these tiles rather than on the multitiles.

\begin{definition}[Lacunary tree]  A \emph{lacunary tree} $\T = (\T,I_\T,\xi_\T)$ is a collection $\T$ of tiles, together with a dyadic time interval
$I_\T \in \D_0$ and a center 
frequency $\xi_\T \in \R$, such that for all $P \in \T$ we have $I_P \subseteq I_\T$ and $\dist(\omega_P, \xi_\T) \sim C_0 |\omega_P|$,
and such that the frequency interval $\omega_P$ of a tile is determined entirely by the length of the time interval, thus if $P, P' \in \T$
and $|I_P| = |I_{P'}|$ then $\omega_P = \omega_{P'}$.  (In particular, this means that distinct tiles in $\T$ have distinct time intervals.)
We say that one lacunary tree $(\T',I_{\T'},\xi_\T)$ is a \emph{subtree} of another $(\T,I_{\T'},\xi_{\T'})$ if $\T' \subseteq \T$ (thus we allow subtrees to have a different time interval and center than the supertree).
We say that two lacunary trees $(\T,I_\T,\xi_\T)$, $(\T',I_{\T'},\xi_{\T'})$ are \emph{strongly disjoint} if $\T \cap \T' = \emptyset$, and whenever $P \in \T$, $P' \in \T'$ are such that $\omega_P \subsetneq \omega_{P'}$, then one has $I_\T \cap I_{P'}=\emptyset$, and similarly with $\T$ and $\T'$ reversed.  We define a \emph{forest} to be any collection $\F$ of lacunary trees such that
any two distinct trees $\T, \T'$ in $\F$ are strongly disjoint.
\end{definition}

Observe from Lemma \ref{quasi} that if $(\T,T,i)$ is an $n$-separated multitile tree, then $(\T_n,I_{T_n},\xi_{T_n,n})$ is a lacunary tree, 
where $\T_n := \{ s_n: s \in \T \}$ is the set of $n$-tiles of the multitile tree $\T$, and $\xi_{T_n,n}$ is a frequency such that
$\dist( \omega_{T_n,n}, \xi_{T_n,n} ) \sim C_0 |\omega_{T_n,n}|$.  Furthermore, if $(\T,T,i)$ and $(\T',T',i)$ are strongly $n$-disjoint, then $(\T_n, I_{T_n}, \xi_{T_n,n})$ and $(\T'_n, T'_n, \xi_{T'_n,n})$ are strongly disjoint.  Thus, we can deduce Theorem \ref{thm:5} from

\begin{theorem}[Maximal Bessel inequality, first reduction]
\label{thm:5a}
Let $\F$ be a forest.  
Let $\mu > 0$ and suppose that $N$ is sufficiently large depending on $\mu$.
For each tile $P$ in $\bigcup_{\T \in \F} \T$, let
$\psi_P$ be a wave packet adapted to $P$, and let $\phi_P$ be the function
$$ \phi_P(x) := 1_{|I_s| > 2^{k(x)}} \psi_P(x).$$
Let $E$ be a finite union of intervals, and let $f \in X_2(E)$ be such that
$$
2^m\le\left(\frac{1}{|I_\T|}\sum_{P \in \T}|\langle f,\phi_P \rangle|^2\right)^{\frac12}\le 2^{m+1}
$$
and
$$
\left(\frac{1}{|I_{T'}|}\sum_{P\in \T\atop{I_P\subset I_{T'}}}|\langle f,\phi_P\rangle|^2\right)^{\frac12}\le 2^{m+1}
$$
 for each $T'\in \T\in\F$.
Then we have
$$
\sum_{\T \in \F} |I_{\T}| \lesssim_\mu 2^{-2m}{\left(\frac{1}{|E|^{1/2}2^m}\right)}^{\frac{2}{\mu}},
$$
for all $\mu > 0$.
\end{theorem}

We will now eliminate the role of the set $E$, replacing it with a certain counting function multiplicity, and also eliminate the role of the size parameter $2^m$.  More precisely, in this section we shall deduce Theorem \ref{thm:5a} from 

\begin{theorem}[Maximal Bessel inequality, second reduction]
\label{thm:3} 
Let $\F$ be a forest.
Let $\mu > 0$ and suppose that $N$ is sufficiently large depending on $\mu$.
Let $\psi_P$, $\phi_P$ be as in Theorem \ref{thm:5a}.
Let $f \in L^2(\R)$ be such that
\begin{equation}
\label{eq:17894hhnt}
1 \leq \left(\frac{1}{|I_\T|}\sum_{P\in \T}|\langle f,\phi_P\rangle|^2\right)^{\frac12} \leq 2
\end{equation}
for each $\T\in \F$, and 
\begin{equation}
\label{eq:17895}
\left(\frac{1}{|I_{T'}|}\sum_{P\in \T\atop{I_P\subset I_{T'}}}|\langle f,\phi_P\rangle|^2\right)^{\frac12} \leq 2
\end{equation}
 for each $T'\in \T$. Let $N_\F$ be the counting function
\begin{equation}\label{nf-def}
N_\F := \sum_{\T \in \F} 1_{I_\T},
\end{equation} 
and let $I_0$ be any interval which contains the support of $N_F$.
Then we have the Bessel-type inequality
\begin{equation}
\label{eq:drac}
\sum_{P \in \bigcup_{\T \in \F} \T}|\langle f,\phi_P\rangle|^2 \lesssim_\mu \|N_{\F}\|_{L^\infty}^{\frac{1}{\mu}} \int |f|^2 \chi_{I_0}^{10}.
\end{equation}
\end{theorem}

We shall prove Theorem \ref{thm:3} in later sections.  For now, we show how it implies Theorem \ref{thm:5}.  The argument is similar to the ``good-$\lambda$'' type estimates used to prove John-Nirenberg $\BMO$ inequalities, and to emphasize this connection (and because this theory may be of some independent interest) we shall proceed in a somewhat abstract manner.

The following observation is trivial, but is still worth recording.

\begin{lemma}[Forest refinement]\label{forestry} 
Let $\F$ be a forest.  For each tree $\T$ in $\F$, let $\F_\T$ be a collection of subtrees
of $\T$ with disjoint time intervals.  Then $\bigcup_{\T \in \F} \F_\T$ is also a forest.
\end{lemma}

Let $\F$ be any forest.  The quantity $\|N_{\F}\|_{L^\infty}$ measures the maximum possible overlap of the
time intervals $I_\T$ of the trees $\T$ in $\F$.  We shall introduce a closely related quantity $\|\F\|_{\BMO}$, defined as
$$ \|\F \|_{\BMO} := \sup_I \frac{1}{|I|} \sum_{\T \in \F: I_\T \subseteq I} |I_\T|,$$ where the supremum is taken over all the 
dyadic intervals $I$.  

\begin{remark} One can relate this $\BMO$-type norm to the genuine (dyadic, vector-valued) $\BMO$ norm by the formula $\|\F\|_\BMO = \| \vec N_\F \|_\BMO^2$, where $\vec N_\F := \sum_{\T \in \F} 1_{I_\T} e_\T$ is a vector-valued counting function, with the $e_\T$ being orthonormal vectors in an abstract Hilbert space.  However, we will not adopt this approach since the theory of vector-valued $\BMO$ is not as familiar as that of ordinary $\BMO$, preferring instead a more direct and elementary approach.
\end{remark}

It is clear that $\|\F\|_{\BMO} \leq \|N_{\F}\|_{L^\infty}$, indeed
\begin{align*}
\frac{1}{|I|} \sum_{\T \in \F: I_\T \subseteq I} |I_\T|
&= \frac{1}{|I|} \int_I \sum_{\T \in \F: I_\T \subseteq I} 1_{I_\T}  \\
&\leq \| \sum_{\T \in \F: I_T \subseteq I} 1_{I_\T} \|_{L^\infty}\\
&\leq \| N_{\F}\|_{L^\infty}.
\end{align*}
While the converse is not quite true, we do expect the $L^\infty$ norm and $\BMO$ norm  to be very close.

Now we obtain some good lambda inequalities for the $\BMO$ norm.  We first observe that to control the $\BMO$ norm of
a collection $\F$ of trees, it suffices to control the $\BMO$ norm of subcollections of trees which are already controlled
in $L^\infty$.

\begin{lemma}  
Let $\F$ be a forest such that
$$ \| \F' \|_{\BMO} \leq B $$
whenever $\F' \subset \F$ is such that $\|N_{\F'} \|_{L^\infty} \leq 2B$.  Then we have
$$ \| \F \|_{\BMO} \leq 2B.$$
\end{lemma}

\begin{proof}  Let $I_0$ be a dyadic interval.  
Call a dyadic interval $J \subseteq I_0$ \emph{heavy} if $|\{ \T \in \F: J \subseteq I_\T \subseteq I_0 \}| > 2B$, and let
$\F'$ be the collection of those trees $\T \in \F$ such that $I_\T \subseteq I_0$ and that $I_\T$ is not heavy.  Then by construction we have
$\|\F'\|_{L^\infty} \leq 2B$, and hence by hypothesis $\|\F'\|_{\BMO} \leq B$.  In particular
$$ \sum_{\T \in \F: I_\T \subseteq I_0; I_\T \hbox{ not heavy}} |I_\T| \leq B |I_0|.$$

Now we deal with the heavy intervals.  If we let $\J$ denote the set of maximal dyadic heavy intervals, then we have
\begin{align*}
\sum_{\T \in \F: I_\T \subseteq I_0; I_\T \hbox{ heavy}} |I_\T|
&\leq \sum_{J \in \J} \sum_{\T \in \F: I_\T \subseteq J} |I_\T| \\
&\leq \sum_{J \in \J} \| \F \|_{\BMO} |J| \\
&= \frac{\| \F \|_{\BMO}}{2B} \int_{\bigcup_{J \in \J} J} 2B \\
&\leq \frac{\| \F \|_{\BMO}}{2B} \int_{\bigcup_{J \in \J} J} \sum_{\T \in \F'} 1_{I_{\T}} \\
&\leq \frac{\| \F\|_{\BMO}}{2B} \int_{I_0} \sum_{\T \in \F'} 1_{I_{\T}} \\
&\leq \frac{\| \F \|_{\BMO}}{2B} |I_0| \| \F' \|_{\BMO} \\
&\leq \frac{\| \F \|_{\BMO}}{2} |I_0|. 
\end{align*}
Summing these two estimates, we obtain
$$ \frac{1}{|I_0|} \sum_{\T \in \F: I_\T \subseteq I_0} |I_\T| \leq B + \frac{\|\F\|_{\BMO}}{2},$$
and then taking supremum over $I_0$ we obtain
$$ \| \F \|_{\BMO} \leq B + \frac{\|\F\|_{\BMO}}{2}.$$
The claim follows.
\end{proof}

Similarly, to control the $L^1$ norm of $N_{\F}$, it suffices to control the $L^1$ norm of subcollections $\F'$ which are controlled in $L^\infty$ by the $\BMO$ norm of $\F$:

\begin{lemma}
Let $\F$ be a forest such that 
$$ \| N_{\F'} \|_{L^1} \leq A$$
whenever $\F' \subset \F$ is such that $\| N_{\F'} \|_{L^\infty} \leq \|\F\|_{\BMO}$.  Then we have
$$ \| N_{\F} \|_1 \leq 2A.$$
\end{lemma}

\begin{proof} Set $B := \|\F\|_{\BMO}$.
As before, we call a dyadic interval $J \subseteq I_0$ \emph{heavy} if $|\{ \T \in \F: J \subseteq I_\T \subseteq I_0 \}| > B$, and let
$\F'$ be the collection of those trees $\T \in \F$ such that $I_\T \subseteq I_0$ and that $I_\T$ is not heavy.  Then by construction we have
$\|N_{\F'}\|_{L^\infty} \leq B$ and hence by hypothesis $\|N_{\F'}\|_1 \leq A$.  Now if we let $\J$ be the collection of maximal heavy intervals, then we have
\begin{align*}
 \| N_{\F \backslash \F'} \|_1 &=  \sum_{J \in \J} \sum_{\T \in \F: I_\T \subseteq J} |I_\T| \\
 &\leq \sum_{J \in \J} \|\F \|_{\BMO} |J| \\
 &= \int_{\bigcup_{J \in \J} J} B \\
 &\leq \int_{\bigcup_{J \in \J} J} \sum_{\T \in \F'} 1_{I_\T} \\
 &\leq \| N_{\F'} \|_1 \\
 &\leq A
\end{align*}
and the claim follows.
\end{proof}

We can of course combine these two lemmas to obtain

\begin{corollary}
Let $\F$ be a forest such that 
$$ \| N_{\F'} \|_1 \leq A \hbox{ and } \| \F' \|_{\BMO} \leq B$$
whenever $\F' \subset \F$ is such that $\| N_{\F'} \|_{L^\infty} \leq 2B$.  Then we have
$$ \| N_{\F} \|_1 \leq 2A \hbox{ and } \| \F \|_{\BMO} \leq 2B.$$
\end{corollary}

A specific case of this is

\begin{corollary}
\label{cor:7}
  Let $\F$ be a forest such that for some $\mu > 1$
$$ \| N_{\F'}\|_1 \leq A  \| N_{\F'}\|_{L^\infty}^{\frac1{\mu}} \hbox{ and } \| \F'\|_{\BMO} \leq B  \|N_{\F'}\|_{L^\infty}^{\frac{1}{\mu}} $$
for all $\F' \subseteq \F$.  Then  we have
$$ \| N_{\F} \|_1 \lesssim_\mu A  B^{\frac{1}{\mu-1}} \hbox{ and }
\| \F \|_{\BMO} \lesssim_\mu B  B^{\frac{1}{\mu-1}}.$$
\end{corollary}

Thus to prove a counting function estimate on $\|N_{\F}\|_1$, we are permitted to lose a small  power of the $\|N_{\F}\|_{L^\infty}$ as long as the argument also works for all subtrees and localizes to a $\BMO$ version as well (with a different constant $B$).

Now we can finally prove Theorem \ref{thm:5a}.

\begin{proof}[of Theorem \ref{thm:5a}]  Let $\F' \subset \F$ be arbitrary.  From Theorem \ref{thm:3} with $f$ replaced by $f/2^m$, and $I_0$ chosen to be so large as to contain all the time intervals arising from $\F'$, we have
\begin{align*}
\| N_{\F'}\|_1 &= \sum_{\T \in \F'} |I_{\T}| \\
&\lesssim
\sum_{P \in \bigcup_{\T \in \F'} \T}|\langle f/2^m,\phi_P\rangle|^2 \\
&\lesssim_\mu \|N_{\F'}\|_{L^\infty}^{\frac{1}{\mu}} \int |f/2^m|^2 \\
&\lesssim 2^{-2m} \|N_{\F'}\|_{L^\infty}^{\frac{1}{\mu}} 
\end{align*}
thanks to the $L^2$ normalization of $f \in X_2(E)$.  If we let $I_0$ be an arbitrary dyadic interval, then by
replacing $\F'$ by $\{ \T \in \F': I_\T \subseteq I_0 \}$ in the above argument we see that
\begin{align*}
\frac{1}{|I_0|} \sum_{\T \in \F': I_\T \subseteq I_0} |I_{\T}|
&\lesssim_\mu \frac{1}{|I_0|} \|N_{\F'}\|_{L^\infty}^{\frac{1}{\mu}} \int |f/2^m|^2 \chi_{I_0}^{10} \\
&\lesssim \|N_{\F'}\|_{L^\infty}^{\frac{1}{\mu}} 2^{-2m} |E|^{-1} 
\end{align*}
thanks to the uniform bound of $|E|^{-1/2}$ on $f \in X_2(E)$.  Taking suprema over $I_0$ we conclude
that $\| \F' \|_{\BMO} \lesssim_\mu \|N_{\F'}\|_{L^\infty}^{\frac{1}{\mu}} 2^{-2m} |E|^{-1}$.
Applying Corollary \ref{cor:7} we conclude that
$$ \sum_{\T \in \F} |I_\T| = \| N_{\F} \|_1 \lesssim_\mu 2^{-2m} (2^{-2m} |E|^{-1})^{\frac{1}{\mu-1}}.$$
Replacing $\mu$ by $\mu+1$ we obtain Theorem \ref{thm:5a}.
\end{proof}

\section{Tileset refinements}

It remains to prove Theorem \ref{thm:3}.  In this section we perform some additional elementary reductions. 
First we eliminate the localizing
weight $\chi_{I_0}^{10}$ and we permit the deletion of those tiles which lie inside a small exceptional set.  Then we sparsify the tile set, and remove some logarithmic pileups of time interval multiplicity.

We begin with the first reduction.  We assert that to prove Theorem \ref{thm:3} it suffices to prove the same assertion
with the weight $\chi_{I_0}^{10}$ not present in \eqref{eq:drac}.  The reason for this is that $\chi_{I_0}^{-10}$ is a polynomial,
and because of this (and the hypothesis that all the tiles have time interval contained in $I_0$) 
$\chi_{I_0}^{-10} \psi_P$ is a wave packet adapted to $P$, except for the trivial change that the exponent of $10N$ in \eqref{psisj}
must be reduced slightly to $10(N-1)$.  But this clearly makes no essential difference to the argument since we are free to take $N$
as large as we wish.  Since $\langle f, \phi_P \rangle = \langle f \chi_{I_0}^{10}, \chi_{I_0}^{-10} \phi_P \rangle$, we thus
see that Theorem \ref{thm:3} with the localizing weight $\chi_{I_0}^{10}$ follows more or less automatically from Theorem \ref{thm:3}
without the localizing weight.

The next step is to eliminate the hypotheses \eqref{eq:17894hhnt}, \eqref{eq:17895} and also give the ability to delete a small
exceptional collection of tiles.

\begin{theorem}[Maximal Bessel inequality, third reduction]
\label{thm:3a} 
Let $\F$ be a forest.
Let $\mu > 0$ and suppose that $N$ is sufficiently large depending on $\mu$.
Let $\psi_P$, $\phi_P$ be as in Theorem \ref{thm:5a}, and let $N_\F$ be the counting function \eqref{nf-def}.  
Then there exists an exceptional set $\P_* \subset \bigcup_{\T \in \F} \T$ of tiles with
\begin{equation}\label{omega-bound}
|\bigcup_{P \in \P_*} I_P| \leq \frac{1}{10} \frac{\| N_\F \|_{L^1}}{\| N_\F\|_{L^\infty}}
\end{equation}
such that we have the Bessel-type inequality
\begin{equation}
\label{eq:drac-2}
\sum_{P \in \bigcup_{\T \in \F} \T \backslash \P_*}|\langle f,\phi_P\rangle|^2 \lesssim_\mu \|N_{\F}\|_{L^\infty}^{\frac{1}{\mu}} \|f\|_{L^2}^2
\end{equation}
 for all $f \in L^2(\R)$.
\end{theorem}

\begin{proof}[of Theorem \ref{thm:3} assuming Theorem \ref{thm:3a}]  Write $\Omega := \bigcup_{P \in \P_*} I_P$.  Then
$$ \sum_{P \in \bigcup_{\T \in \F} \T: I_P \not \subseteq \Omega}|\langle f,\phi_P\rangle|^2\le \sum_{P \in \bigcup_{\T \in \F} \T \backslash \P_*}|\langle f,\phi_P\rangle|^2.$$ 
To prove \eqref{eq:drac}, it thus suffices in view of \eqref{eq:drac-2} to show that
$$ \sum_{P \in \bigcup_{\T \in \F} \T: I_P \subseteq \Omega}|\langle f,\phi_P\rangle|^2 \leq \frac{1}{2} 
\sum_{P \in\bigcup_{\T \in \F} \T}|\langle f,\phi_P\rangle|^2.$$
From \eqref{eq:17894hhnt}, it thus suffices to show that
$$ \sum_{P \in \bigcup_{\T \in \F} \T: I_P \subseteq \Omega}|\langle f,\phi_P\rangle|^2 \leq \frac{1}{2} \| N_\F \|_{L^1}.$$
For each tree $\T$ in $\F$,
consider the tile set $\{ P \in \T: I_P \subseteq \Omega \}$.  If $Q$ is any tile in this set with $I_Q$ maximal with respect to set inclusion, then $I_Q \subseteq \Omega$ and from \eqref{eq:17895} we have
$$ \sum_{P \in \T: I_P \subseteq I_Q \subseteq \Omega}|\langle f,\phi_P\rangle|^2 \leq 4 |I_Q|.$$
Summing this over all such $Q$ (noting that the $I_Q$ are disjoint by dyadicity and maximality) we conclude
$$ \sum_{P \in \T: I_P \subseteq \Omega}|\langle f,\phi_P\rangle|^2 \leq 4 |I_\T \cap \Omega| = 4 \int_\Omega 1_{I_\T}.$$
Summing this over all $\T \in \F$ we obtain
$$ \sum_{P \in \bigcup_{\T \in \F} \T: I_P \subseteq \Omega}|\langle f,\phi_P\rangle|^2 \leq 4 \int_\Omega N_\F
\leq 4 |\Omega| \| N_\F \|_{L^\infty}$$
and the claim follows from \eqref{omega-bound}.
\end{proof}

We still have to prove Theorem \ref{thm:3a}.  The next step will be to sparsify the collection of tiles.  Recall the three dyadic grids 
${\mathcal D}_0$,
${\mathcal D}_1$,
${\mathcal D}_2$ from \eqref{dyad}.  One can easily verify that for every interval $J$ (not necessarily dyadic) there exists a $d \in \{0,1,2\}$ and a shifted 
dyadic interval $J' \in {\mathcal D}_d$ such that $J \subseteq J' \subseteq 3J$; we will say that $J$ is \emph{$d$-regular}.

Let $A \geq 1$, and let $d \in \{0,1,2\}$.  We shall say that a collection of $\I \subset {\mathcal D}_0$ of time intervals 
is \emph{$(A,d)$-sparse} if we have the following properties:
\begin{itemize}
\item[(i)] If $I, I' \in \I$ are such that $|I| > |I'|$, then $|I| \geq 2^{100A} |I'|$.
\item[(ii)] If $I, I' \in \I$ are such that $|I| = |I'|$ and $I \neq I'$, then $\dist(I, I') \geq 100A |I'|$.
\item[(iii)] If $I \in \I$, then $A I$ is $d$-regular, thus there exists an interval $I_{A} \in {\mathcal D}_d$ such that
$AI \subseteq I_{A} \subseteq 3AI$.  We refer to $I_A$ as the \emph{$A$-enlargement} of $I$.
\end{itemize}

If $\I$ is an $(A,d)$-sparse set of time intervals and $P$ is a tile whose time interval $I_P$ lies in $\I$, we write $I_{P,A}$ for the
$A$-enlargement of $I_P$.  Similarly if $\T$ is a tree whose time interval $I_\T$ lies in $\I$, we write $I_{\T,A}$ for the $A$-enlargement of
$I_\T$.  

To prove Theorem \ref{thm:3a}, it suffices to prove a variant for $(A,d)$-sparse sets of tiles.  More precisely, we can reduce to

\begin{theorem}[Maximal Bessel inequality, fourth reduction]
\label{thm:3b} Let $A, D, \nu > 1$, and suppose that $N$ is sufficiently large depending on $\nu$.
Let $\F$ be a forest with $\|N_\F\|_{L^\infty} \leq D$.  Let $\P := \bigcup_{\T \in \F} \T$, and suppose that the time intervals
$$ \{ I_P: P \in \P \} \cup \{ I_\T: \T \in \F \}$$
are $(A,d)$-sparse.  Let $\psi_P$, $\phi_P$ be 
as in Theorem \ref{thm:5a}.
Then there exists an exceptional set $\P_* \subset \bigcup_{\T \in \F} \T$ of tiles with
\begin{equation}\label{omega-bound-A}
|\bigcup_{P \in \P_*} I_P| \lesssim_\nu (A^{-\nu} + D^{-\nu}) \sum_{\T \in \F} |I_\T|
\end{equation}
such that we have the Bessel-type inequality
$$
\sum_{P \in \P \backslash \P_*}|\langle f,\phi_P\rangle|^2 \lesssim_\nu 
((\log(2 + AD))^{10} + A^{10-\nu} D^{10}) \|f\|_{L^2}^2
$$
 for all $f \in L^2(\R)$.
\end{theorem}

\begin{proof}[of Theorem \ref{thm:3a} assuming Theorem \ref{thm:3b}]  Let $A, \nu$ be chosen later, and set $D := \|N_\F\|_{L^\infty}$.  We need the following lemma:

\begin{lemma}[Sparsification] 
Let $\I$ be a collection of time intervals.  Then we can split $\I = \I_1 \cup \ldots \cup \I_L$ with $L = O(A^2)$ such that
each $\I_l$ for $1 \leq l \leq L$ is $(A,d)$-sparse for some $d=0,1,2$.
\end{lemma}

\begin{proof}  By pigeonholing the scale parameter into cosets of $100A \Z$, we can partition $\I$ into $100A$ subcollections,
such that on each subcollection we have the scale separation property (i) from the definition of $(A,d)$-sparseness.
Similarly if we partition the position parameter at each fixed scale into cosets of $100A$, we can partition further into
$(100A)^2$ subcollections on which we also have the position separation property (ii).  Finally, we make the elementary
observation that for each dyadic $I \in \D_0$ there exists $d=0,1,2$ such that there exists $I_A \in {\mathcal D}_d$ with
$AI \subseteq I_A \subseteq 3AI$.  A final pigeonholing based on the $d$ parameter concludes the claim.
\end{proof}

We apply this lemma to the set $\I := \{ I_P: P \in \P \}$ to split $\I$ into $\I_1,\ldots,\I_L$ for some $L = O(A^2)$.  Then we have
$\P = \P_1 \cup \ldots \cup \P_L$, where $\P_l := \{ P \in \P: I_P \in \I_l \}$.  Observe that if $\T$ is a lacunary tree in $\F$, then
$\T \cap \P_l$ is also a lacunary tree.  The time interval $I_\T$ of this tree need not lie in $\I_l$, however one can partition $\T \cap \P_l$ into
subtrees with this property.  More precisely, if we let $I$ be any interval in $\{ I_P: P \in \T \cap \P_l \}$ which is maximal with
respect to set inclusion, then $( \{ P \in \T \cap \P_l: I_P \subseteq I\}, I, \xi_\T )$ is a lacunary tree whose the time interval $I$
also lies in $\I_l$.  Let $\F_l$ be the collection of all the trees obtained in this manner for fixed $l$, as $\T$ varies over $\F$ and
$I$ varies over the maximal intervals in $\{ I_P: P \in \T \cap \P_l \}$, thus $\P_l = \bigcup_{\T \in \F_l} \T$.  
Since for each fixed $\T$ the intervals $I$ are disjoint,
one easily verifies the pointwise estimate $N_{\F_l} \leq N_\F$, and hence $\| N_{\F_l} \|_{L^\infty} \leq D$.  
Applying Theorem \ref{thm:3b} (if $N$ is large depending on $\nu$), one can then obtain an exceptional set $\P_{l,*} \subset \P_l$ obeying \eqref{omega-bound-A} such that
$$ \sum_{P \in \P_l \backslash \P_{l,*}}|\langle f,\phi_P\rangle|^2 \lesssim_\nu ((\log(1+AD))^{10} + A^{10-\nu} D^{10}) \|f\|_{L^2}^2.$$
Setting $\P_* := \bigcup_{1 \leq l \leq L} \P_{l,*}$ we thus conclude
$$ |\bigcup_{P \in \P_*} I_P| \lesssim_\nu A^2 (A^{-\nu} + D^{-\nu}) \sum_{\T \in \F} |I_\T|$$
and
$$ \sum_{P \in \P \backslash \P_*}|\langle f,\phi_P\rangle|^2 \lesssim_\nu 
A^2 ((\log(2 + AD))^{10} + A^{10-\nu} D^{10}) \|f\|_{L^2}^2.$$
If we then set $\nu := 100 +400\mu$ and $A := C_\mu D^{1/4\mu}$ for a large constant $C_\mu$ we obtain the claim.
\end{proof}

The hypothesis in Theorem \ref{thm:3b} is currently assuming some control on the quantity
$\| N_\F \|_{L^\infty} = \| \sum_{\T \in \F} 1_{I_\T} \|_{L^\infty}$.  In the arguments which follow, it is more convenient to assume
control on the larger quantity $\| \sum_{\T \in \F} \M 1_{I_\T} \|_{L^\infty}$, where of course $\M$ is the Hardy-Littlewood maximal function.  It is not necessarily the case that control of the former implies control of the latter, due to ``logarithmic pile-ups''
such as those where the intervals $I_\T$ are lacunary around a fixed origin; this is also related to the failure of the Fefferman-Stein vector-valued maximal inequality \cite{FS} at this endpoint.  Nevertheless, by removing all the tiles in a small
set it is possible to control the latter from the former.  More precisely, we have

\begin{lemma}
\label{lem:16373}
Let $\I$ be a finite set of intervals in ${\mathcal D}_d$ for some $d=0,1,2$
such that $\|\sum_{I\in\I}1_{I} \|_{L^\infty} \leq D$ for some $D$.
Then $\I$ can be split into two collections $\I=\I^{\sharp}\cup
\I^{\flat}$ such that
$$
\|\sum_{I\in\I^{\sharp}}(\M1_I)^2\|_{L^\infty}\lesssim_\nu D^{3}
$$
and
\begin{equation}\label{flatt}
|\bigcup_{I\in\I^{\flat}}I|\lesssim_\nu D^{-\nu} \sum_{I \in \I} |I|.
\end{equation}
\end{lemma}
\begin{proof}
See ~\cite[Lemma 3.14]{La}.
\end{proof}

As a consequence, we can reduce Theorem \ref{thm:3b} to

\begin{theorem}[Maximal Bessel inequality, fifth reduction]
\label{thm:3c} Let $A, M, \nu > 1$, and suppose that $N$ is sufficiently large depending on $\nu$.
Let $\F$ be a forest with 
\begin{equation}\label{mit}
\|\sum_{\T \in \F} \M 1_{I_\T}\|_{L^\infty} \leq M.
\end{equation}
Let $\P := \bigcup_{\T \in \F} \T$, and suppose that the time intervals
$$ \{ I_P: P \in \P \} \cup \{ I_\T: \T \in \F \}$$
are $(A,d)$-sparse.  
Suppose also that we have the technical condition
\begin{equation}\label{tile-cluster}
 \sup_{x \in I_P} \dist(x, \partial I_\T) \geq A^{-\nu} |I_\T|
\end{equation}
for all $P \in \P$ and $\T \in \F$ (this ensures that tiles do not cluster near the edges of trees).
Let $\psi_P$, $\phi_P$ be 
as in Theorem \ref{thm:5a}.
Then we have the Bessel-type inequality
\begin{equation}
\label{eq:drac-2-B}
\sum_{P \in \P}|\langle f,\phi_P\rangle|^2 \lesssim_\nu 
((\log(2 + AM))^{10} + A^{1-\nu} M^{2}) \|f\|_{L^2}^2
\end{equation}
 for all $f \in L^2(\R)$.
\end{theorem}

\begin{proof}[of Theorem \ref{thm:3b} assuming Theorem \ref{thm:3c}]  Apply
Lemma \ref{lem:16373} to the collection $\I := \{ I_\T: \T \in \F\}$ to create the partition $\I = \I^\sharp \cup \I^\flat$
with the desired properties.  Set
$$ \P_* := \bigcup_{\T \in \F: I_\T \in \I^\flat} \T \cup \bigcup_{\T \in \F} \{ P \in \P:
\sup_{x \in I_P} \dist(x, \partial I_\T) < A^{-\nu} |I_\T| \}.$$
Observe that
$$ \bigcup_{P \in \P_*} I_P \subseteq \bigcup_{I \in \I^\flat} I
\cup \bigcup_{\T\in \F} \{ x \in I_P:  \dist(x, \partial I_\T) < A^{-\nu} |I_\T| \}$$
and hence by \eqref{flatt}
$$ \bigcup_{P \in \P_*} I_P \lesssim_\nu (D^{-\nu} + A^{-\nu}) \sum_{\T \in I_\T} |I_\T|.$$
Now since the intervals $I_\T$ have multiplicity at most $D$, we have
$$ \|\sum_{\T \in \F: I_\T \not \in \I^\flat} \M 1_{I_\T}\|_{L^\infty} \leq D \|\sum_{I \in \I^\sharp} 
\M 1_{I}\|_{L^\infty} \lesssim_\nu D^4.$$
Applying Theorem \ref{thm:3b} with $M \sim_\nu D^4$ (and all the
trees with spatial interval in $\I^\flat$ have been completely removed from $\P \backslash \P_*$) we obtain
$$ \sum_{P \in \P \backslash \P_*} |\langle f,\phi_P\rangle|^2 \lesssim_\nu 
((\log(2 + AD^4))^{10} + A^{1-\nu} D^{8}) \|f\|_{L^2}^2$$
and the claim follows.
\end{proof}

It remains to prove Theorem \ref{thm:3c}.
We may dualize \eqref{eq:drac-2-B}, observing that it is equivalent to the estimate
$$ \| \sum_{P \in \P} a_P \phi_P \|_{L^2}^2 \lesssim_\nu 
((\log(2 + AM))^{10} + A^{1-\nu} M^{2}) \| a \|_{l^2}$$
for any sequence $a = (a_P)_{P \in \P}$ of complex numbers.  By definition of $\phi_P$,
it thus suffices to show the maximal Bessel-type inequality
\begin{equation}\label{ptarg}
 \| \sup_k |\sum_{P \in \P: |I_P| > 2^{k}} a_P \psi_P| \|_{L^2} \lesssim_\nu 
((\log(2 + AM))^{10} + A^{1-\nu} M^{2}) \| a \|_{l^2}.
\end{equation}

At this point we shall pause to sketch the general strategy we shall employ to prove \eqref{ptarg}, following
\cite{La}.  First we shall split the tile set $\P$ into layers $\P = \P_1 \cup \ldots \cup \P_J$.  Roughly speaking, 
the idea is to arrange these layers so that the time intervals of $\P_{j'}$ tend to be 
(locally) wider than those of $\P_{j}$ when $j' < j$.  Since $\psi_P$ is essentially concentrated in $I_P$ (or more accurately
$I_{P,A}$), this heuristically gives rise to an estimate of the form
$$ |\sup_k |\sum_{P \in \P: |I_P| > 2^{k}} a_P \psi_P| \leq
\sup_j [ |\sum_{P \in \bigcup_{j' < j} \P_{j'}} a_P \psi_P| + \sup_k |\sum_{P \in \P_j: |I_P| > 2^k} a_P \psi_P| ].$$
To deal with the former expression we shall use the Radamacher-Menshov inequality and a non-maximal Bessel inequality (which is essentially \eqref{ptarg} without the supremum in $k$, and with somewhat fewer logarithmic losses on the right-hand side).
To deal with the second term we replace the supremum in $j$ by a square function, and reduce to
controlling the contribution of a localized expression over a single generation $\P_j$ (which will ultimately reduce to
a certain maximal inequality of Bourgain \cite{Bo1}).  

For technical reasons it turns out that one needs to treat the ``boundary'' of the layers $\P_j$ separately from the rest of
the $\P_j$, in order to improve the separation properties between layers.  As such we will have to execute the above strategy
twice, once for the boundary tiles and once for the interior tiles.

We now turn to the details, beginning with the selection of the layers.  Introduce the sets $\I \subset {\mathcal D}_0$ and $\I_A \subset {\mathcal D}_d$ by
$$ \I := \{ I_\T: \T \in \F \}; \quad \I_A := \{ I_{\T,A}: \T \in \F \}.$$
Observe that the $(A,d)$-sparseness of $\I$ ensures that the map $I \mapsto I_A$ is a bijection from $\I$ to $\I_A$ which preserves the set inclusion
relation.
Since $AI_\T \subseteq I_{\T,A} \subseteq 3AI_\T$ we see that
$$ 1_{\T,A} \subseteq 10A 1_{\T}$$
and hence by \eqref{mit} we have the multiplicity bound
\begin{equation}\label{tamult}
\| \sum_{I \in \I_A} 1_I \|_{L^\infty} \leq 10AM.
\end{equation}
We then partition
$$ \I_A = \I_{A}^{(1)} \cup \I_{A}^{(2)} \cup \ldots \cup \I_{A}^{(10AM)}$$
recursively by defining $\I_{A}^{(j)}$ to be those intervals in $\I_A \backslash \bigcup_{i < j} \I_{A}^{(i)}$ which are maximal with respect
to set inclusion, thus $\I_A^{(j)}$ is a collection of disjoint intervals in the dyadic grid ${\mathcal D}_d$.  
Observe that for $1 < j \leq 10AM$, each interval in $\I_A^{(j)}$ is contained in exactly one interval in $\I_A^{(j-1)}$;
since $\| \sum_{I \in \I} 1_I \|_{L^\infty} \leq 10AM$, we conclude that $\I_A^{(1)},\ldots,\I_A^{(10AM)}$ do indeed partition $\I_A$.
Using the bijection between $\I$ and $\I_A$, we thus induce a partition $\I = \I^{(1)} \cup \ldots \cup \I^{(10AM)}$ of $\I$.

Let $1 \leq j \leq 10AM$.  For each $I \in \I^{(j)}$, let $\P_I$ denote the tiles with  time interval $I$:
$$ \P_I := \{ P \in \P: I_P = I \}.$$
Observe that each tree $\T$ in $\F$ contributes at most one tile to $\P_I$, by definition of a lacunary tree, and if $\T$ does contribute a tile then $1_I \leq 1_{\T}$.  By \eqref{mit} we thus have
\begin{equation}\label{pi-card}
\# \P_I \leq M \hbox{ for all } I \in \I.
\end{equation}
We also introduce the tileset $\P_{<I}$ for $I \in \I^{(j)}$ by
$$ \P_{<I} := \{ P \in \P: I_P \subsetneq I; I_P \not \subseteq J \hbox{ for all } J \in \bigcup_{i > j} \I^{(i)} \};$$
thus $\P_{<I}$ is the collection of tiles whose time interval is strictly contained in the interval $I \in \I^{(j)}$, but is not contained
in any interval from a later layer of $\I$.  Since every tile $P \in \P$ has its time interval $I_P$ contained in some interval in $\I$
(because $P$ is contained in some tree $\T \in \F$, and hence $I_P \subseteq I_\T \in \I$) we see that we have the partition
$$ \P = \bigcup_{1 \leq j \leq 10AM} \bigcup_{I \in \I^{(j)}} (\P_I \cup \P_{<I}).$$
To prove \eqref{ptarg}, it thus suffices by the triangle inequality to prove the estimates
\begin{equation}\label{ptop}
 \| \sup_k |\sum_j \sum_{I \in \I^{(j)}} \sum_{P \in \P_I: |I_P| > 2^k} a_P \psi_P| \|_{L^2} \lesssim_\nu 
( (\log(2 + AM))^{10} + A^{1-\nu} M^2) \| a \|_{l^2}.
\end{equation}
and
\begin{equation}\label{pnotop}
 \| \sup_k |\sum_j \sum_{I \in \I^{(j)}} \sum_{P \in \P_{<I}: |I_P| > 2^{k}} a_P \psi_P| \|_{L^2} \lesssim_\nu 
( (\log(2 + AM))^{10} + A^{1-\nu} M^2) \| a \|_{l^2}.
\end{equation}

The estimate \eqref{ptop} is easier and is proven in Section \ref{ptop-est}.  The estimate \eqref{pnotop} is more difficult, relying in particular on a certain inequality of Bourgain, and is proven in Section \ref{pnotop-est}.  To conclude this section, we present two tools which
will be used to prove both \eqref{ptop} and \eqref{pnotop}.  The first is a non-maximal Bessel inequality, and more precisely
the bound
\begin{equation}\label{nonmax}
 \| \sum_{P \in \P} a_P \psi_P \|_{L^2}^2 \lesssim \log(1+M) \| a \|_{l^2}. 
\end{equation}
This inequality may be of some independent interest and is proven in Section \ref{nonmax-sec}.  
Secondly, we will rely on the following form of the standard Radamacher-Menshov inequality, whose proof we include for
sake of completeness.  We observe first the trivial bound
\begin{equation}\label{square}
\| \sup_i |f_i| \|_{L^2} \leq \| (\sum_i |f_i|^2)^{1/2} \|_{L^2} = (\sum_i \|f_i\|_{L^2}^2)^{1/2}
\end{equation}
valid for any finite collection of $L^2$ functions $f_i$.  This bound is usually too crude for applications, as the summation in $i$ usually creates
an undesirable polynomial loss in the estimates, however one can refine this polynomial loss to a logarithmic loss in the following way.

\begin{theorem}[Radamacher-Menshov]
\label{thm:Be}  
Let $(f_l)_{l=1}^L$ be a sequence of functions in $L^2(\R)$ which are almost orthogonal in the sense that there exists a  constant $B$, such that for each finite sequence $\epsilon_1,\ldots,\epsilon_L \in \{-1,+1\}$ of signs we have
$$\left\|\sum_{l=1}^L \epsilon_lf_l\right\|_{L^2}\le B.$$
Then we have the maximal inequality
$$\left\|\sup_{L'\le L}|\sum_{l=1}^{L'}f_l|\right\|_{L^2}\lesssim B \log(2+L).$$
\end{theorem}

\begin{proof} We may take the $f_l$ to be real-valued.
By adding dummy $f_l$ if necessary, we may assume that $L=2^m$ for some integer $m\ge 1$. For each set $I \subseteq \{1,\ldots,L\}$ let
$f_I := \sum_{l \in I} f_l$.  For each $0 \leq m' \leq m$ let $\I_{m'}$ denote the collection of sets of the form
$\{ 2^{m'} j + 1, \ldots, 2^{m'} j + 2^{m'} \}$ for $j=0,\ldots,2^{m-m'}-1$.  For each fixed $m'$, the sets in $\I_{m'}$ 
partition $\{1,\ldots,L\}$, and thus by hypothesis we have
$$ \| \sum_{I \in \I_{m'}} \epsilon_I f_I \|_{L^2} \le B$$
for all signs $\epsilon_I = \pm 1$.  If we square this inequality we obtain
$$ \sum_{I \in \I_{m'}} \|f_I\|_{L^2}^2 + \sum_{I,J \in \I_{m'}: I \neq J} \epsilon_I \epsilon_J \langle f_I, f_J \rangle
\leq B^2.$$
If we then set $\epsilon_I$ to be independent random signs and take expectations, we conclude
$$ \sum_{I \in \I_{m'}} \|f_I\|_{L^2}^2 \leq B^2.$$
By \eqref{square} this implies that
$$ \| \sup_{I \in \I_{m'}} |f_I| \|_{L^2} \leq B.$$
By representing $L'$ in binary and using the triangle inequality we have the pointwise estimate
$$ |\sum_{l=1}^{L'} f_l| \leq \sum_{0 \leq m' \leq m} \sup_{I \in \I_{m'}} |f_I| $$
for all $L' \le L$.  Taking suprema over all $L'$, taking $L^2$ norms, and applying the triangle inequality, the claim follows.
\end{proof}

\section{Proof of \eqref{ptop}}\label{ptop-est}

We first prove the estimate \eqref{ptop}, which is relatively easy, and serves as a model for the more complicated
estimate \eqref{pnotop}.  

Intuitively, the contribution of the wave packets $\psi_P$ for $P \in \P_I$ should be localized to
the time interval $I_A$.  To exploit this we introduce the tail error
$$ E(x) := \sum_j \sum_{I \in \I^{(j)}: x \not \in I_A} \sum_{P \in \P_I} |a_P| |\psi_P(x)|.$$

This error is small:

\begin{lemma}[Tail estimate]\label{exterior-estimate}  We have
$$
 \|E\|_{L^2} \lesssim_\nu A^{-\nu} M^{1/2} \| a \|_{l^2}.
$$
\end{lemma}

\begin{proof}  
From \eqref{psisj} one easily verifies the pointwise estimates
$$ |\psi_P| \lesssim |I|^{-1/2} \M 1_I$$
and the $L^1$ bound 
$$ \| |\psi_P| (1 - 1_{I_A}) \|_{L^1} \lesssim A^{-10N+10} |I|^{1/2}$$
whenever $P \in \P_I$.
The former bound and \eqref{mit} implies the estimate
$$
\| \sum_j \sum_{I \in \I^{(j)}: x \not \in I_A} \sum_{P \in \P_I} 
|I|^{1/2} |a_P| |\psi_P(x)| \|_{L^\infty} \lesssim
M \| a \|_{l^\infty}$$
while the latter bound and the triangle inequality implies the bound
$$
\| \sum_j \sum_{I \in \I^{(j)}: x \not \in I_A} \sum_{P \in \P_I} 
|I|^{-1/2} |a_P| |\psi_P(x)| \|_{L^1} \lesssim
A^{-10N+10} \| a \|_{l^1}.$$
The claim then follows from interpolation (or from Cauchy-Schwarz), since we assume $N$ sufficiently large depending on $\nu$. 
\end{proof}

To exploit this tail estimate we use the following pointwise inequality:

\begin{lemma}\label{pointlemma}  For almost every $x$ we have
$$ \sup_k |\sum_j \sum_{I \in \I^{(j)}} \sum_{P \in \P_I: |I_P| > 2^k} a_P \psi_P(x)|
\lesssim
\sup_{j_0} |\sum_{j \leq j_0} \sum_{I \in \I^{(j)}} \sum_{P \in \P_I} a_P \psi_P(x)|+ E(x).$$
\end{lemma}

\begin{proof} We may assume that $x$ is not the endpoint of any dyadic interval.  It suffices to show that for every $k$ and $x$ there
exists a $j_0$ such that
$$ |\sum_j \sum_{I \in \I^{(j)}} \sum_{P \in \P_I: |I_P| > 2^k} a_P \psi_P(x)|
\lesssim |\sum_{j \leq j_0} \sum_{I \in \I^{(j)}} \sum_{P \in \P_I} a_P \psi_P(x)|+ O(E(x)).$$
Since $I_P = I$, we can write the left-hand side as
$$ |\sum_j \sum_{I \in \I^{(j)}: |I| > 2^k} \sum_{P \in \P_I} a_P \psi_P(x)|.$$
By definition of $E$ and the triangle inequality, we can bound this by
$$ |\sum_j \sum_{I \in \I^{(j)}: |I| > 2^k; x \in I_A} \sum_{P \in \P_I} a_P \psi_P(x)| + O(E(x)).$$
For each $1 \leq j \leq 10AM$, we know that there is at most one interval $I_j \in \I^{(j)}$ whose dilate $I_{j,A}$
contains $x$, and furthermore these intervals are decreasing in $j$ (adopting the convention that $I_j = \emptyset$ if no interval in $\I^{(j)}$ contains $x$).  Thus if we let $j_0$ be the largest $j$ for which $|I_{j_0}| > 2^{k(x)}$ (with $j_0=0$
if no such $j$ exists), then we see that if $1 \leq j \leq 10AM$ and $I \in \I^{(j)}$ are such that $x \in I_A$, then
$|I|> 2^{k(x)}$ if and only if $j \leq j(x)$.  Thus we can bound the preceding expression by
$$ |\sum_{j \leq j_0} \sum_{I \in \I^{(j)}: x \in I_A} \sum_{P \in \P_I} a_P \psi_P(x)| + O(E(x)).$$
One can then remove the constraint $x \in I_A$ by definition of $E(x)$ and the triangle inequality.
\end{proof}

In light of the above two lemmas, we see that to prove \eqref{ptop} it suffices to show that
\begin{equation}\label{potato}
 \|\sup_{j_0} |\sum_{j \leq j_0} \sum_{I \in \I^{(j)}} \sum_{P \in \P_I} a_P \psi_P| \|_{L^2} \lesssim
(\log(2 + AM))^{10} \| a \|_{l^2}.
\end{equation}
Applying the Radamacher-Menshov inequality (Theorem \ref{thm:Be}), it suffices to show that
$$
 \|\sum_j \epsilon_j \sum_{I \in \I^{(j)}} \sum_{P \in \P_I} a_P \psi_P \|_{L^2} \lesssim 
(\log(2 + AM))^{9} \| a \|_{l^2}.
$$
for all choices $\epsilon_1,\ldots,\epsilon_{10AM} \in \{-1,+1\}$ of signs.  
But this follows from the non-maximal Bessel inequality \eqref{nonmax} (with some room to spare), since the $\P_I$ are disjoint in $\P$.
This concludes the proof of \eqref{ptop}.
 
\section{Proof of \eqref{pnotop}}\label{pnotop-est}

Now we prove \eqref{pnotop}.  We shall argue as in the proof of \eqref{ptop}, although the details shall be more technical, and we shall also
rely crucially on a maximal inequality of Bourgain.

In the previous section we localized the contribution of $\P_I$ to the interval $I_A$.
It turns out (because of the $(A,d)$-sparseness hypothesis)
that the contribution of $\P_{<I}$ can be localized even further, to the interval $I$ itself.  To formalize this we again introduce a tail
error
$$ \tilde E(x) := \sum_{I \in \I^{(j)}: x \not \in I} \sum_{P \in \P_{<I}} |a_P| |\psi_P(x)|.$$ 

\begin{lemma}[Tail estimate]\label{exterior-estimate-2}  We have
$$
 \| \tilde E \|_{L^2} \lesssim_\nu 
A^{1-\nu} M^{2} \| a \|_{l^2}
$$
\end{lemma}

\begin{proof} Since there are only $10AM$ values of $j$, it suffices by the triangle inequality to show that
\begin{equation}\label{pip}
 \| \sum_{I \in \I^{(j)}: x \not \in I} \sum_{P \in \P_{<I}} |a_P| |\psi_P(x)| \|_{L^2} \lesssim_\nu 
A^{-\nu} M \| a \|_{l^2}
\end{equation}
for each $j$, which we now fix.  Suppose for the moment that we could show the pointwise estimate
\begin{equation}\label{pointpip}
 \sum_{P \in \P_{<I}} |a_P| |\psi_P(x)| \lesssim_\nu A^{-\nu} M c_I |I|^{-1/2} \M 1_I(x)^2
\end{equation}
for each $I$ and $x \not \in I$, where $c_I := (\sum_{P \in \P_{<I}} |a_P|^2)^{1/2}$.  Then the left-hand side of
\eqref{pip} is bounded by
$$ \lesssim_\nu A^{-\nu} M \| \sum_{I \in \I^{(j)}} c_I |I|^{-1/2} \M 1_I(x)^2 \|_{L^2}.$$
Applying the Fefferman-Stein maximal inequality \cite{FS}, which among other things asserts that
$$ \| \sum_i \M f_i^2 \|_{L^2} = \| ( \sum_i \M f_i^2)^{1/2} \|_{L^4}^2 \lesssim 
\| ( \sum_i |f_i|^2)^{1/2} \|_{L^4}  = \| \sum_i |f_i|^2\|_{L^2},$$
we can bound the left-hand side of \eqref{pip} by
$$ \lesssim_\nu A^{-\nu} M \| \sum_{I \in \I^{(j)}} c_I |I|^{-1/2} 1_I(x)^2 \|_{L^2}.$$
 Since the intervals in $\I^{(j)}$ are disjoint, this expression is bounded by
$$ A^{-\nu} M (\sum_{I \in \I^{(j)}} |c_I|^2)^{1/2} \leq A^{-\nu} M \| a \|_{l^2}$$
as desired.

It remains to prove \eqref{pointpip}.  By Cauchy-Schwarz it suffices to verify the estimates
\begin{equation}\label{pp1}
 \sum_{P \in \P_{<I}} |I_P| |\psi_P(x)| \lesssim M^2 |I|^{1/2}
 \end{equation}
and
\begin{equation}\label{pp2}
 \sum_{P \in \P_{<I}} |a_P|^2 |I_P|^{-1} |\psi_P(x)| \lesssim_\nu A^{-2\nu} c_I^2 |I|^{-3/2} \M 1_I(x)^4.
 \end{equation}
To prove \eqref{pp1}, we break $\P_{<I}$ up into $\P_{<I} \cup \T$, where $\T$ ranges over the forest $\F$.  Observe that
$\P_{<I} \cup \T$ is empty unless $I \subseteq I_\T$; by \eqref{mit} we thus see that there are at most $M$ trees $\T$ for which
$\P_{<I} \cup \T$ is non-empty.  Thus it suffices to show that
$$  \sum_{P \in \P_{<I} \cap \T} |I_P| |\psi_P(x)| \lesssim |I|^{1/2}.$$
But from \eqref{psisj} we have $|I_P| |\psi_P(x)| \lesssim |I_P|^{1/2} (\M 1_{I_P}(x))^{100}$ (say).  Since the $I_P$ are dyadic subintervals of $I$
and each interval can occur at most once in $\T$, the claim follows.

It remains to prove \eqref{pp2}.  From the definition of $c_I$ and the triangle inequality it suffices to prove that
$$ |I_P|^{-1} |\psi_P(x)| \lesssim_\nu A^{-2\nu} |I|^{-3/2} \M 1_I(x)^4$$
for each $P \in \P_{<I}$ and $x \not \in I$.  But from the $(A,d)$-sparseness hypothesis we see that
$I_P \subsetneq I$ and $|I_P| \leq 2^{-100A} |I|$, while from \eqref{tile-cluster} (recalling that $I$ is the time interval of some tree $\T$)
we have $\sup_{x \in I_P} \dist(x, \partial I) \geq A^{-\nu} |I|$.  The claim now follows from \eqref{psisj}, the exponential gain of $|I|/|I_P| \geq 2^{100A}$ being more than sufficient to compensate for any polynomial losses in $A$ or in $|I|/|I_P|$.
\end{proof}

The analog of Lemma \ref{pointlemma} is

\begin{lemma}\label{pointlemma-2}  For almost every $x$ we have
\begin{align*} 
\sup_k |\sum_j \sum_{I \in \I^{(j)}} \sum_{P \in \P_{<I}: |I_P| > 2^k} a_P \psi_P(x)|
&\lesssim
\sup_{j_0} |\sum_{j < j_0} \sum_{I \in \I^{(j)}} \sum_{P \in \P_{<I}} a_P \psi_P(x)| \\
&\quad + 
\sup_{I \in \I} \sup_k |\sum_{P \in \P_{<I}: |I_P| > 2^k} a_P \psi_P(x)|\\
&\quad + \tilde E(x).
\end{align*}
\end{lemma}

\begin{proof} We again may assume that $x$ is not the endpoint of a dyadic interval.  We fix $k$; it would suffice 
to find a $j_0$ and an $I_0 \in \I$ such that
\begin{align}
|\sum_j \sum_{I \in \I^{(j)}} \sum_{P \in \P_{<I}: |I_P| > 2^k} a_P \psi_P(x)|
&\leq
|\sum_{j < j_0} \sum_{I \in \I^{(j)}} \sum_{P \in \P_{<I}} a_P \psi_P(x)| \nonumber\\
&\quad + |\sum_{P \in \P_{<I_0}: |I_P| > 2^k} a_P \psi_P(x)|\label{i0}\\
&\quad + O(\tilde E(x)).\nonumber
\end{align}
By definition of $\tilde E(x)$, we have
$$
|\sum_j \sum_{I \in \I^{(j)}} \sum_{P \in \P_{<I}: |I_P| > 2^k} a_P \psi_P(x)|
\leq |\sum_j \sum_{I \in \I^{(j)}; x \in I} \sum_{P \in \P_{<I}: |I_P| > 2^k} a_P \psi_P(x)| + O( \tilde E(x) ).$$
Let $j_0$ be the largest $j$ for which there exists an interval in $\I^{(j_0)}$ which contains $x$ and has length greater than $2^k$.  
There is only one such interval; call it $I_0$.  We can thus estimate the contribution of the $j=j_0$ term by
\eqref{i0}, and reduce to showing that
$$ |\sum_{j < j_0} \sum_{I \in \I^{(j)}; x \in I} \sum_{P \in \P_{<I}: |I_P| > 2^k} a_P \psi_P(x)|
\lesssim 
|\sum_{j < j_0} \sum_{I \in \I^{(j)}} \sum_{P \in \P_{<I}} a_P \psi_P(x)| + O(\tilde E(x)).$$
But if $I \in \I^{(j)}$ and $x \in I$ then $I$ and $I_0$ overlap.  Since $I_0$ belongs to a later layer $\I^{(j_0)}$ than $I$ we must
have $I_0 \subsetneq I$, and thus $|I| > |I_0| > 2^k$.  Hence the constraint $|I_P| > 2^k$ is redundant and can be removed. The claim
now follows from the triangle inequality.
\end{proof}

In light of the above two lemmas, to prove \eqref{pnotop} it would suffice to show that
$$ \| \sup_{j_0} |\sum_{j < j_0} \sum_{I \in \I^{(j)}} \sum_{P \in \P_{<I}} a_P \psi_P| \|_{L^2}
\lesssim (\log(2 + AM))^{10} \| a \|_{l^2}$$
and
$$ \| \sup_{I \in \I} \sup_k |\sum_{P \in \P_{<I}: |I_P| > 2^k} a_P \psi_P(x)| \|_{L^2}
\lesssim
(\log(2 + AM))^{10}  \| a \|_{l^2}.$$
The first inequality is proven in exactly the same way as \eqref{potato} and is omitted, so we now turn
to the second inequality.  By \eqref{square} it would suffice to show that
$$ (\sum_{I \in \I} 
\| \sup_k |\sum_{P \in \P_{<I}: |I_P| > 2^k} a_P \psi_P(x)| \|_{L^2}^2)^{1/2}
\lesssim (\log(2 + AM))^{10}  \| a \|_{l^2},$$
which in turn would follow from the estimate
$$\| \sup_k |\sum_{P \in \P_{<I}: |I_P| > 2^k} a_P \psi_P(x)| \|_{L^2}
\lesssim (\log(2 + AM))^{10} (\sum_{P \in \P_{<I}} |a_P|^2)^{1/2}$$
for each fixed $I$.

Let $\T_1, \T_2, \ldots, \T_J$ be all the trees in $\F$ which intersect $\P_{<I}$; the time interval of such trees must contain $I$, and so 
from \eqref{mit} we have $J \leq M$. 
We can then write
$$ \sum_{P \in \P_{<I}: |I_P| > 2^k} a_P \psi_P(x) = \sum_{j=1}^J \sum_{P \in \P_{<I} \cap \T_j: |I_P| > 2^k} a_P \psi_P(x).$$
Let $\xi_1, \ldots, \xi_J$ be the base frequencies of $\T_1,\ldots,\T_J$.
Since $\T_j$ is a lacunary tree, we see that if $P \in \P_{<I} \cap \T_j$ then $\psi_P$ has Fourier support in an interval of width $|I_P|^{-1}$
and distance $\sim C_0 |I_P|^{-1}$ from $\xi_j$.  By the strong disjointness of the $\T_j$ we see that these intervals must be disjoint.
This implies that
$$ \sum_{j=1}^J \sum_{P \in \P_{<I} \cap \T_j: |I_P| > 2^k} a_P \psi_P(x) = \Pi_k \sum_{j=1}^J \sum_{P \in \P_{<I} \cap \T_j} a_P \psi_P(x)$$
where $\Pi_k$ is a Fourier projection to the union of $J$ intervals centered at $\xi_1,\ldots,\xi_J$, each of radius $\sim C_0 2^{-k}$.
We now invoke a deep maximal inequality of Bourgain \cite[Lemma 4.11]{Bo1}, which asserts in our notation that
$$ \| \sup_k |\Pi_k f| \|_{L^2} \lesssim \log(2+J)^2 \|f\|_{L^2}.$$
Using this, we reduce to showing that
$$\| \sum_{j=1}^J \sum_{P \in \P_{<I} \cap \T_j} a_P \psi_P| \|_{L^2}
\lesssim (\log(2 + AM))^{8} (\sum_{P \in \P_{<I}} |a_P|^2)^{1/2}.$$
But this follows from the non-maximal Bessel inequality \eqref{nonmax}.
This concludes the proof of \eqref{pnotop}.

\section{Proof of \eqref{nonmax}}\label{nonmax-sec}

We now prove \eqref{nonmax}.  We shall in fact prove the slightly more general statement, which may have some independent interest:

\begin{proposition}[Nonmaximal Bessel inequality]\label{nmprop}  Let $\F$ be a forest, let $\P := \bigcup_{\T \in \F} \T$, and for each tile $P \in \P$ let $\psi_P$ be a wave packet adapted to $P$.  Suppose also that $\| \sum_{\T \in \F} 1_\T \|_{L^\infty} \leq M$.  Then we have
$$ \| \sum_{P \in \P} a_P \psi_P \|_{L^2} \lesssim \log(2+M) \| a \|_{l^2}$$
for any sequence $a = (a_P)_{P \in \P}$ of complex numbers.
\end{proposition}

\begin{remark} By duality and the $TT^*$ method, this inequality is also equivalent to the assertion that
$$(\sum_{P \in \P} |\langle f, \psi_P \rangle|^2)^{1/2} \lesssim \log(2+M) \| f \|_{L^2}$$
or that
$$ \| \sum_{P \in \P} \langle f, \psi_P \rangle \psi_P \|_{L^2} \lesssim \log(2+M)^2 \| f \|_{L^2}$$
for all $f \in L^2$. The logarithmic loss can probably be lowered to $\log(2+M)^{1/2}$ but cannot be removed
entirely; see \cite{BL}.
\end{remark}

We prove Proposition \ref{nmprop} in stages.  The most important step is to establish a restricted version of the proposition without the logarithmic loss in $M$.

\begin{proposition}[Restricted Bessel inequality]\label{nm-restrict}  Let $\F, \P, \psi_P,M$ be as in Proposition \ref{nmprop}.  Suppose that
$a = (a_P)_{P \in \P}$ obeys the Carleson condition $\sum_{P \in \T'} |a_P|^2 \lesssim 2^{2m} |I_{\T'}|$ for all $\T \in \F$ and all subtrees $\T'$ of $\T$, where $m$ is a fixed integer.  Then we have
$$ \| \sum_{P \in \P} a_P \psi_P \|_{L^2} \lesssim 2^{m} (\sum_{\T \in \F} |I_\T|)^{1/2}.$$
\end{proposition}

\begin{proof} See \cite[Lemma 6.6]{MTT2}.  The main idea is to square both sides, use standard
estimates on the inner products $|\langle \psi_P, \psi_Q \rangle|$, and exploit the strong disjointness of the trees $\T$ in the forest $\F$.
\end{proof}

Next, we establish restricted $L^p$ type estimates with a polynomial loss in $M$.

\begin{proposition}[Crude Bessel inequality]\label{nm-crude}  Let $\F, \P, \psi_P,M, m, a$ be as in Proposition \ref{nm-restrict}.  
Then for any $1 < p < \infty$ we have
$$ \| \sum_{P \in \P} a_P \psi_P \|_{L^p} \lesssim_p 
2^{m} M (\sum_{\T \in \F} |I_\T|)^{1/p}.$$
\end{proposition}

\begin{remark} One can improve the factor of $M$ here by interpolation with Proposition \ref{nm-restrict}, and at the endpoint $p=1$ one
can remove the loss in $M$ entirely.  However for our purposes any polynomial factor in $M$ will suffice.
\end{remark}

\begin{proof} First observe that we can partition the forest $\F$ into forests $\F_1 \cup \ldots \cup \F_M$, with each $\F_j$ having multiplicity one
in the sense that $\| \sum_{\T \in \F_j} 1_\T \|_{L^\infty} \leq 1$.  Indeed one could set $\F_M$ to be a maximal collection of trees in $\F$ whose time intervals are distinct and are maximal with respect to set inclusion, remove $\F_M$ from $\F$ (dropping the multiplicity by $1$), and induct;
we leave the details to the reader.  From the triangle inequality we see that it thus suffices to verify the claim when $M=1$. We may also normalize
$m=0$.  

Let $\T$ be a tree in $\F$.  We can partition the dyadic interval $I_\T$ into four equally sized dyadic 
sub-intervals $I_{\T,1}, I_{\T,2}, I_{\T,3}, I_{\T,4}$, from left to right.  Let $\T_l, \T_r \subset \T$ be the trees
$\T_l := \{ P \in \T: I_P \subseteq I_{\T,1} \}$ and $\T_r := \{ P \in \T: I_P \subseteq I_{\T,4} \}$ with spatial intervals
$I_{\T,1}$ and $I_{\T,4}$ respectively, and let $\F'$ be the forest formed by these trees $\T_l$ and $\T_r$, and $\P' := \bigcup_{\T \in \F'} \T$.  Observe that this forest also
has multiplicity one, and that $\sum_{\T \in \F'} |I_\T| = \frac{1}{2} \sum_{\T \in \F} |I_\T|$.  It thus suffices by the obvious recursion argument
to prove the Bessel inequality with $\P$ replaced by $\P \backslash \P'$ (conceding a factor of $\frac{1}{1 - 2^{-p}} \sim p$ 
in the implicit constant).  The practical upshot of this reduction is that for any tree $\T$ in the forest $\F$,
we may assume without loss of generality that none of the tiles in $\T$ have time interval contained in the left quarter $I_{\T,l}$ or right
quarter $I_{\T,r}$ of the tree.

From the Carleson condition we have the crude bound $|a_P| \lesssim |I_P|^{1/2}$ for all $P \in \P$.  From this, \eqref{psisj}, and the above reduction on the trees $\T$ one easily verifies the pointwise estimate
$$ (1 - 1_{I_\T}(x)) |\sum_{P \in \T} a_P \psi_P(x)| \lesssim \M 1_{I_\T}(x)^{10}$$
(say) for all $\T \in \F$ and $x \in \R$.  From the Fefferman-Stein maximal inequality \cite{FS} and the multiplicity one nature of $\F$ we thus have
$$ \| \sum_{\T \in \F} (1 - 1_{I_\T}(x)) |\sum_{P \in \T} a_P \psi_P(x)| \|_{L^p} \lesssim
\| \sum_{\T \in \F} \M 1_{I_\T}(x)^{10} \|_{L^p} \lesssim_p (\sum_{\T \in \F} |I_\T|)^{1/p}$$
and hence by the triangle inequality it will suffice to show that
$$ \| \sum_{\T \in \F} 1_{I_\T} |\sum_{P \in \T} a_P \psi_P| \|_{L^p} 
 \lesssim_p (\sum_{\T \in \F} |I_\T|)^{1/p}.$$
From the disjointness of the intervals $I_\T$ it thus suffices to show that
$$ \| \sum_{P \in \T} a_P \psi_P \|_{L^p(I_\T)} 
 \lesssim_p |I_\T|^{1/p}$$
for each tree $\T$.  By shifting the frequency dyadic grid if necessary we may 
assume that $\xi_T = 0$; this essentially turns the wave packets $\psi_P$ into wavelets.
The Carleson condition on the $a_P$ and standard almost orthogonality estimates then give the $L^2$ estimate
$$ \| \sum_{P \in \T} a_P \psi_P \|_{L^2} \lesssim |I_\T|^{1/2}$$
and the $\BMO$ estimate
$$ \| \sum_{P \in \T} a_P \psi_P \|_{\BMO} \lesssim 1$$
from which the claim follows by the John-Nirenberg inequality.
\end{proof}

The idea is now to combine the above two propositions via some sort of real interpolation method to obtain Proposition \ref{nmprop}.
It may well be possible to use one of the existing real interpolation theorems in the literature to obtain this conclusion, but we will
use a more explicit argument, based on the following decomposition of an arbitrary $l^2$ sequence $a$ into Carleson sequences.

\begin{lemma}[Stopping time algorithm]\label{nm-stop}  Let $\F, \P, \psi_P,M$ be as in Proposition \ref{nmprop}.  Suppose that
$a = (a_P)_{P \in \P}$ obeys the Carleson condition $\sum_{P \in \T'} |a_P|^2 \leq 2^{2m} |I_{\T'}|$ for all $\T \in \F$ and all subtrees $\T'$ of $\T$, where $m$ is a fixed integer.  Then we can partition $\P = \bigcup_{\T \in \F_1} \T \cup \bigcup_{\T \in \F_2} \T$, where $\F_1$ is a collection
of subtrees of trees in $\F$ such that
\begin{equation}\label{sumest}
2^{2m} \sum_{\T \in \F_1} |I_\T| \sim \sum_{\T \in \F_1} \sum_{\P \in \T} |a_P|^2 
\end{equation}
and $\F_2$ is a collection of subtrees of trees in $\F$ such that
$\sum_{P \in \T'} |a_P|^2 \leq 2^{2(m-1)} |I_{\T'}|$ for all $\T \in \F_2$ and all subtrees $\T'$ of $\T$. Furthermore we have
$\|\sum_{\T \in \F_j} 1_{I_\T}\|_{L^\infty} \leq M$ for $j=1,2$.
\end{lemma}

\begin{proof}  It suffices to establish this lemma in the case when the forest $\F$ consists of a single tree, $\F = \{\T\}$, with $M=1$,
since the general case then follows by applying the lemma to each tree separately and taking unions, using Lemma \ref{forestry}, as well
as the observation that the contribution to $\sum_{\T \in \F_j} 1_{I_\T}$ arising from a single tree $\T_0$ in $\F$ will be bounded pointwise
by $1_{\T_0}$.

Let $\I$ be the set of all dyadic intervals $I$ in $I_\T$ such that $\sum_{P \in \T: I_P \subseteq I} |a_P|^2 > 2^{2(m-1)} |I|$, and such that $I$
is maximal with respect to set inclusion among all such intervals with the property; thus the intervals in $\I$ are disjoint and lie inside $I_\T$.  We then let $\F_1$ be the forest consisting of trees $\T_I = (\T_I, I, \xi_T)$ of the form $\T_I := \{ P \in \T: I_P \subseteq I \}$, where $I$ ranges
over $\I$.  By construction it is clear that $\F_1$ is indeed a forest, and that $\sum_{\P \in \T'} |a_P|^2 \sim 2^{2m} |I_{\T'}|$ for all
$\T' \in \F_1$; summing over all $\T'$ we obtain \eqref{sumest}.  If we let $\F_2$ consist of the single tree $\T_2 = (\T_2, I_{\T}, \xi_\T)$
consisting of all the tiles not covered by $\F_1$, thus $\T_2 := \T \backslash \bigcup_{\T' \in \F_1} \T'$, then we see from construction that
$\sum_{P \in \T'} |a_P|^2 \leq 2^{2(m-1)} |I_{\T'}|$ for all subtrees $\T'$ of $\T_2$.  The claim follows.
\end{proof}

Iterating this lemma in the usual manner, starting with $m$ extremely large and exploiting the fact that the forest $\F$ contains only finitely many
tiles, we obtain
  
\begin{corollary}[Iterated stopping time algorithm]\label{nm-stop-iter}  Let $\F, \P, \psi_P,M$ be as in Proposition \ref{nmprop}.  Then there
exist forests $\F_m$ for each integer $m$, together with a tile set $\P_{-\infty}$, such that we have the partition
$$ \P = \bigcup_m \bigcup_{\T \in \F_m} \T \cup \P_{-\infty},$$
such that we have the Carleson condition $\sum_{P \in \T'} |a_P|^2 \leq 2^{2m} |I_{\T'}|$ for all $m$, all $\T \in \F_m$ and all 
subtrees $\T'$ of $\T$, we have the bound
\begin{equation}\label{msum}
\sum_m 2^{2m} \sum_{\T \in \F_m} |I_\T| \sim \sum_{P \in \P} |a_P|^2 
\end{equation}
and such that $a_P = 0$ for all $P \in \P_{-\infty}$.  Finally we have
$\|\sum_{\T \in \F_m} 1_{I_\T}\|_{L^\infty} \leq M$ for all $m$. 
\end{corollary}

Of course, all but finitely many of the $\F_m$ will be empty.

We can now prove Proposition \ref{nmprop}.  We apply Corollary \ref{nm-stop-iter}.  The tiles in $\P_{-\infty}$ yield no contribution and can be discarded.
We reduce to establishing that
$$ \| \sum_m F_m \|_{L^2} \lesssim \log(2+M) \| a \|_{l^2}$$
where $F_m := \sum_{\T \in \F_m} \sum_{P \in T} a_P \psi_P$.  If we let $L$ be the first integer larger than $100 \log(2+M)$, it suffices by the triangle inequality to show that
$$ \| \sum_{m: m = l \mod L} F_m \|_{L^2} \lesssim \| a \|_{l^2}$$
for all residue classes $l \mod L$.  Squaring this and using symmetry it suffices to show that
$$ \sum_{m: m = l \mod L} \|F_m\|_{L^2}^2 + \sum_{m, m': m, m' = l \mod L; m' > m} |\langle F_m, F_{m'} \rangle| \lesssim \|a\|_{l^2}^2.$$
Note that if $m,m' = l \mod L$ and $m' > m$ then $m' \geq m+L$.
Introduce the quantities $A_m := 2^{-2m} \sum_{\T \in \F_m} |I_\T|$; from \eqref{msum} it suffices to show that
$$ \sum_m \|F_m\|_{L^2}^2 + \sum_{m,m': m' > m+L} |\langle F_m, F_{m'} \rangle| \lesssim \sum_m A_m.$$
From \eqref{nm-crude} we have $\|F_m\|_{L^2} \lesssim A_m^{1/2}$, and so we reduce to showing that
$$ \sum_{m,m': m' > m+L} |\langle F_m, F_{m'} \rangle| \lesssim \sum_m A_m.$$
We now use Proposition \ref{nm-crude} to obtain
$$ \| F_m\|_{L^4} \lesssim 2^{m} M (\sum_{\T \in \F_m} |I_\T|)^{1/4} = 2^{m/2} M A_m^{1/4}$$
and
$$ \| F_{m'}\|_{L^{4/3}} \lesssim 2^{m'} M (\sum_{\T \in \F_m} |I_\T|)^{3/4} = 2^{-m'/2} M A_{m'}^{3/4}$$
and hence by H\"older's inequality
$$ |\langle F_m, F_{m'} \rangle| \lesssim 2^{-(m'-m)/2} M^2 A_m^{1/4} A_{m'}^{3/4} \lesssim 2^{-(m'-m)/2} M^2 (A_m + A_{m'}).$$
Summing this and using the geometric series formula we conclude
$$ \sum_{m,m': m' > m+L} |\langle F_m, F_{m'} \rangle| \lesssim 2^{-L/2} M^2 \sum_m A_m$$
and the claim follows from the definition of $L$.  This concludes the proof of Proposition \ref{nmprop}, and \eqref{nonmax} follows.
The proof of Theorem \ref{thm:ct1} (and hence Corollary \ref{cor:ct2}) is now complete.
\endprf

\section{Appendix: a correspondence principle}\label{correspond-sec}

The purpose of this appendix is to prove the following correspondence principle.

\begin{proposition}\label{correspondence}
 Let $A$ be an $(n-1) \times m$ matrix with integer entries.  In addition to the operators $T_{A,\R}^*$ and $T_{A,\X}^*$ defined in
\eqref{eq:avmax} and \eqref{tax-ergodic}, we introduce the operator $T_{A,\Z}^*$ defined on functions $\phi_i: \Z \to \R$ of compact support
\begin{equation*}
T_{A,\Z}^{*}(\phi_1,\ldots,\phi_{n-1})(l):=\sup_{N>0}\frac{1}{(2N+1)^{m}}\sum_{|n_1|,\ldots,|n_m|\le N }\prod_{i=1}^{n-1}|\phi_i(l+\sum_{j=1}^{m}a_{i,j}n_j)|.
\end{equation*}
Let $1 < p_1,\ldots,p_{n-1} \leq \infty$ and $p'_n$ be such that $1/p_1 + \ldots + 1/p_{n-1} = 1/p'_n$.
Then the following claims are equivalent.
\begin{itemize}
\item[(i)] $T_{A,\R}^*$ maps $L^{p_1}(\R) \times \ldots \times L^{p_{n-1}}(\R)$ to $L^{p'_n}(\R)$.
\item[(ii)] $T_{A,\Z}^*$ maps $l^{p_1}(\Z) \times \ldots \times l^{p_{n-1}}(\Z)$ to $l^{p'_n}(\Z)$.
\item[(iii)] For every dynamical system $\X$, 
$T_{A,\X}^*$ maps $L^{p_1}(\X) \times \ldots \times L^{p_{n-1}}(\X)$ to $L^{p'_n}(\X)$, with a bound uniform in $\X$.
\end{itemize}
\end{proposition}

\begin{proof} We first show that (i) implies (ii).  Let $\phi_1,\ldots,\phi_{n-1}: \Z \to \R$ have finite support.  For each such $\phi_i$ define $f_i:\R\to\R$ in such a way that $f_i(x)=\phi_i(l)$ if  $x\in [l-\frac13,l+\frac13]$ for some $l\in \Z$, and 0 otherwise. Note that for each $x\in [l-\frac16,l+\frac16]$ and $N\ge 1$
\begin{align*}
\frac{1}{(2N+1)^{m}}&\sum_{|n_1|,\ldots,|n_m|\le N }\prod_{i=1}^{n-1}|\phi_i(l+\sum_{j=1}^{m}a_{i,j}n_j)|\lesssim \\&\lesssim \frac{1}{(2N+1)^{m}}\int_{|t_1|,\ldots,|t_m|\le N+1}\prod_{i=1}^{n-1}|f_i(x+\sum_{j=1}^{m}a_{i,j}t_j)|d\vec{t}\\&\lesssim T_{A}^{*}(f_1,\ldots,f_{n-1})(x).
\end{align*}
From the hypothesis (i) we thus conclude (ii).

Now we show that (ii) implies (i).  Without loss of generality we may take $f_1,\ldots,f_{n-1}$ to be smooth, positive  and compactly supported.
Approximating an integral by the Riemann sum, we obtain
$$
\| T_{A,\R}(f_1,\ldots,f_{n-1}) \|_{L^{p'_n}(\R)}=$$ $$= \lim_{\eps \to 0}
\eps^{-1/p'_n}
\| \sup_{N > 0} \frac{1}{(2N+1)^m} \sum_{|n_1|,\ldots,|n_{m-1}| \leq N} \prod_{i=1}^{n-1} f_i(\eps (l+\sum_{j=1}^{m}a_{i,j}n_j) ) \|_{l^{p'_n}(\Z)}.$$
Applying the hypothesis (ii) we obtain
$$
\| T_{A,\R}(f_1,\ldots,f_{n-1}) \|_{L^{p'_n}(\R)}\lesssim$$ 
$$ \lesssim \limsup_{\eps \to 0}
\eps^{-1/p'_n} \prod_{i=1}^{n-1} \| f_i(\eps \cdot) \|_{l^{p_i}(\Z)}.$$
Approximating integrals by Riemann sums again and using the scaling hypothesis $1/p_1 + \ldots + 1/p_{n-1} = 1/p'_n$ we obtain (i) as desired.

Now we show that (ii) implies (iii).
Define $M:=\max \{\sum_{j=1}^{m}|a_{i,j}|:1\le i\le n-1\}$. Let $f_i\in L^{p_i}({\bf X})$,  let $L\ge 1$ be an arbitrary  number and let $x\in X$  also be arbitrary. By applying the hypothesis (ii) to the functions $\phi_i$ defined by $\phi_i(l)=f_i(S^lx)$ if $|l|\le (M+1)L$ and $\phi_i(l)=0$ otherwise, we get that 
$$\sum_{|l|\le L}\left(T_{A,{\bf X},L}^{*}(f_1,\ldots,f_{n-1})(S^lx)\right)^{p_n'}\lesssim \prod_{i=1}^{n-1}\left(\sum_{|l|\le L}|f_i|^{p_i}(S^lx)\right)^{\frac{p_n'}{p_i}}, 
$$
with an implicit constant independent on $x$ and $L$. The quantity $T_{A,{\bf X},L}^{*}(f_1,\ldots,f_{n-1})(x)$ denotes the maximal operator over averages with $N\le L$. Integration with respect to $x$ and H\"older's inequality imply that
$$\|T_{A,{\bf X},L}^{*}(f_1,\ldots,f_{n-1})\|_{L^{p_n'}({\bf X})}\lesssim \prod_{i=1}^{n-1}\|f_i\|_{L^{p_i}({\bf X})}.$$ By letting $L\to\infty$
we obtain (iii).

To show that (iii) implies (ii), we specialize (iii) to the finitary dynamical system $X = \Z/N\Z$ with the standard shift $Sx := x+1$ and
the uniform probability measure.  Letting $N \to \infty$ (taking advantage of the uniformity of the bounds in (iii) in $N$) and renormalizing
the probability measure to be counting measure (taking advantage of the scaling condition) we obtain (ii); we omit the details.
\end{proof}


\begin{thebibliography}{99}
\bibitem{As:1} I. Assani, {\em Pointwise convergence of ergodic averages along
cubes}, preprint.
\bibitem{As:2} I. Assani, {\em  Multiple recurrence and almost sure convergence for weakly mixing dynamical systems}  Israel J. Math.  \textbf{103}  (1998), 111--124.
\bibitem{BL} J. Barrionuevo and M. Lacey {\em A weak-type orthogonality principle}  Proc. Amer. Math. Soc.  131  (2003),  no. 6, 1763--1769. 
\bibitem{Bo1} J. Bourgain, {\em Pointwise ergodic theorems for arithmetic sets},  Publ. Math. IHES \textbf{69} (1989), 5-45.
\bibitem{Bo2} J. Bourgain, {\em Double recurrence and almost sure convergence}, J. Reine Angew. Math. \textbf{404} (1990), 140-161.
\bibitem{CH}  M. Christ,{\em On certain elementary trilinear operators},  Math. Res. Lett.  \textbf{8}  (2001),  no. 1-2, 43-56.
\bibitem{De} C. Demeter,  {\em Divergence of  combinatorial averages}, preprint.\bibitem{De1} C. Demeter, {\em Pointwise convergence of the ergodic bilinear Hilbert transform}, accepted for publication to the Illinois Journal of Mathematics.
\bibitem{DTT0} C. Demeter, T. Tao and C. Thiele, {\em A trilinear maximal function via arithmetic combinatorics}, work in progress.
\bibitem{DLTT} C. Demeter, M. Lacey, T. Tao and C. Thiele, {\em Breaking the duality in the Return Times Theorem}, preprint.
\bibitem{FS} C. Fefferman and E. M. Stein, {\em Some Maximal inequalities},  Amer. J. Math.  {\bf 93} (1971), 107--115
\bibitem{Fu} H. Furstenberg, {\em Ergodic behavior of diagonal measures and
a theorem of Szemerdi on arithmetic progressions}, J. Analyze Math. {\bf 31}
(1977), 204--256.
\bibitem{GT} B. J. Green and T. Tao, {\em The primes contain arbitrarily long arithmetic progressions}, preprint.
\bibitem{KH} B. Host and B. Kra, {\em Nonconventional ergodic averages and
nilmanifolds},  Ann. of Math. (2)  161  (2005),  no. 1, 397--488. 
\bibitem{La} M. Lacey, {\em The bilinear maximal functions map into $L\sp
p$ for $2/3<p\leq1$}, Ann. of Math. (2) \textbf{151} (2000), no. 1, 35--57.
\bibitem{Le:1} E. Lesigne, {\em Sur la convergence ponctuelle de certaines moyennes ergodiques},  C. R. Acad. Sci. Paris S\'er. I Math.  \textbf{298}  (1984),  no. 17, 425--428.
\bibitem{Le:2} J-M. Derrien and E. Lesigne, {\em Un th\'eor\'eme ergodique polynomial ponctuel pour les endomorphismes exacts et les $K$-syst\'emes},  Ann. Inst. H. Poincar\'e Probab. Statist.  \textbf{32}  (1996),  no. 6, 765--778.
\bibitem{MTT1} C. Muscalu, T. Tao and C. Thiele, {\em Multilinear
operators given by singular multipliers}, J. Amer. Math. Soc. \textbf{15} (2002),no. 2, 469--496.
\bibitem{MTT3} C. Muscalu, T. Tao and C. Thiele, {\em $L^p$ estimates for
the "Biest" I.  The Walsh case},  Math. Ann.  329  (2004),  no. 3, 401--426.
\bibitem{MTT2} C. Muscalu, T. Tao and C. Thiele, {\em $L^p$ estimates for
the "Biest" II.  The Fourier model},  Math. Ann.  329  (2004),  no. 3, 427--461.
\bibitem{stein:maximal}
E.~M. Stein, \emph{On limits of sequences of operators}, Ann. of Math.
\textbf{74} (1961): 140-170.
\bibitem{S} E. Szemer\'edi, {\em  On sets of integers containing no $k$ elements in arithmetic progression},   Acta Arith. {\bf 27}  (1975), 199--245.
\bibitem{Z} T. Ziegler, {\em Universal characteristic factors and Furstenberg averages}, preprint.
\end{thebibliography}
\end{document}